\documentclass{amsart} \usepackage{amssymb,enumerate,amscd,color}
\theoremstyle{plain}
\newtheorem{theo}[equation]{Theorem}
\newtheorem{lem}[equation]{Lemma}
\newtheorem{cor}[equation]{Corollary}
\newtheorem{prop}[equation]{Proposition}
\newtheorem{Con}[equation]{Conjecture}
\theoremstyle{definition}
\newtheorem{Def}[equation]{Definition}
\newtheorem{Rem}[equation]{Remark} 
\newtheorem{Not}[equation]{Notation}

\newcommand{\bL}{{\mathbb L}}

\newcommand{\F}{{\mathbb F}}
\newcommand{\bH}{{\mathbb H}}
\newcommand{\Q}{{\mathbb Q}}
\newcommand{\R}{{\mathbb R}}
\newcommand{\T}{{\mathbb T}}
\newcommand{\Z}{{\mathbb Z}}
\newcommand{\Kappa}{\boldsymbol \kappa}
\newcommand{\Mu}{\boldsymbol \mu}
\newcommand{\Ml}{{\boldsymbol \mu}_{\gl}}
\newcommand{\ga}{{\mathfrak a}}

\newcommand{\gb}{{\mathfrak b}}
\newcommand{\gc}{{\mathfrak c}}
\newcommand{\gC}{{\mathfrak C}}

\newcommand{\gf}{{\mathfrak f}}
\newcommand{\gG}{{\mathfrak G}}
\newcommand{\gH}{{\mathfrak H}}
\newcommand{\gI}{{\mathfrak I}}
\newcommand{\gl}{{\mathfrak l}}

\newcommand{\gm}{{\mathfrak m}} 
\newcommand{\gn}{{\mathfrak n}}
\newcommand{\go}{{\mathfrak o}}
\newcommand{\gol}{\mathfrak{o}_\mathfrak{l}}

\newcommand{\gp}{{\mathfrak p}}

\newcommand{\gq}{{\mathfrak q}}

\newcommand{\gS}{{\mathfrak S}}
\newcommand{\gW}{{\mathfrak W}}
\newcommand{\cA}{{\mathcal A}}
\newcommand{\be}{\bar{e}}
\newcommand{\cB}{{\mathcal B}}
\newcommand{\cC}{{\mathcal C}}
\newcommand{\cD}{{\mathcal D}}
\newcommand{\cE}{{\mathcal E}}
\newcommand{\cF}{{\mathcal F}}

\newcommand{\cH}{{\mathcal H}}
\newcommand{\cI}{{\mathcal I}}

\newcommand{\cK}{{\mathcal K}}
\newcommand{\cL}{{\mathcal L}}
\newcommand{\cM}{{\mathcal M}}

\newcommand{\cO}{{\mathcal O}}
\newcommand{\cP}{{\mathcal P}}

\newcommand{\cS}{{\mathcal S}}
\newcommand{\cT}{{\mathcal T}}
\newcommand{\cU}{{\mathcal U}}
\newcommand{\cV}{{\mathcal V}}
\newcommand{\cW}{{\mathcal W}}
\newcommand{\cX}{{\mathcal X}}
\newcommand{\cY}{{\mathcal Y}}
\newcommand{\cZ}{{\mathcal Z}}
\newcommand{\Zl}{{\mathcal Z}_{\gl}}
\newcommand{\Fl}{{\mathbb F}_{\gl}}
\newcommand{\uA}{{\underline{A}}}
\newcommand{\uD}{{\underline{D}}}
\newcommand{\dE}{{\Dot{E}}}

\newcommand{\tc}{\tilde{c}}
\newcommand{\bt}{{\bf t}}
\newcommand{\w}{{\bf w}}
\newcommand{\wA}{\widehat{A}}

\newcommand{\m}{{\bf m}}
\newcommand{\M}{{\bf M}}
\newcommand{\ann}{\operatorname{ann}}
\newcommand{\ch}{\operatorname{char}}
\newcommand{\Image}{\operatorname{Image} \, }
\newcommand{\gr}{\operatorname{{\bf gr}} }
\newcommand{\Gal}{\operatorname{Gal}}

\newcommand{\GL}{\operatorname{GL}}
\newcommand{\SL}{\operatorname{SL}}
\newcommand{\GS}{\operatorname{GSp}}
\newcommand{\SP}{\operatorname{Sp}}
\newcommand{\res}{\operatorname{res}}
\newcommand{\ord}{\operatorname{ord}}
\newcommand{\rad}{\operatorname{rad}}
\newcommand{\rank}{\operatorname{rank}}

\newcommand{\len}{\operatorname{length}}
\newcommand{\rk}{\operatorname{rank}}

\newcommand{\End}{\operatorname{End}}
\newcommand{\Hom}{\operatorname{Hom}}
\newcommand{\Ext}{\operatorname{Ext}}
\newcommand{\Leg}[2]{\left(\begin{smallmatrix} \underline{#1} \\ #2 \end{smallmatrix}\right)}   
 
\makeatletter
\renewcommand*\l@subsection{\@tocline{2}{0pt}{30pt}{0pt}{}}
\makeatother
\begin{document}

\title{Paramodular abelian varieties of odd conductor}
\subjclass[2010]{Primary 11G10; Secondary 14K15, 11F46}
\keywords{Abelian variety, finite flat group scheme, polarization, division field, paramodular group} 
\author[A. Brumer]{Armand Brumer}
\address{Department of Mathematics, Fordham University, Bronx, NY 10458}
\email{brumer@fordham.edu}
\author[K. Kramer]{Kenneth Kramer}
\address{Department of Mathematics, Queens College, Flushing, NY 11367; The Graduate Center CUNY, New York, NY 10016}
\email{kkramer@qc.cuny.edu}

\thanks{Research of the second author partially supported by NSF grant DMS 0739346}

\begin{abstract}
A precise and testable modularity conjecture for rational abelian surfaces $A$ with trivial endomorphisms,  $\End_\Q A = \Z$, is presented. It is   consistent with our examples, our non-existence results and recent work of C.\! Poor and D.\! S.\! Yuen on weight 2 Siegel paramodular forms.  We obtain fairly precise information on $\ell$-division fields of {\em semistable} abelian varieties,  mainly when  $A[\ell]$ is reducible, by considering  extension problems for group schemes of small rank. 
\end{abstract}

\maketitle 

\tableofcontents

\section{Introduction}\numberwithin{equation}{section}
The Langlands philosophy suggests that the $L$-series of an abelian surface $A$ over $\Q$ might be that associated to a Siegel cuspidal eigenform of weight 2 with rational eigenvalues, for some unspecified group commensurable with $\SP_4(\Z).$  We recall from \cite{Shi} that the ring of endomorphisms of   abelian surfaces $A$ can be either $\Z$ or an order in a quadratic number field $k$. The latter $A$ are of $\GL_2$-type, as defined by  \cite{Rib4}. It is thus a consequence of the work of Khare and Winterberger \cite{KhW} on the Serre conjecture that they are classically modular and their $L$-series are products of two  $L$-functions attached to newforms on $\Gamma_0(N)$ or $\Gamma_1(N),$ depending on whether $k$ is  real or not. All  examples of modularity known to us (\cite{SMT}, \cite{Oka}, \cite{Sma}) involve surfaces of $\GL_2$-type and depend on the lift from classical forms to Siegel modular forms created by  Yoshida  \cite{Yos} for this purpose. Deep work by Tilouine (\cite{Til1},\cite{Til2})  and   Pilloni (\cite{Pil})  uses Hida families   to obtain overconvergent p-adic modular forms associated to  certain abelian surfaces under strong assumptions. 

Our long term project,  originally a study of genus two curves of prime conductor provoked  by the thesis of Jaap Top \cite{Top}, became a search for  a precise and {\em testable}  modularity conjecture for {\em all} abelian surfaces $A$ defined over $\Q$ and  {\em not} of $\GL_2$-type, that is those for which $\End_\Q(A)=\Z$.  We believe we have found one.

The appropriate modular forms are on the paramodular group (\cite{Ibu}, \cite{Gri1}) of level $N$, namely $K(N)= \gamma \hbox{M}_4(\Z)\gamma^{-1}\cap \SP_4(\Q),$ with $\gamma=\text{diag}[1,1,N,1].$  Explicitly: 
\begin{eqnarray*}
K(N)&=&\left\{ g\in \SP_4(\Q)\ \Bigg\vert\ g=
\left[\begin{matrix}*&*&*/N&*\cr N*&*&*&*\cr
N*&N*&*&N*\cr N*&*&*&*\end{matrix}\right]
 \right\}, 
\end{eqnarray*}
where $*$ is an integer. The quotient of the Siegel upper half space $\mathfrak{H}_2$ by   $K(N)$ is the coarse  moduli space of abelian surfaces with $(1,N)$-polarization \cite{BiLa}.   In order to study its Kodaira dimension, Gritsenko (\cite{Gri1}, \cite{Gri2})  introduced a Hecke-equivariant lift from classical Jacobi forms  $J_{k,N}$ to paramodular forms of weight $k$ on $K(N),$  as a variant of the Saito-Kurokawa lift. These lifts violate the Ramanujan bounds and their  $L$-series have poles, and so appear to be of no interest. In fact, they play a crucial role in the construction of the desired paramodular cuspforms \cite{PoYu1}.  In \cite{RS2}, {\em newforms} on $K(N)$ are defined  as  Hecke eigenforms perpendicular to the images of operators from paramodular forms on lower levels. We shall refer to those cuspidal newforms $f$ of weight 2 perpendicular to the Gritsenko lifts as {\em non-lifts} on $K(N)$. Write $L(f,s)$ for the associated degree 4 $L$-function, often denoted $L(f,s,\text{spin})$. 

Let  $\displaystyle{\T_\ell(A)=\lim_{\leftarrow} A[\ell^n]}$ be the Tate module of $A$.
Motivated by results of \cite{And}, \cite{RS1}, \cite{Tay1} and by the compatibility with standard conjectures \cite{Ser2} on  the Hasse-Weil L-series $L(A,s)$,  we  propose the following hypothesis.   

\begin{Con} \label{wcon} There is a one-to-one correspondence between isogeny classes of  abelian  surfaces $A_{/\Q}$  of conductor $N$ with $\End_{\Q} A=\Z $ and weight 2 {\em non-lifts} $f$ on $K(N)$ with rational eigenvalues, up to scalar  multiplication. Moreover, the L-series  of $A$ and $f$ should agree and  the $\ell$-adic representation of  $\T_\ell(A)\otimes\Q_\ell$ should be  isomorphic to those associated to $f$   for any $\ell$ prime to $N.$ 
\end{Con}

\noindent See section \ref{Mod}, added in April 2018, for a modification of Conjecture \ref{wcon}.

\vspace{5 pt}

In contrast to Shimura's classical construction from elliptic newforms,  no  known method  yields an abelian surface from a Siegel eigenform. 

\medskip 

It is  difficult to determine the number of non-lift newforms $f$ on $K(N)$ and more than a few Euler factors of $L(f,s).$ Counting isogeny classes of surfaces of given conductor is even less accessible. However,  it is  easy to compute as many Euler factors  as desired for an explicitely  known abelian surface.  

\medskip

Throughout this paper, $\go$ denotes the ring of integers of a totally real number field of degree $d$ and  $\gl$ a prime of $\go$ above $\ell$ with residue field $\Fl.$  The reader might, on first reading, keep in mind the most important case, namely $d=1$ and so $\go=\Z.$  Our excuses for handling the more general situation are threefold. Our arguments  applied without much extra effort, except slightly more notation. As explained below, they provide  more evidence toward our conjectures. Finally, it seems an entertaining challenge to produce examples, for instance, of abelian fourfolds $A$ with $\End_\Q A$ an order in a quadratic number field. 

\begin{Def} 
The abelian variety $A_{/\Q}$  is of {\em $\go$-type} if $\End_\Q A \, \simeq \, \go.$ Its conductor has the shape $N_A=N^d$, cf. Lemma \ref{sqrfree}, and the reduced conductor is $N_A^0 := N$.  When  $\dim A=2d,$ we say  $A$ is $(\go,N)$-{\em paramodular}.     
\end{Def}

\begin{Rem}\label{rem} \hfill \begin{enumerate}[i)]
\item An $\go$-type abelian variety with $\dim A=d$ has  real multiplication and so, by \cite{KhW}, is a quotient of $J_0(N)$, the Jacobian of the modular curve $X_0(N),$ where $N$ is the {\em reduced} conductor.
\item  An $(\go,N)$-{\em paramodular} abelian variety  is $\Q$-simple, is {\em not} of $\GL_2$-type and its Rosati involution acts trivially on $\go.$
\item A surface $A_{/\Q}$  is paramodular exactly when $\End_\Q A=\Z$, but $\End_{ \overline{\Q}} A$ might be larger. If $K$ is a quadratic field and $E_{\! /\! K}$ is  an elliptic curve,  not $K$-isogenous to its conjugate, then the Weil restriction $A=R_{K/\Q} E$ is paramodular.
\end{enumerate}
\end{Rem}

Guided by the case of abelian varieties over $\Q$ with real multiplications, a more  optimistic conjecture generalizing our earlier one is the following.
\begin{Con}\label{scon}
Let   $f$ be a weight 2 non-lift  for $K(N)$. Let $\go$ be the maximal order in the totally real number field $k_f$ generated by the Hecke eigenvalues of $f.$    Then  there is an  $(\go,N)$-paramodular abelian variety $A_f$ with  $L(A_f,s)=\prod_{\sigma} L(f^{\sigma},s),$ where $\sigma$ runs through the embeddings of $k_f$ into $\R.$ Conversely, an abelian variety $A$ of $(\go,N)$-paramodular type  should be isogenous  to $A_f$ for a weight 2 non-lift newform $f$ on $K(N)$. 
\end{Con}

Current technology might verify our conjecture for Weil restrictions of elliptic curves (cf.\! \cite{Tay2}) and surfaces with $\End_{\overline{\Q}} A \supsetneq \Z.$  In fact, \cite{JLR}\footnote{We thank Brooks Roberts for sending us this preprint upon receipt of our manuscript.} implies that if an  elliptic curve over a real quadratic field is ``Hilbert modular," then its  Weil restriction is paramodular of the predicted level (see Appendix B).  It is conceivable that our precise paramodular  conjecture could be proved, assuming that the $L$-series of a paramodular surface is the $L$-series of a cuspidal automorphic representation of $\GS_4(\mathbb{A}_\Q).$  

To support the conjectures on the arithmetic side, we must prove that  no $(\go,N)$-paramodular variety $A$   exists when no paramodular non-lift exists and produce a member  of each isogeny class for each non-lift that does exist.   For comparison, few results on non-existence or counts of elliptic curves of a given conductor $N$ were known before modularity was proved, even though the issue reduces to S-integral points on the discriminant elliptic curves $c_4^3-c_3^2=1728\Delta,$ for $\Delta $ involving only primes of $N.$ For abelian surfaces $A$, there is no analogous diophantine equation and the problem is exacerbated by the plethora of group schemes available as constituents of $A[\ell]$, as  illustrated by Appendix \ref{DATA}.  The  profusion of intricate lemmas reflects the existence of varieties  satisfying conditions close to the ones we impose.  We mention some of  the subtleties encountered below.

\begin{enumerate}[{\rm i)}] 
\item  When we show that there is no abelian surface of conductor $N$, our proof actually shows that there  is no semistable $(\go,N)$-paramodular abelian variety $A$ with $|\F_\gl|=2$.  In some cases, Conjecture \ref{scon}  and likely paramodular forms suggest  that such $(\go,N)$-paramodular $A$'s of dimension at least  four do exist, explaining why $N$ was not eliminated. 
\item  There are  non-semistable abelian surfaces $B$ such that $B[\ell]\simeq A[\ell]$ as Galois modules, for a putative  semistable surface $A$. 
\item The conductor  of $A[\ell]$ can be a proper divisor of the conductor of the abelian surface $A.$ This makes it difficult to rule out divisors and multiples of conductors of existing abelian surfaces. 
\end{enumerate}

This paper mainly treats non-existence.  The abelian variety $A_{/\Q}$ is semistable if, for each prime $p$, the connected component of the special fiber of its N\'eron model $\cA$ fits into  
\begin{equation}\label{TB}
0 \to \cT_p \to \mathcal {A}^0_p \to \mathcal{B}_p \to 0
\end{equation} 
with $\cB_p$ an abelian variety and $\cT_p$ a torus.  
When $A$ is semistable, deep  results of Grothendieck  \cite{Gro} and Fontaine \cite{Fo}, imply that the  number fields $\Q(A[\ell^n])$  have such tightly controlled ramification that their non-existence rules out certain conductors.  For this, reason we  specialize to semistable  abelian varieties.

\begin{Rem}\label{rem2} \hfill \begin{enumerate}[i)]
\item  Lemma \ref{sqrfree} shows that an $(\go,N)$-paramodular variety  is semistable if $N$ is {\em squarefree}. This is not necessary: the conductor of a semistable paramodular surface  with totally toroidal reduction at $p$ is divisible by $p^2.$ 
\item  All endomorphisms of $A$ are defined over $\Q$ when $A$ is semistable \cite{Rib1}. In particular, a semistable paramodular variety is absolutely simple.
\end{enumerate}
\end{Rem}

No algorithm to find all abelian surfaces of given conductor  is  known,  even less those  not $\Q$-isogenous to a Jacobian, so  we looked for surfaces  by whatever method we could. We include non-principally polarized surfaces and  Jacobians $J(C)$ of conductor $N$ such that $C$ has bad reduction consisting of two genus one curves meeting in one point  at some primes $p\nmid N$ (see Appendix B).    As a special case of Theorem \ref{PPAV}, based on \cite{KhW} and \cite{How}, a paramodular abelian surface of {\em prime} conductor is $\Q$-isogenous to a Jacobian.  

\medskip

As a concrete numerical application of our general results:

\begin{prop} Suppose $A$ is  a semistable abelian  surface of  {\em odd} non-square  conductor $N.$ 
\begin{enumerate}[{\rm i)}] \item If $N\le 500$, then $N$ can only be {\rm 249,\,277,\,295,\,349,\,353,\,389,\,427,\,461} for which examples are known or  {\rm 415,\,417},  which should {\em not}  occur. 
\item See Tables {\rm\ref{BadCond}} and {\rm\ref{EX}} for the data obtained for odd conductors $N < 1000.$
\end{enumerate} 
\end{prop}

The work of Poor and Yuen provides support for our conjecture.  Tables of  cusp forms of weight 2 on $K(N)$ for {\em primes} $N \le 600$ are in \cite{PoYu1}.  For {\em all}  conductors $N \le 1000$  and a few other values, further evidence  will be in \cite{PoYu2}. 

We compare our results with theirs, including some still unpublished data.  There are at least as many known or suspected paramodular non-lift newforms of weight two with rational eigenvalues as known isogeny classes of paramodular surfaces,  including those not semistable or of even conductor.  For almost all non-lifts $f$, we found an  abelian surface $A$ of the same conductor whose  Euler factors  agree with those of $f$  at {\em very small} primes.   Also, the parity of the rank of $A(\Q)$ matches that predicted by the $\epsilon$-factor of $f$.  When we showed that no abelian surface of  a given conductor $N<1000$ exists, their data suggest that all weight two paramodular newforms  with rational eigenvalues are Gritsenko lifts.  
 
Suppose  $A$ has a polarization of degree prime to $q$ and a torsion point of  order $q.$ Then there is   a filtration on $A[q]$ with  a  subgroup and, by  duality,  a quotient of order $q.$ Thus, the characteristic polynomial of Frobenius   at a prime $\ell$ of good  reduction is congruent to  $H_{\ell}(x)= (1-x)(1-{\ell}x)(1-a_{\ell} x+{\ell}x^2) \bmod{q}.$  By Serre's conjecture \cite{Ser3,KhW}, there is an eigenform $g$ of weight 2 on  $\Gamma_0(N)$ with Euler polynomial at $\ell$  congruent to $(1-a_\ell x+\ell x^2) \bmod{\mathfrak{q}}$ for some $\mathfrak{q}\,|\, q.$ Since $H_\ell$ is the Euler factor of the Gritsenko lift $G(g)$ of a Jacobi form attached to $g$, this suggests the possibility of a  congruence  mod $\mathfrak{q}$ between the Fourier series of the non-lift $f$ associated to $A$ by our conjecture and $G(g)$.  Such  matching congruences were  found in \cite{PoYu1} and lend   further supporting evidence.

While the data is  far from complete, it seems convincing enough for publication and dissemination of the conjecture, at least as a challenge.

Recall that the Langlands dual group of GSp$_{2g}$ is GSpin$(2g+1)$.   When interpreted on the split orthogonal group SO$(g+1,g)$ and its associated homogenous space,  the groups $\Gamma_0(N)$ for $g=1$ and $K(N)$ for $g=2$, which at first sight seem so different, are   both instances of  similar  subgroups. Let $L$ be an integral lattice with inner product. There is a natural map $\tau_L\!: \, {\rm SO}(L)\to {\rm O}(\hat{L}/L)$, where $\hat{L}$ is the dual lattice of $L$. The stable orthogonal group $\widetilde{\rm O}(L)$    is the intersection  of  the kernel of $\tau_L$ and the spinor map \cite{Nik,Gri3}. Let $\bH$ be the hyperbolic plane and consider the lattice $\bL_g(N)=\bH^g\perp \langle 2N \rangle$, where $\langle 2N \rangle$ is spanned by a vector of length $2N$. By  \cite[Prop.\! 1.2 and p.\! 485]{Gri3}, $K(N)/\langle \pm 1\rangle$ corresponds to $\widetilde{\rm O}(\bL_2)$ under the identification of GSp$_4$ with SO$(3,2).$ Similarly, $\Gamma_0(N)/\langle \pm 1\rangle$ corresponds to  $\widetilde{\rm O}(\bL_1)$ in the identification of PSL$_2$ with SO$(2,1)$.  Upon  learning this at the June 2010 Conference in his honor, B.\! Gross  immediately generalized our conjecture to one for symplectic  motives in a letter to Serre \cite{Gro}. 

\medskip

\noindent {\bf Acknowledgments.} We are grateful for the opportunity to lecture on preliminary aspects of this work at Edinburgh, Essen, MSRI,  Tokyo, Kyoto, Osaka,  Irvine, Rome, Banff, Shanghai, Beijing and New York.  We were inspired by Ren\'e Schoof, who kindly provided us with preprints.  His hospitality and support to the first author during a visit to Roma III in May 2005  helped  this project along. The contributions  of  Brooks Roberts and Ralf Schmidt as well as those of Cris Poor and David S. Yuen were decisive to our main Conjecture. We thank them heartily for that as well as for useful conversations and correspondence. We also wish to thank the referee for many useful suggestions to improve  the exposition. 

\section{Overview of the paper} \label{Over}
To avoid excessive repetition, we adhere to some conventions for the whole paper. For any finite set  $S$ of primes, let $\Z_{S}=\Z[\{p^{-1} \, |\, p\in S\}]$ and  $\ell$ always be a prime not in $S.$ Write $p_v$ for the prime in $\Z$ below the valuation $v$. The constant group scheme of order $\ell$ over $\Z_S$ is denoted $\cZ_{\ell}=\Z/{\ell}\Z$ and its Cartier dual is $\Mu_{\ell}$. We use $\go$ for the Dedekind ring of integers in a totally real number field of degree $d$ and $\gl$ for a prime ideal of $\go$ over $\ell$. An $\go$-module scheme \cite[p.\! 148]{Tat2} is an abelian group scheme $\cW$ with a homomorphism from $\go$ to $\End \cW$.   The associated Galois module is the $\go$-module of points $W = \cW(\overline{\Q})$.  We have the one-dimensional $\F_\gl$-module schemes $\Mu_\gl = \Mu_\ell \otimes_{\F_{\ell}} \Fl$ and $\cZ_\gl = (\Z/\ell\Z) \otimes_{\F_{\ell}}  \Fl,$ defined in \cite[p. 46]{Maz}.  
We reserve $\cZ$ (resp. $\cM$) for an \'etale (resp. multiplicative) $\gl$-primary $\go$-module scheme over $\Z_S$ all of whose simple constituents are isomorphic to $\Zl$ (resp. $\Mu_\gl$).  We shall shorten this to ``filtered by $\cZ_\gl$'s" or ``filtered by $\Mu_\gl$'s."

We  often abbreviate ``abelian variety" to ``variety" since, aside from curves, they are  the only varieties we consider. We henceforth assume that {\bf all abelian varieties  are semistable of $\go$-type and isogenies are $\go$-linear}, unless the contrary is explicitly stated.  We denote by $N_A$ the conductor of the abelian variety $A$.   

For the cases originally studied in \cite{Fo, Sch1, BK1}, the discriminants of the fields encountered were small enough to ensure that the only simple group schemes occurring were   $\Mu_\ell$ and  $\Z/\ell\Z,$ with their extensions being split. To prove our non-existence results we needed to consider other simple groups schemes and   non-split extensions.  

Accordingly, we generalize an important category introduced by Schoof in \cite{Sch1}. Fix a set $S$ of primes and a prime $\gl$ in $\go$ above $\ell$ not in $S$.  Let $\underline{A}$ be the category of finite flat $\gl$-primary module schemes $\cW$ over $\Z_S$.  Let $\uD$ be the full subcategory of those such that $(\sigma-1)^2=0$ on the associated Galois module $W$ for all $\sigma$ in the inertia groups of the places over $S.$ Clearly, $\cZ_{\gl}$ and  $\Mu_\gl$ belong to $\uD.$

Let $A$ be a semistable abelian variety of $\go$-type with good reduction outside $S$. As in  \cite{Sch1}, Grothendieck's semistable reduction theorem \cite{Gro} implies  that $A[\gl^n]$ and its subquotients belong to $\uD.$ The  {\em exceptional} $\F_\gl$-module schemes are the simple constituents of $A[\gl]$ {\em not} isomorphic to $\cZ_{\gl}$ or $\Mu_\gl,$ if any. The associated Galois modules, also called {\em exceptional}, are thus the irreducibles whose $\F_\gl$-dimension is at least two. Let $\gS_\gl^{all}(A)$ be the multiset of simple $\F_\gl[G_{\Q}]$-modules in a composition series for $A[\gl]$ and $\gS_\gl(A)$  the multiset of exceptionals,  each with its multiplicity. By Proposition \ref{JH1}, $\gS_\gl^{all}(A)$ and $\gS_\gl(A)$ are isogeny invariants.  

An $\gl$-primary $\go$-module scheme in $\uD$ is {\em prosaic} if all its simple constituents are one-dimensional $\F_\gl$-module schemes.   To account for the obstruction to switching adjacent simple constituents in a composition series, the concept of a {\em nugget} is developed in \S \ref{IntNug}.  A {\em prosaic nugget} is an $\go$-module scheme $\cW$ such that $0\subsetneq\cZ\subsetneq\cW$, with $\cZ$ filtered by $\cZ_\gl$'s and  $\cW/\cZ=\cM$   filtered by $\Mu_\gl$'s, and no increasing filtration of $\cW$  has a $\Mu_\gl$ occurring before a $\cZ_\gl.$   See \S \ref{IntNug} for the  more delicate notion and properties of a nugget  with an exceptional subquotient. 

Theorem \ref{trivial} constrains the number of one-dimensional constituents of $A[\gl]$. Put $\Omega(n)=\sum_p \ord_p(n)$ and $\Omega_{\ell}(n)=\sum_{S_{\ell}} \ord_p(n)$, where $S_{\ell}=\{\text{primes }p \equiv \pm 1 \, \bmod{\tilde{\ell}}\}$ with $\tilde{\ell} = 8$ if $\ell=2$, $\tilde{\ell}= 9$ if $\ell=3$ and $\tilde{\ell}= \ell$  otherwise.   Then Corollary \ref{nilpgenbd} gives
$$  
2 \, \dim A \le \Omega(N_A) + \Omega_\ell(N_A) 
$$ 
when $\Q(A[\gl])$ is an $\ell$-extension of $\Q(\Mu_\ell).$ There are  hyperelliptic Jacobians of small dimension for which $A[\ell]$ is prosaic
and  the upper bound is attained for $\ell=2$. 

\begin{Not}
Let $\gI_A$ be the category of abelian varieties $\Q$-isogenous to $A$, with isogenies as morphisms.  If $A$ is of $\go$-type, $\gI^\gl_A$ is the subcategory of abelian varieties of $\go$-type whose morphisms are $\go$-isogenies with $\gl$-primary kernels.  
\end{Not}

In \S \ref{PQR}, we introduce the concept of a {\em mirage}.  A mirage $\gC$  associates to each $B$ in  $\gI_A^\gl$ a set $\gC(B)$ of certain  $\F_\gl$-module subschemes of $B[\gl]$, with natural maps induced by isogenies.  As an example, $\gC(B)$ might be the set of $\F_\gl$-subschemes of $B[\gl]$ filtered by $\Mu_\gl$'s.  Other choices depend on Grothendieck's filtration of the Tate module at semistable primes of bad reduction.   We say that $B$ is  {\em obstructed}  (with respect to $\gC$) if $\gC(B)=\{0\}$ and that $\gC$ is {\em unobstructed} if no $B$ is obstructed.  

Proposition \ref{Faltings} shows  that if  $\gC$ is  unobstructed, then there is a $B$ isogenous to $A$ and  a filtration 
$
0 \subset \cW_1 \subset \dots \subset \cW_s = B[\gl^r],
$
with $\cW_{i+1}/\cW_i$ in $\gC(B/\cW_i)$ for all $i$.  We  choose our mirages so that such a filtration  cannot exist and then exploit the special properties that the  obstructed members of $\gI_A^\gl$ satisfy.  As one illustration, when $\Q(A[\gl])$ is a 2-extension for some $\gl\,|\,2$ and all primes dividing $N_A$ are $3 \bmod{4},$ we prove in  Theorem \ref{pqrthm} that $2 \dim A \le \Omega(N_A).$

Some of the criteria in \S 4--6 depend on arithmetic invariants of  extensions of exceptional module schemes $\cE$ in $\gS_\gl(A)$. These invariants depend on the arithmetic of the number fields generated by the points of such extensions. Those number fields are typically large Galois extensions which, in our applications, are the Galois closures of small cyclic extensions of tractable number fields with well controlled conductors.  In \S 7, we estimate the latter when $\dim_{\F_\gl}E = 2$, to the point that Magma can be invoked.  Similarly, when $\dim_{\F_2}A[\gl]= 4,$ enough  information on $\Q(E)$ is obtained that we could rely on the  Bordeaux tables \cite{BT}.  Finally, our data on paramodular varieties are summarized in the Appendices.

\section{Preliminaries} \label{Prelims} \numberwithin{equation}{subsection}
\subsection{Basics} 
Let $\F$ be a finite field of characteristic $\ell,$ $G$ a finite group and $V$ a finite $\F[G]$-module. The contragredient $\widehat{V}=\Hom_{\F}(V,\F)$ is an $\F[G]$-module via the action  on $V$ and the trace ${\rm Tr}_{\F/\F_\ell}$ induces an isomorphism $\widehat{V} \simeq \Hom_{\F_\ell}(V,\F_\ell)$. 
Now, let $G$ be a quotient of $G_{\Q}=\Gal(\overline{\Q}/\Q)$ by an open subgroup.  Write $\F(1) = \F \otimes \omega$ for the Tate twist by the mod-$\ell$ cyclotomic character $\omega$ and let $V^*=\Hom_\F(V,\F(1)).$     A non-degenerate additive pairing 
$ 
[ \,\; , \,\; ]: \; V \times V \to \F_\ell(1)
$  
satisfying $[g(x),g(y)] = \omega(g)[x,y]$ for all $g$ in $G$ and $[\alpha x, y] = [x, \alpha y]$ for all $\alpha$ in $\F$ is equivalent to an $\F[G]$-isomorphism  $V^* \simeq V$. We say that $V$ is a symplectic Galois module if, in addition, the pairing is alternating.   Then $\dim_\F V = 2n$ is even and, upon the choice of a symplectic basis, $V$ yields a Galois representation into 
$$
{\rm R}_{2n}(\F) := \{ g \in \GS_{2n}(\F) \, | \; [gx,gy] = \omega(g)[x,y] \text{ for all } x,y \in V \}.
$$

\noindent If $W$ is an $\go[G]$-module, let $W^G$ be the submodule  fixed pointwise by $G.$  When  $V$ is simple, $\m_V(W)$ is the multiplicity of $V$ in any composition series for  $W.$  The annihilator of an $\go$-module $M$ will be written as $\ann_\go M$.  

\medskip

We use a capital calligraphic letter for a finite flat group scheme and the corresponding capital Roman letter for its Galois module of $\overline{\Q}$-points, e.g.\! $\cV$ and $V$ respectively.  We write $\Q(V)$ for the field defined by the points of $V=\cV(\overline{\Q}).$ The Cartier dual of $\cV$ is $\cV^D=\Hom(\cV,{\mathbb G}_m)$ and its Galois module is $V^{*}.$   

What we need about abelian schemes and their polarizations over Dedekind domains may be found in the first few pages of \cite{Oda} and \cite{FC}. Under our standing assumption that $A_{/\Q}$ is of $\go$-type, with good reduction outside $S$, the group $A[\mathfrak{a}]$ of $\mathfrak{a}$-division points is an $\go$-module scheme over $\Z_S$ for any ideal $\mathfrak{a}$ of $\go$ prime to $S$.

The following result of Raynaud (\cite{Con1}, \cite{Ray1}) allows us to treat group schemes that occur as subquotients of known group schemes via their associated Galois modules.  In essence, the generic fiber functor induces an isomorphism between the lattice of finite flat closed $R$-subgroup schemes of $\cV$ and finite flat closed $K$-subgroup schemes of $\cV_{|K}$, where $K$ is the field of fractions of $R$.

\begin{lem}\label{CoRay} 
Let $R$ be a Dedekind domain with quotient field $K$ and $\cV$ a finite flat group scheme over $R$ with generic fiber $V=\cV_{|K}.$  If $W=V_2/V_1$ is a subquotient of $V$, for closed immersions of finite flat $K$-group schemes $V_1 \hookrightarrow V_2 \hookrightarrow V$, there are unique closed immersions of finite flat $R$-group schemes $\cV_1 \hookrightarrow \cV_2\hookrightarrow \cV$, such that $V_i=\cV_{i \, |K}$, and there is a unique isomorphism $\cV_2/\cV_1\simeq \cW$ compatible with $(\cV_2/\cV_1)_{|K}\simeq W.$ 
\end{lem}  

If $\cV$ is an $\go$-module scheme, then $\go$-module scheme subquotients of $\cV$ correspond to $\go[G_\Q]$-module subquotients of $V$ by this lemma. 
While we  depend on  this lemma, the reader could instead rely on the  Mayer-Vietoris sequence of \cite[Prop 2.4]{Sch1}.  

\medskip

Consider a strictly increasing filtration of $\go$-module schemes over $\Z_S$, 
\begin{equation} \label{cF}
\cF=\{0=\cW_0 \subset \cW_1\subset \dots \subset\cW_s=\cW\},
\end{equation}
where $\cW$ is $\gl$-primary and the inclusions are closed immersions.  We denote the list of successive quotients by $\gr \cF=[\dots,\cW_i/\cW_{i-1},\dots]$, often writing $\gr \cW$ without explicitly naming the filtration from which it arose. When $\cF$ is a composition series, the multiset of module schemes appearing in $\gr \cF$ may depend on the choice of $\cF$.  By the Jordan-H\"older theorem, the corresponding multiset of irreducible Galois modules does not.

We say that $\cW$ or $\cF$ is {\em prosaic} if all the composition factors are isomorphic to   $\Zl$ or $\Ml$, i.e.\! their associated Galois modules are one-dimensional over $\F_{\gl}$.  

\begin{Not}  Let $\cV$ be an $\gl$-primary $\go$-module scheme over $\Z_S.$  By standard abuse, we write $\cV^{et}$ for the maximal \'etale quotient of $\cV_{|\Z_{\ell}}$ and $\cV^{m}= ((\cV^D)^{et})^D$ for the maximal multiplicative subgroup of $\cV_{|\Z_{\ell}}.$ Similarly, $\cV^0=(\cV_{|\Z_{\ell}})^0$ will denote the connected component and $\cV^{b}=\cV^0/\cV^m$  the biconnected subquotient.  Once a place $\lambda$ over $\ell$ is chosen, with decomposition group $\cD_{\lambda},$ we use the symbols $V^{et}$, $V^m$,  $V^0$ and $V^{b}$ for the corresponding $\cD_{\lambda}$-module.  \end{Not}

We  have the important result of Fontaine, as formulated in  \cite[Thm.\! 1.4]{Maz} and stated here for finite flat $\gl$-primary $\go$-module schemes $\cV_1, \cV_2$ over $\Z_S$.

\begin{lem} \label{Mazur} 
If  $\ell$ is odd and  $V_1 \simeq V_2$ as Galois modules, then $\cV_1\simeq\cV_2.$ 
This holds for $\ell=2$ if, in addition,  $\cV_1^{et} = \cV_2^{et} = 0$ or $\cV_1^m = \cV_2^m = 0.$\end{lem} 

We recall some information about Cartier duality of $\go$-module schemes over $\Z_S.$ 

\begin{lem}\label{selfdual}
Let $\cW\subseteq\cV$ be finite flat $\go$-module schemes over $\Z_S.$  Any
isomorphism $f:\cV^D\simeq\cV$ induces a pairing on the  Galois module $V.$ The submodule scheme of $\cV$ corresponding to $W^{\perp}$
is  $\cW^{\perp}=f((\cV/\cW)^D)$  and $\cW^D$ is isomorphic to $\cV/\cW^{\perp}.$  If $W$ is totally isotropic, $\cW^{\perp}/\cW$ is isomorphic to its Cartier dual.
\end{lem}

\proof
The dual of  the exact sequence $ 0\to \cW \to \cV \to \cV/\cW\to 0$ is
$$
 0\to (\cV/\cW)^D \to \cV^D \to \cW^D\to 0.
$$
Hence the Galois  module corresponding to  $f((\cV/\cW)^D)$ is
$W^{\perp}.$  If $\cW_1\subset \cW_2,$  then
$$
0\to (\cV/\cW_2)^D \to (\cV/\cW_1)^D \to (\cW_2/\cW_1)^D\to 0.
$$
Apply to $\cW_1=\cW$ and $\cW_2=\cW^{\perp}$ to verify the last
claim. \qed

\begin{lem} \label{gensym}
Let $W$ be a symplectic $\go[G_\Q]$-module; $V$ an  irreducible submodule.
\begin{enumerate}[{\rm i)}]
\item Then $V$ is annihilated by some prime $\gl$ of $\go$ and one of  
the following holds:

\noindent {\rm a)} $V$ is nonsingular and $W=V\perp V^{\perp}$,  or

\noindent {\rm b)} $V$ is totally isotropic, $V^{\perp}/V$ is  
nonsingular and  $W/V^{\perp}\simeq V^* = \Hom(V,\Mu_{\ell})$.

\item If $V$ is cyclic as an $\go$-module then $V$ is totally isotropic.
\item If $W$ is semisimple, all irreducible submodules of $W$ are  
nonsingular precisely when $W$ contains no non-zero totally isotropic  
submodule.
\end{enumerate}
\end{lem}

\proof
If $V$ is irreducible, $V \cap V^\perp = 0$ or $V$ and (i) easily  
follows.   If $V$ also is cyclic as an $\go$-module, it is  one-dimensional over $\F_\gl = \go/\gl$.  For any $a,b$ in $\F_\gl$, we  
can solve $c^2 + d^2 = ab$ in $\F_\gl$.  Then the alternating pairing  
on $W$ satisfies $\langle ax,bx \rangle = \langle cx,cx \rangle +  
\langle dx,dx \rangle = 0$, proving  (ii).  Suppose that $W$ is  
semisimple.  If $W$ has no totally isotropic submodule, then every  
irreducible submodule is non-singular by (i).  The converse in (iii)  
is clear. \qed

\begin{lem} \label{polar1}
Let $\cW$ be a  self-dual $\F$-module scheme over $\Z_S$  whose  $\F[G_\Q]$-module $W$ is symplectic and has a unique  simple constituent $E$ such that $\dim E\ge 2$. Then  there is a self-dual subquotient $\cE$ of $\cW$ with  Galois module $E$. 
\end{lem}
\proof
Use induction on the size of $W$. Let $\cX$ be a simple submodule scheme of $\cW$. If $\cX$ is  one-dimensional, then $X$ is isotropic and  by induction applied to $\cV = \cX^\perp/\cX$, we may  suppose there is {\em no}  one-dimensional Galois submodule in $W$. Thus $X\simeq E$.  If $X$ is isotropic, then $\cX^D\simeq \cW/\cX^\perp,$ a contradiction. If $X$ is nonsingular, then $X^\perp$   has no  one-dimensional Galois submodule and so  $W=X$ and we are done. \qed

\medskip

\noindent {\bf Warning:}   The subgroup scheme corresponding to a self-dual Galois submodule $W$ is not necessarily isomorphic to its Cartier dual.

\subsection{Tate module and conductor} \label{TMC} Let $A$ be semistable of $\go$-type and dimension $g$. Fix a prime $\lambda$ in $\overline{\Q}$  over the prime $\ell$ of good reduction. 
Let  $\ell\go=\prod_{\gl} \gl^{e_{\gl}}$ and  $f_{\gl} = [\F_{\gl}\!:\!\F_{\ell}].$ Then  $\go_{\ell}:=\go\otimes \Z_{\ell}=\prod_{\gl} \go_{\gl}.$  Let  $\T_\ell(A)$ be the Tate module  and $\displaystyle{\T_{\gl}(A)=\lim_{\longleftarrow}A[\gl^n]}$.   The  actions of $\go_{\ell}$ and  Galois commute.  
\begin{lem} \label{met}
We have $\rank_{\go_{\gl}} \T_{\gl}(A) = 2g/d.$   For fixed $\lambda$, $\T_\gl(A)^m$ and $\T_\gl(A)^{et}$ are pure free $\gol$-submodules of the same rank, which may vary with $\gl$.
\end{lem}

\proof  We know from \cite{Rib2} that $\T_{\ell}(A)=\prod_{\gl} \T_{\gl}(A)$ is a free $\go_{\ell}$-module of rank $2g/d.$  From the canonical isomorphism to the Tate module of the reduction,  $\T_\gl(A)^{et}$ is a free $\go_{\gl}$-module.  As a free quotient, $\T_\gl(A)^{et}$ is a direct summand, and so  is pure.   By Cartier duality, $\T_\gl(\widehat{A})^m$ is free  of the same rank.  Corresponding to any $\go$-polarization, there is an isogeny $A \to \widehat{A}$ preserving the multiplicative component and so $\T_\gl(A)^m$ and $\T_\gl(\widehat{A})^m$ also have the same rank.  To show that $\T_\gl(A)^m$ is pure, one may use the fact that it is the submodule of $\T_\gl(A)$ orthogonal to $\T_\gl(\widehat{A})^0$. \qed

\medskip 

Since  $\go$ acts by functoriality on the connected component of the special fiber of the N\'eron model of $A,$ the dimensions of $\mathcal{T}_p$  and $\mathcal{B}_p$ in (\ref{TB}) are multiples of $d.$

\begin{Not} \label{torp}
Write $t_p=\dim \mathcal{T}_p$ for the toroidal dimension at $p$ and $\tau_p = t_p/d$.  By semistability,  the reduced conductor of $A$ is $N_A^0 = \prod_p p^{\tau_p}$. 
\end{Not}

We review some results of  Grothendieck (cf.\! \cite{Gro}, \cite{BK1}). 
Since $L_\infty = \Q(A[\ell^\infty])$ depends only on the isogeny class of $A,$ the dual variety $\widehat{A}$ has the same $\ell^\infty$-division field. Let $v$ be a place over $p$ and $\cD_v$ its decomposition group inside $\Gal(L_\infty/\Q).$ The inertia group $\cI = \cI_v \subseteq \cD_v$ acts on $A[\ell^{\infty}]$ and   $A[\gl^{\infty}]$ through its maximal tame quotient, a pro-$\ell$ cyclic group $\langle\sigma_v\rangle,$ whose generator satisfies  $(\sigma_v-1)^2=0.$   The fixed space $M_f(A,v,\ell)=\T_\ell(A)^{\cI}$  is a pure $\go_\ell$-submodule of $\T_\ell(A)$.  The toric space $M_t(A,v,\ell)$  is  the $\go_{\ell}$-submodule of $\T_\ell(A)$ orthogonal to $M_f(\widehat{A},v,\ell)$ under the natural pairing of $\T_\ell(A)$ with $\T_\ell(\widehat{A})$.  
Moreover, $(\sigma_v-1)\T_\ell(A)$ has  finite index in $M_t(A,v,\ell)$. Define $M_f(A,v,\gl)$ and $M_t(A,v,\gl)$ either analogously or  by tensoring  with $\go_{\gl}.$ 
Our earlier remarks together with the   $\go_{\ell}[\cD_v]$-isomorphisms  $M_t(A,v,\ell) \simeq \T_\ell(\cT_p)$ and $M_f(A,v,\ell)/M_t(A,v,\ell) \simeq \T_\ell(\cB_p)$   imply that
 \begin{equation}\label{rank}
\rk_{\go_{\gl}}M_t(A,v,\gl)=\rk_{\go_{\gl}} (\sigma_v-1)\T_{\gl}(A)=\frac{t_p}{d}=\tau_p. \end{equation}
The restriction of $\sigma_v$ to $\Gal(\Q(A[\ell])/\Q)$ generates a subgroup of order 1 or  $\ell.$ 

\begin{Rem}\label{Groth}
The image $\overline{M}_t$  of $M_t(A,v,\gl)$  in $A[\gl]$ is an  $\F_\gl[\cD_v]$-submodule such that $\dim_{\F_\gl} \overline{M}_t = \tau_p$, even if $\sigma_v$ acts trivially on $A[\gl]$.  Hence, $\tau_p$ is bounded from below by the least dimension of any simple      $\F_\gl[\cD_v]$-constituent  of $A[\gl].$  
\end{Rem}

Write $\gf_p(V)$ for the Artin conductor exponent at $p$ of the finite $\gol[G_\Q]$-module $V$ and $N_V$ for its global Artin conductor.  If $\cI$ acts tamely, $\gf_p(V) = {\rm length}_{\gol} V/V^{\cI}.$ Denote by {\rm Frob}$_v$ a choice of arithmetic Frobenius.   

\begin{lem} \label{super1}
Let $0 \to V_1 \to V \stackrel{\pi}{\longrightarrow} V_2 \to 0$ be an exact sequence of finite $\gol[\cD_v]$-modules, with $v$  a prime above $p \ne \ell.$  Suppose  $\cI_v$ acts on $V$ via a { pro}-$\ell$ cyclic group  $\langle \sigma \rangle$ and $(\sigma-1)^2(V) = 0$.  Let $M_i = (\sigma-1)(V_i)$ and $\widetilde{V}_i = V_i^{\langle \sigma \rangle}/M_i$.  Then
\begin{enumerate}[{\rm i)}]
\item there is a well-defined $\gol[\Phi]$-map $\overline{\delta}\!: \, \widetilde{V}_2 \to \widetilde{V}_1(-1)$,  where  $\Phi = {\rm Frob}_v,$ and  \vspace{2 pt}
\item $\gf_p(V) = \gf_p(V_1)+ \gf_p(V_2) + {\rm length}_{\gol}{\rm Im}( \overline{\delta})$.
\end{enumerate}
\end{lem}

\proof 
By the snake lemma, we have the exact sequence of $\Phi$-modules  
\begin{equation} \label{snaker} 
0 \to V_1^{\langle \sigma \rangle} \to V^{\langle \sigma \rangle} \to V_2^{\langle \sigma \rangle} \stackrel{\delta}{\to} V_1/M_1, 
\end{equation} 
where $\delta$ is induced by $y \rightsquigarrow (\sigma-1)(x),$ with $y = \pi(x).$  Since $(\sigma-1)^2(V) = 0,$  we see that $\delta(M_2) \equiv 0 \pmod{M_1}$ and  we obtain the $\go_\gl$-map  $\overline{\delta}.$ Then  (ii) follows. 

To see that $\overline{\delta}$ is a $\Phi$-map, note that $\Phi$ raises to the $p^{\rm th}$ power on $\cI_v$ and  that $\sigma^{p-1}+ \dots +1$ is multiplication by $p$ on $(\sigma-1)(V)$ to obtain $\Phi\overline{\delta}=p \, \overline{\delta}\Phi.$ \qed
 
\begin{lem} \label{super2} 
Let $\cI_v$ act on the  $\F_\gl[\cD_v]$-module $V$ via $\langle \sigma \rangle$ with $(\sigma-1)^2(V) = 0.$  Then $\gf_p(V^*)=\gf_p(V).$  
If $\dim_{\F_{\gl}} V \le 3$ and $V$ is ramified at $v$,  then $\gf_p(V) = 1.$
\end{lem}

\proof
 In the natural pairing $\widehat{V} \times V \to \F_\gl$, we have $\widehat{V}^{\langle \sigma \rangle} = ((\sigma-1)V)^{\perp}$, since $\sigma$ is trivial on $\mu_{\ell}$.  Hence
$
\gf_p(V) = \dim V/V^{\langle \sigma  \rangle} = \dim (\sigma-1)V = \dim \widehat{V}/\widehat{V}^{\langle\sigma \rangle} = \gf_p(\widehat{V}).
$
The last claim follows from $(\sigma-1) V \subseteq V^{\sigma}.$     \qed
\medskip

The  inclusion $M_f(A,v,\gl)/\gl^rM_f(A,v,\gl)\hookrightarrow A[\gl^r]^{\cI}$ implies
 that 
 \begin{eqnarray} \label{lr}
\gf_p(A[\gl^r])&=&\len_{\gol}(A[\gl^r]/A[\gl^r]^{\cI}) \nonumber\\
&\le& r\len_{\gol}A[\gl]-r \rk_{\go_{\gl}} M_f(A,v,\gl) \\
&\le& r \frac{2g}{d} -r\left(\frac{2g}{d} - \frac{t_p}{d}\right)=r \frac{t_p}{d} = r \tau_p . \nonumber
\end{eqnarray} 

\begin{lem}\label{sqrfree} Let $A$  be a  $\Q$-simple abelian variety, not necessarily semistable, with  $\go\subseteq \End_{\Q}A. $ Then the conductor $N_A = N^{d}$ for some integer $N.$
\begin{enumerate}[{\rm i)}]
\item If $N$ is squarefree, then $A$ is semistable and the quotient field of $\go$ is a maximal commutative subfield of $\End_\Q^0\!A =\End_\Q\!A\otimes\Q.$ 
\item If $A$ is semistable and $g=\dim A$ is prime, then either $\End A=\Z$ or $A$ is classically modular. 
\end{enumerate}
\end{lem}

\proof    Since $A$ is $\Q$-simple, $D = \End_\Q^0\!A$ is a division algebra with center a CM field $K$. Let $\dim_{K} D=m^2$ and $[K\!:\!\Q]=r,$ so that a maximal commutative subfield of $D$ has degree $mr$ over $\Q$. The conductor formula,  applied to the $\ell$-adic representation as in \cite{Ser3}, with $\ell$ sufficiently large,   shows that the exponents in the  conductor  must  be multiples of $mr$ and of $d$. Because $\go$ is a maximal order, were $\End_{\Q}A$ to contain  $\go$ properly, the conductor exponent would be a  multiple of $d.$ Similarly,  if $A$ is not semistable at $p,$ the conductor exponent of $A$ at $p$ is at least $2 d$ by \cite[\S 4]{Gro}. This proves (i). 

For (ii), semistability implies $\End^0\!A=\End_\Q^0\!A$ by \cite{Rib1}. Since the invariant differentials form a $D$-module,  $g$ is a multiple of $r m^2,$ so  $D=K$.  If $K$ is not $\Q$, then $[K\!:\!\Q]=g$. When $g$ is odd,  $K$ is totally real and $A$ has RM.  The same holds when $g=2$ because Shimura \cite{Shi} showed that for a surface, $\End A$ cannot be an order in a {\em complex} quadratic number field. Finally, $A$ is a  simple factor of $J_0(N)^{new}$, with $N_A = N^g$, by \cite{KhW}.  \qed

\begin{prop}\label{JH1} If $A$ and $B$ are $\go$-isogenous $\go$-type abelian varieties and $\gl$ is a prime ideal of $\go$,  then $\gS_\gl^{all}(B) = \gS_\gl^{all}(A)$.  
\end{prop}

\proof
By the Jordan-Holder theorem, $\gS_\gl^{all}(A)$ does not depend on the  
choice of composition series for $A[\gl]$.  We use induction on the  
order of the kernel $U$ of the $\go$-isogeny $f: \, A \to B$.  If $U
[\gl]$ is trivial, $f$ induces an isomorphism of $A[\gl]$ to $B[\gl].
$  If not, let $\alpha$ be an element of $\go$ with $\ord_\gl(\alpha)  
= 1$ and consider a composition series
$$
0 \subset V_1 \subset \dots \subset V_r \subset \dots \subset V_n=A
[\ell] \subset V_{n+1} \dots \subset V_{n+r} \subset \dots \subset V_
{2n}=A[\ell^2],
$$
chosen so that $V_r=U[\gl]$ and $\alpha \, V_{n+i}=V_i$  for $i \le n.
$  Visibly, for $C=A/V_r$, we have $\gS_\gl^{all}(A)=\gS_\gl^{all}(C).$    
Moreover, $\gS_\gl^{all}(C) = \gS_\gl^{all}(B)$ by induction hypothesis, since  
the kernel of the induced isogeny $C \to B$ is $U/V_r$. Hence $\gS_\gl^{all}
(B) = \gS_\gl^{all}(A)$.  \qed

\subsection{Ramification}\label{RamSect}
We recall Serre's convention \cite[Ch. IV]{Ser1} for the ramification numbering.  Let $L/K$ be a Galois extension of $\ell$-adic fields with Galois group $G$.  Denote the ring of integers of $L$ by $\cO_L$ and a prime element by $\lambda_L$.  Set
$$
G_n = \{ \sigma \in G \, | \,
         \ord_{\lambda_L}(\sigma(x)-x)  \ge n+1 \; \text{ for all } x \in \cO_L \},
$$
so that $G_0$ is the inertia group and $[G_0\!:\!G_1]$ is the degree of tame ramification.  Recall the Herbrand function: if $m \le u \le m+1$, then
\begin{equation} \label{Herb} 
\varphi_{L/K}(u) =\frac{1}{|G_0|}(\,|G_1|+\dots+|G_m|+(u-m)|G_{m+1}|\,).
\end{equation}

We restate the famous result of Abrashkin \cite{Abr} and Fontaine \cite{Fo} on ramification groups, but using the upper numbering of Serre, namely $G^{m} = G_n$, with $m = \varphi_{L/K}(n)$.  Fontaine's numbering is larger by 1.

\begin{lem}\label{FontaineBound} 
Let $\cV$ be a finite flat group scheme of exponent $\ell$ over $\Z_{\ell}$,  $L = \Q_\ell(V)$ and $G=\Gal(L/\Q_{\ell})$.  If $\alpha > 1/(\ell-1)$,  then  $G^\alpha$ acts trivially on $V$.   Moreover, the root discriminant $r_L$ of $L/\Q_\ell$ satisfies
$$
r_L:=|d_{L/\Q_\ell}|^{\frac{1}{[L:\Q_\ell]}} < \ell^{1+\frac{1}{\ell-1}}.
$$
\end{lem}

We now return to the global situation, with $\cV$ an $\go$-module scheme over $\Z_S$ and $\lambda$ a place of $\overline{\Q}$ over $\ell \not\in S$ .  The set $T_V$ of {\em bad primes} of $V$ consists exactly  of those
dividing   $N_V,$ namely the  finite primes $p \ne \ell$ that ramify in $\Q(V).$ 

\begin{Def}\label{acc} Fix a finite set $S$ of primes and $\gl$ a prime of $\go$ over $\ell$ not in $S.$
An $\gl$-primary $\go$-module scheme over $\Z_S$ is {\em acceptable} if it is a subquotient of $A[\gl^n]$ for some semiabelian scheme $A$ over $\Z$ with good reduction outside $S$ and an action  $\iota:\go\hookrightarrow \End_\Q A$.  
\end{Def}

\begin{Def}\label{GF} 
An $\F_{\gl}[G_{\Q}]$-module $V$ is {\em semistable} if $L=\Q(V)$ satisfies
\begin{enumerate}[{\rm i)}]
\item the inertia group $\cI_\lambda(L/\Q)^\alpha = 1$ for each $\lambda$ over $\ell$ and all $\alpha > 1/(\ell-1)$, and
\vspace{3pt}
\item $\cI_v(L/\Q) =\langle \sigma_v\rangle$, with $(\sigma_v-1)^2(V) = 0$ for each place $v$ dividing $N_V.$
\end{enumerate}
\end{Def}

\noindent The Galois module of an acceptable $\F_\gl$-module scheme  is semistable. 

\begin{Rem} \hfill

i) Acceptable $\gl$-primary $\go$-module schemes form a full subcategory $\uA$
of  the category $\uD$ of Schoof, as generalized in \S \ref{Over} above. The definition of the Baer sum shows that $\Ext^1_\uA(\cU_1,\cU_2)$ is a subgroup of $\Ext^1_\uD(\cU_1,\cU_2)$ when $\cU_i$ are acceptable. 

ii)  Our notion makes  Lemma \ref{CoRay} available.  It seems difficult to check when an object in $\uD$ is acceptable. 
In our applications, the extensions considered are subquotients of $A[\gl^n]$ for  a fixed abelian variety and so are themselves acceptable. 
\end{Rem}

\begin{Rem}\label{MAZ} Let $\cV$ be an acceptable $\F_{\gl}$-module scheme with $\dim V=1.$  The ramification degree of primes over $S$ in $\Q(V)/\Q$ divides $\ell$, but $\Gal(\Q(V)/\Q)$ is a subgroup of $\F_{\gl}^\times$, so  abelian of order prime to $\ell$. Thus $\Q(V)\subseteq\Q(\Mu_\ell)$    and  we conclude that $\cV$ is isomorphic to $\Zl$ or its Cartier dual $\Ml$ by \cite[Prop.\! 1.5]{Maz}.
\end{Rem}

\begin{Not}
Write  $\rad(m) $ for the product of the distinct prime factors of $m.$   
\end{Not}

\begin{prop} \label{simple}
Let $\cW$ be an acceptable $\gl$-primary $\go$-module scheme with $\gl \, \vert \, 2$.  If $\rad(N_W)$ divides one of the integers in
$
\mathfrak{T}_0=\{13,15,17,21,39,41,65\},
$
then $\cW$ is prosaic.  {\rm (GRH is assumed for $\rad(N_W) \ge 39$.)} 
\end{prop}

\begin{proof} Let $F$ be the field generated by $i$ and the square roots of the primes dividing $N_W.$ Any simple constituent $V$ of $W$ is a semistable $\F_\gl[G_\Q]$-module.  Use the refined Odlyzko bounds in \cite[Table 2]{BK3} to conclude that $F(V)/F$ is an abelian extension whose conductor is bounded by \cite[Lemma 5.8]{BK3}. Magma ray class group computations show that $\Gal(\Q(V)/\Q)$ is a 2-group and thus  $\dim_{\F_\gl} V = 1$.  By Remark \ref{MAZ}, $\cV \simeq \cZ_\gl$ or $\Mu_\gl$  so  $\cW$ is filtered by $\Zl$'s and $\Ml$'s, as claimed. 
\end{proof}

\begin{Def}\label{exc}    An irreducible semistable  $\F_{\gl}[G_{\Q}]$-module   $E$   is  {\em exceptional} if  $\dim_{\F_{\gl}} E \ge 2.$  The  $\F_\gl$-module scheme $\cE$ is {\em exceptional} if its  generic fiber is.  When considering a specific exceptional $\cE$, we write $F=\Q(E) ,$ $\Delta = \Gal(F/\Q)$  and $T_E$ for the exact set of primes dividing the conductor $N_E$ of $E.$
\end{Def}

By convention, $\cZ$ and $\cM$ are acceptable $\go$-module schemes over $\Z_S$, filtered by $\cZ_\gl$'s and $\Mu_\gl$'s respectively, with $\ell$ not in $S$.   Clearly then $\Q(Z)/\Q$ and $\Q(M)/\Q(\Mu_\ell)$ are $\ell$-extensions unramified outside $S$. 
  
\begin{lem} \label{ram}
Let $0 \to \cZ \to \cV \to \cX \to 0$ be an exact sequence of $\go$-module schemes over $\Z_S$, $T$ the set of bad primes of $X$ and $L = \Q(V)$.  Then:
\begin{enumerate}[{\rm i)}]
  \item $L/\Q(X)$ is unramified at places $\lambda$ over $\ell$;
  \item $\lambda$ splits in $L/\Q(Z,X)$ if $\cX$ is connected over $\Z_\ell$; 
  \item $L/\Q(X)$ is unramified outside $(S-T) \cup \{\infty\}$ if $\gl \, \cV = 0$; 
  \item $L/\Q(\Mu_\ell)$ is an $\ell$-extension if $\cX=\cM$;
 \item $\Q(Z)=\Q$ and $\cZ$ is constant if $N_Z=1$ and in particular if $N_V=N_X$.
\item If  $0 \to \cV \to \cW \to \cM \to 0$ is exact and $N_V=N_W$ 
        then  $\cM^D$ is constant.
\end{enumerate}
\end{lem}

\proof
Any $\sigma$ in $\cI_\lambda(L/\Q(X))$ acts trivially on $X$ and trivially on $V^{et} = V/V^0$, so $(\sigma-1)(V) \subseteq Z$ and $(\sigma-1)(V) \subseteq V^0$.  Since $\cZ$ is \'etale at $\lambda$ and $\cV^0$ is connected, we have $(\sigma-1)(V) \subseteq Z \cap V^0 = \{0\}$.   This proves (i).  

In (ii), the exact sequence defining $\cV$ splits over $\Z_\ell$, since  $\cV^0 \simeq \cX$ and so the primes over $\ell$ split in $L/\Q(X,Z)$.  In (iii), the ramification degree of each $p$ in $S$ divides $\ell$.  Hence $\Q(X)$ already accounts for all the ramification over each $p$ in $T$.  In (iv), $\Gal(L/\Q(\Mu_\ell))$ is an extension of  $\ell$-groups and therefore is an $\ell$-group. 

In (v),    $\cZ$ is \'etale locally at $\ell$  and $N_Z=1,$ so  $\Q(Z)$ is unramified everywhere.  Thus, $\Q(Z)=\Q$ and $\cZ$ prolongs to a constant $\F_{\gl}$-module scheme over $\Z$, as in \cite[Prop.\! 1.5, Prop.\! 3.1]{Maz}.  If $N_V=N_X$, then $N_Z=1$ by Lemma  \ref{super1}. 

By Cartier duality, (vi) holds.  \qed

\begin{lem} \label{genG} 
Let $V$ be a semistable $\F_\gl[G_{\Q}]$-module, with $\gl \, \vert \, 2$, $L = \Q(V)$ and $G=\Gal(L/\Q)$.  For {\em each} bad prime $p$ of $V ,$ pick {\rm one} place $v$ and a generator $\sigma_v$ of $\cI_v(L/\Q).$  Let $U$ consist of these involutions  and $\sigma_\infty$  a complex conjugation.   In general, $G$ is generated by the conjugates of $U$ and simply by $U$ when $G$ is a $2$-group.
\end{lem}

\proof By semistability, $\sigma_v$ is an involution. Since the fixed field of the  conjugates of $U$ is $\Q$, they generate $G.$  If $G$ is a 2-group and $U$ does not generate $G$, then $U$ lies in a subgroup of index 2 whose fixed field is $\Q(\sqrt{2})$,  violating Lemma  \ref{FontaineBound}.  \qed

\begin{Rem} \label{slam}
The Artin symbol $(-1,\Q_2^{ab}/\Q_2)$ is trivial on $\Q_2^{nr}$ and inverts 2-power roots of unity.  Let $W$ be a Galois submodule of $\T_\gl(A)$ or of $A[\gl^r]$, with $\gl | 2$.  Fix a place $\lambda$ over 2 in $L = \Q(W)$ and let $L_0$ be the maximal abelian subfield of the completion $L_\lambda.$ Let $\sigma_\lambda$ in $\cD_\lambda(L/\Q)$ extend  $(-1,L_0/\Q_2)$ to $L_\lambda$.   
Then $\sigma_\lambda$ acts by inversion on $W^m$ and trivially on $W^{et}$.  If  $W^b = 0$, then  $W^{et} \simeq W/W^m$ and $\sigma_\lambda^2 = 1$ in $\cD_\lambda(L/\Q).$   In the previous lemma, $\sigma_\infty$ may be replaced by $\sigma_\lambda$.  The next lemma shows how  $\sigma_\lambda$ detects ramification in $W[\gl]$. 
\end{Rem}

\begin{lem} \label{slamlem}
Let $\tilde{L}_0$ be the maximal abelian subfield of a 2-extension $\tilde{L}/\Q_2$ with $\Gal(\tilde{L}/\Q_2)^\alpha = 1$ for all $\alpha >1$.   The Artin symbol $a = (-1, \tilde{L}_0/\Q_2)$ is trivial if and only if $\tilde{L}$ is unramified.  
\end{lem}

\proof
By restriction, $\Gal(\tilde{L}_0/\Q_2)^\alpha = 1$ for $\alpha >1$, so $U^{(2)}=1+4\Z_2$ is contained in $N_{\tilde{L}_0/\Q_2}(\tilde{L}_0^{\times})$ by \cite[\!XV,\S2,Cor.\! 1]{Ser1}.  We have $a=1$ if and only if $-1\in N_{\tilde{L}_0/\Q_2}(\tilde{L}_0^{\times}).$   If so, all units of $\Q_2$ are  norms and $\tilde{L}_0/\Q_2$ is unramified.  Then $\tilde{L}_0/\Q_2$ is cyclic, so $\tilde{L}= \tilde{L}_0$ by Burnside's theorem. The converse is proved similarly.  \qed

\subsection{Polarizations}\label{PolSect}  We extend here some results on polarizations from \cite{Mil}, \cite
{Wil} and \cite{How}. We say that $(A,\varphi)$ is $\go$-polarized if  
$\End A = \go$ and the polarization on $A$ induces an $\go$-linear  
isogeny $\varphi: A \to \wA$.  Thus $\Kappa = \ker \varphi$ is a  
Cartier self-dual group scheme whose points form an $\go[G_\Q]$-module.  Throughout this section, $\gn$ is an $\go$-ideal annihilating $\Kappa$.  If the positive integer $n$ is contained in $\gn$, we have the Weil pairing $\be_\gn: \, A[\gn] \times \wA[\gn] \to \Mu_n$ with $\be_\gn(\theta a, a') = \be_\gn(a, \theta a')$  
for all $\theta$ in $\go$.  Define $ \be_\gn^\varphi(a,a') = \be_\gn(a,\varphi(a'))$ to obtain an alternating pairing $A[\gn] \times A[\gn] \to  \Mu_n$
which induces a perfect  induced on $A[\gn]/\Kappa.$   We also have a perfect alternating pairing $\be^\varphi: \, \Kappa \times \Kappa \to \Mu_n $ such that
\begin{equation} \label{kpairing}
\be^\varphi(a,a') = \be_\gn(a,\varphi(\alpha'))\quad\text{whenever}
\quad a' = n\alpha',
\end{equation}
independent of the choices (cf.\! \cite[p.\! 135]{Mil}).  The order of $\Kappa$ is the degree  
of the polarization and the square of its Pfaffian.

\begin{Def}
An $\go$-isogeny $f: \, A \to B$ acts on the $\go$-polarization $\phi
$ of $B$  by $f^*\phi = \hat{f} \phi f$.  Write $(A,\varphi) \succ (B,
\phi)$ if $\varphi=f^{*} \phi$ and $f$ is not an isomorphism.  Say   $
(A,\varphi)$ is {\em minimally $\go$-polarized} if it is minimal with  
respect to this ordering.
\end{Def}

The next lemma is essentially a restatement of \cite[Prop.\! 16.8]{Mil}.

\begin{lem} \label{lowerpol}
Suppose that $\Lambda$ is a proper, totally isotropic $\go[G_\Q]$-submodule of $ \Kappa$ and let $f: A \to B=A/\Lambda$ be the  
canonical map.  Then there is an $\go$-polarization $\phi$, such that  
$(A,\varphi) \succ (B,\phi)$.  Moreover, $\ker{\phi} = f(\Lambda^
\perp) \simeq \Lambda^{\perp}/\Lambda$, where $\Lambda^\perp$ is the  
orthogonal complement of $\Lambda$ with respect to $\be^\varphi$, and  
$|\ker \phi| = |\Kappa|/|\Lambda|^2$.
\end{lem}

\begin{prop}\label{WH}
Let $(A,\varphi)$ be minimally $\go$-polarized.  Then $\Kappa = \ker \varphi$ is an orthogonal direct sum of simple $\go[G_{\Q}]$-modules, symplectic for $\be^\varphi$,  on which $G_{\Q}$ acts  non-trivially.  Further,  the annihilator ideal $\ga = \ann_\go (\Kappa)$ is squarefree.
\end{prop}

\proof
By Lemma  \ref{lowerpol}, $\Kappa$ contains no  $\go[G_\Q]$-submodule totally isotropic for the $\be^\varphi$ pairing.  Then Lemma  \ref{gensym} implies that each irreducible submodule $V$ of $\Kappa$ is symplectic, with non-trivial Galois action.  Use $\Kappa = V \perp V^{\perp}$ to continue by induction.  Since $\Kappa$ is semisimple, its $\go$-annihilator is squarefree. \qed

\begin{cor} \label{JH2} \hfill
\begin{enumerate}[{\rm i)}] \item For each $V$ in $\gS_\gl^{all}(A)$, we have $\m_V(A[\gl]) = \m_{V^*}(A[\gl])$.
\item Let $\gS_\gl(A)=\{E\}$ with $E$ remaining  irreducible as a $\cD_\lambda$-module. Then some subquotient  $\cE$ of $A[\gl^\infty]$ is Cartier selfdual and  biconnected.
\end{enumerate}
\end{cor} 
\proof (i)  Thanks to \ref{JH1}, we may assume that $(A,\varphi)$ is minimally $\go$-polarized.  
If $\ga = \ann_\go(\ker \varphi)$ is prime to $\gl,$ then $A[\gl]$ is  
its own Cartier dual.  If not, $\gl$ exactly divides $\ga$ by Proposition  
\ref{WH} and then $\Kappa[\gl]$ and $A[\gl]/\Kappa[\gl]$ are  
isomorphic to their own Cartier duals.  Since each constituent $V$ of a  
symplectic module $W$ satisfies $\m_V(W) = \m_{V^*}(W)$, we have $\m_V(A
[\gl]) = \m_{V^*}(A[\gl])$. 

 (ii) Let $B$ be minimally $\go$-polarized in the isogeny class of $A.$ Then  $B[\gl]$ is Cartier selfdual or we have an exact sequence $0\to\Kappa\to B[\gl]\to \Kappa'\to 0$ with $\Kappa$ and $\Kappa'$ both self-dual.   Lemma  \! \ref{polar1} yields a self-dual subquotient $\cE$ of $A[\gl]$.  The filtration of $\cE_{|\Z_\ell}$ by multiplicative, biconnected and \'etale subquotients proves our claim. \qed

\begin{Rem}\label{havemu}
The field obtained by adjoining all irreducible $\go[G_\Q]$-constituents of $A[\gl]$ is an isogeny invariant. Taking $(A,\varphi)$ to be minimally $\go$-polarized as above,  either $A[\gl]$ or $(\ker \varphi)[\gl]$  is Cartier self-dual and so $\Mu_\ell\subseteq \Q(A[\gl])$.   \end{Rem}

\begin{lem} \label{flip}
Let  $(A,\varphi)$ be minimally $\go$-polarized, $\Kappa = \ker \varphi$ and $\ga=\ann_\go(\Kappa).$  Let $\theta$ be a totally positive element of $\go$ dividing $\ga$ and write $\ga = \theta \gb$.
\begin{enumerate}[{\rm i)}]
\item There is an $\go$-polarization $\phi$ on $B = A/\Kappa[\theta]
$, such that $\ker \phi$ is  isomorphic to $\left(A[\theta]/\Kappa
[\theta]\right) \oplus \Kappa[\gb]$.
\item If $|\Kappa|$ is minimal for $\gI_A^\gl
$ and  $\gl = (\theta)$ is a prime dividing $\ga$, then $2 \dim_{\F_
{\gl}}\!\Kappa[\gl] \le \dim_{\F_{\gl}}\!A[\gl].$ Further,  some composition factor of $A[\gl]/\Kappa[\gl]$ is 
symplectic.
\end{enumerate}
\end{lem}

\proof
By Proposition \ref{WH}, $\ga$ is squarefree, so $\gb$ is prime to $\theta$. Since $\theta$ is totally positive, it is a sum of squares in the fraction field of $\go$. It follows that $\psi=\varphi \theta$ is a polarization on $A$ whose kernel obviously contains $\Kappa$.  Let $\Lambda = \Kappa[\theta]$ and let $\Lambda^\perp\subseteq \ker \psi$ be its orthogonal complement under the $\be^\psi$ pairing (\ref{kpairing}).   Given $a$ in $\Lambda$ and $a'$ in $\ker \psi$, write $a' = \theta\alpha'$.  By the definitions of the pairings, we have
$$
\be^\psi(a,a') = \be_\gn(a,\psi(\alpha')) = \be_\gn(a,\varphi(\theta \alpha')) = \be^\varphi(a,\theta a').
$$
Since the orthogonal complement of $\Lambda$ with respect to the $\be^{\varphi}$ pairing on $\Kappa$ is $\Kappa[\gb]$, we find that $\Lambda^\perp = \{ a' \in \ker \psi \, | \, \theta a' \in \Kappa[\gb]  \, \}$.   But multiplication by $\theta$ is an isomorphism
on $\Kappa[\gb]$.  Hence $\Lambda^\perp = A[\theta] + \Kappa[\gb] \supseteq \Lambda$ and so $\Lambda$ is totally isotropic for the $\be^\psi$ pairing.  By Lemma  \ref{lowerpol}, there is an induced polarization  $\phi$ on $B$, such that
$ \ker \phi \simeq \Lambda^\perp/\Lambda = (A[\theta] +
\Kappa[\gb])/\Kappa[\theta].$  This proves (i).

Now let $(C,\gamma) \preccurlyeq (B,\phi)$ be $\go$-minimal.   Then any irreducible submodule  of $\ker \gamma$ is nonsingular by Proposition \ref{WH} and (ii) follows from minimality of $|\Kappa|$.  \qed

\begin{lem}
Let $(A,\varphi)$ and $(B,\phi)$ be $\go$-polarized, with $B$ in $
\gI_A^\gl.$ If $\Kappa = \ker \varphi$ has minimal order for $\gI_A^{\gl},$ then $\ga = \ann_\go(\Kappa)$ divides $\ga' = \ann_\go(\ker \phi)$ provided
\begin{enumerate}[{\rm i)}]
\item each prime factor $\gl$ of $\ga$ has a totally positive generator and
\item $\m_E(A[\gl])=1$ whenever $E$ in $\gS_\gl(A)$ has non-trivial Galois action and admits a symplectic pairing.
\end{enumerate}
The second condition is fulfilled when the reduced conductor of $A$ is squarefree.
\end{lem}

\proof
We argue as in \cite[p.\! 213ff]{How}.   Let $\psi = f^*\phi=\hat
{f} \phi f$ be the polarization on $A$ induced by the $\go$-isogeny  
$f\!: \, A \to B$.  Write $\Phi = \ker f$, so $\ker \hat{f}$ is  
isomorphic to $\Phi^* = \Hom(\Phi,{\mathbb G}_{{\rm m}})$.  We can find $\alpha$ and $
\beta$ in $\go$ such that $\beta \varphi = \alpha \psi$.  By Proposition  
\ref{WH}, the $\gl$-primary part of $\Kappa$ is semisimple and  
annihilated by $\gl$.  Furthermore, any simple component $E$ is  
symplectic, with non-trivial Galois action.  For $V$ in $\gS_\gl(A)$,  
the multiplicity $\m_V$ is additive on short  
exact sequences and $\m_V(A[\gl])$ is  isogeny invariant.  It follows  
that
\begin{equation}  \label{mults}
\m_V(\Phi) + \m_V(\Phi^*) +\m_V(\ker \phi) =
\ord_{\gl}(\beta/\alpha)\m_V(A[\gl])+\m_V(\Kappa[\gl]).
\end{equation}
Suppose $\gl$ does not divide $\ga'$.  Put $V = E$  in (\ref{mults}) to show that $\ord_\gl(\beta/\alpha)$ is odd.  By Lemma  \ref{flip},  $A[\gl]/\Kappa[\gl]$ has an irreducible symplectic constituent $E'$ with non-trivial Galois action.  By assumption, $\m_{E'}(A[\gl]) = 1$ so $\m_{E'}(\Kappa[\gl]) = 0$ and then  $\m_{E'}(\ker \phi)$ is odd.  Hence $\gl$ divides $\ga'$.  Indeed, $\ga$ divides $\ga'$ because $\ga $ is squarefree.  \qed

\begin{cor} \label{polar}
There is at most one symplectic module in $\mathfrak{S}_{\gl}(A)$ if one of the following holds, with $N$ the reduced conductor of $A$ and $p$ a prime:
\begin{enumerate}[{\rm i)}]
 \item  $N = p$ and $\ell \le 19$, or
 \item $\ell=2,$  $N = m p$ and $\rad(m)$  divides a $Q$
in $\mathfrak{T}_0$ of {\rm Proposition \ref{simple}}.
\end{enumerate}
\end{cor}

\proof
Under (i) or (ii), Lemma  \ref{simple} shows that $p$ must appear in the conductor of any member of $\gS_{\gl}(A)$.  But $\gf_p(A[\gl])\le 1$  by (\ref{lr}). \qed

\begin{theo}\label{PPAV}
Let $A$ be a semistable $(\go,mp)$-paramodular abelian variety, with $\rad(m) \le 10$ and prime $p \ge 11$.  If the strict ideal class group of $\go$ is trivial, then $A$ is $\go$-isogenous over $\Q$ to a principally polarized abelian variety. In particular, $A$  is $\Q$-isogenous  to a Jacobian if it is a  surface.
\end{theo}

\proof
Assume that $A$ is an $\go$-linear polarization $\varphi$  whose kernel $\Kappa\ne 0$ has minimal order  for the isogeny class.  Let $\gl$ be a prime of $\go$ dividing $\ga = \ann_\go(\Kappa)$ and let $\ell$ lie below $\gl$  in $\Z.$ Proposition \ref{WH} implies that both $W = \Kappa[\gl]$ and $W' = A[\gl]/\Kappa[\gl]$ have irreducible constituents of dimension at least 2 over $\F = \go/\gl$.  By Lemma  \ref{met}, $A[\gl]$ is 4-dimensional over $\F$.  Thus $W$ and $W'$ admit odd {\em irreducible} Galois representations into $\GL_2(\F).$   By semistability and Lemma  \ref{super1}, the conductors of $W$ and $W'$ are squarefree. 

We consider three cases.  If $\ell$ is prime to $pm$, the group schemes $\cW$ and $\cW'$ are finite flat over $\Z_\ell$. The conductor exponents at $p$ satisfy
$\gf_p(W) + \gf_p(W') \le \gf_p(A[\gl]) \le 1$, thanks to (\ref{lr}).  Hence the conductor $n_0$ of one of them is prime to $p$, and so $n_0$ divides $\rad(m)$.  By \cite{KhW}, the corresponding Galois module gives rise to a cusp form of weight 2 on $\Gamma_0(n_0)$, contradicting $n_0\le 10.$

If $\ell=p,$ we work over $\Z_p.$ The action of  $\go$ extends to the N\'eron model and $\tau_p(A)=1.$  We infer from \cite[Ch.\! 4]{MB}\footnote{We thank J. de Jong for this reference.} that there is  a finite flat $\F$-module subscheme $\cF$ of $A[\gl]$ of dimension 3.   Since $\cW$ is irreducible it must be contained in $\cF$.  The conductor of $W$ is prime to $p$ by definition.  We conclude as above.

When $\ell$ divides $m$, the Galois module  $W$ or $W'$ whose conductor $n_0$  is prime to $p$ either has weight 2 or has weight $\ell+1,$ but can be twisted to have weight 2 and  level $\ell n_0$. Since $J_1(\rad(m))$ has genus 0, we again have a contradiction. 

When $A$ is a surface, its conductor precludes  being isogenous to a product of two elliptic curves.  Since we may assume  $A$ to be principally polarized, we conclude from \cite {Weil} that it is a Jacobian.  \qed

\subsection{Cohomology} \label{CohomSec}
If $X$ and $Y$ are finite $\F[G]$-modules, then $\Ext^1_{\F[G]}(X,Y)$ is isomorphic to $H^1(G,\Hom_{\F}(X,Y)),$ as in  \cite{CR,Asch}.  The cohomology class corresponding to an exact sequence of $\F[G]$-modules 
\begin{equation} \label{YVX}
0 \to Y \to V \stackrel{\pi}{\to} X \to 0,
\end{equation}
is represented by a 1-cocycle $c,$ with $c_g(x) = (gi-i)(x) = g(i(g^{-1}x))-i(x)$, where $g \in G$ and $i \in \Hom_\F(X,V)$ is any section of $\pi$.

\begin{lem} \label{makesplit}
Let $X' = \pi(V^G)$ and $Y' = Y \cap \ga_G V$, where $\ga_G$ is the augmentation ideal in $\F[G]$.  The following sequences are $\F[G]$-split exact:
\begin{enumerate}[{\rm i)}]
\item $0 \to Y \to \pi^{-1}(X')  \to X' \to 0$,  \hspace{2 pt} with \hspace{2 pt} $\dim_\F X^G/X' \le \dim_\F  H^1(G, Y)$,
\vspace{5 pt}
\item $0 \to Y/Y' \to V/Y' \to X \to 0$, \hspace{2 pt} with \hspace{2 pt} $\dim_\F Y'/ \ga_G Y \le \dim_\F H^1(G,\widehat{X})$.
\end{enumerate}
\end{lem}

\proof
The cohomology sequence
$
0 \to Y^G \to V^G \stackrel{\pi}{\longrightarrow} X^G \stackrel{\delta}{\to} H^1(G,Y)
$
implies our bound on $\dim X^G/X'$.    If $W$ is a complement in $V^G$ for  $Y^G,$  then $\pi^{-1}(X')=Y+V^G=Y\oplus W$ provides  a splitting in  (i).
  
The homology sequence
$
H_1(G,X) \stackrel{\partial}{\to} Y/\ga_G Y \to  V/\ga_G V \to X/\ga_G X \to 0
$
shows that $\dim Y'/\ga_G Y \le \dim H_1(G,X) = \dim H^1(G,\widehat{X})$, by duality.  If $W$ is an  $\F$-subspace  of $V$ containing $\ga_G V,$ then $W$ is a Galois module.  We may choose $W$ so that $V/\ga_GV = ((Y+\ga_GV)/\ga_GV)\oplus (W/\ga_GV).$   Then $V/Y' = (Y/Y') \oplus (W/Y')$ is a  splitting in (ii).   Alternatively, observe that (i) and (ii) are dual statements.  \qed

\medskip
Let $\{G_i\}$ be a collection of subgroups of $G$.  For each $i$, let $\overline{Y}_i$ be a $G_i$-quotient of $Y$ and $X_i$ a $G_i$-submodule of $X$.    Put 
$
\Kappa_1 = \ker\{H^1(G,Y) \to \prod H^1(G_i,\overline{Y}_i)\}
$ 
and
$\Kappa_2 = \ker\{H^1(G,\widehat{X}) \to \prod H^1(G_i,\widehat{X}_i)\}$,
where the maps are induced by restriction.

\begin{cor} \label{splitcor}
Let $Y^{\rm ker}_i = \ker\{Y \to \overline{Y}_i\}$ and $\overline{V}_i = V/Y^{\rm ker}_i$.  
\begin{enumerate}[{\rm i)}]
\item If the sequences $0 \to \overline{Y}_i \to \overline{V}_i \to X \to 0$ are $\F[G_i]$-split exact for all $i$, then $\dim X^G/X' \le \dim \Kappa_1$, where $X' = \pi(V^G)$.
\item If the sequences 
$
0 \to Y \to \pi^{-1}(X_i) \to X_i \to 0
$
are $\F[G_i]$-split exact for all $i$, then $\dim Y'/\ga_G Y \le \dim \Kappa_2$, where $Y' = Y \cap \ga_G Y$.
\end{enumerate}
\end{cor}

\proof
i) By the splitting, $\pi$ induces a surjection $\overline{V}^{G_i} \twoheadrightarrow X^{G_i}$.  It follows from the diagram below that the image of $\partial$ is contained in $\Kappa_1$.  But $X^G/X' \simeq \Image \partial$.
$$
\begin{CD}
       V^G          @>\pi>>      X^G    @>\partial>>  H^1(G,Y) \\
     @VVV                     @VVV                     @VresVV     \\  
\overline{V}_i^{G_i} @>\pi>> X^{G_i} @>0>> H^1(G_i,\overline{Y}_i) \; .
\end{CD}
$$
 ii) Similarly, the diagram below shows that $\delta$ vanishes on $C = \sum \, {\rm cores} \, H_1(G_i,X_i)$.
$$
\begin{CD}
H_1(G_i,\pi^{-1}(X_i)) @>>> H_1(G_i,X_i) @>0>> Y/\ga_{G_i} Y \\
              &    &    @VcoresVV                  @VVV                             \\
              &    &             H_1(G,X) @>\delta>>  Y/\ga_G Y   
\end{CD}
$$
Hence
$
\dim Y'/\ga_G Y = \dim \Image \delta \le \dim H_1(G,X)/ C = \dim \Kappa_2,
$
where the last equality holds by duality. \qed

\begin{Rem} \label{ForD1}
Let $G_0$ be the subgroup of $G$ acting trivially on both $X$ and $Y$ in (\ref{YVX}) and let $\Delta = G/G_0$.   Assume $\ch(\F) = \ell$ and write $\gH =  \Hom_\F(X,Y)$.  Then the image of $[c]$ under the restriction map
$$
H^1(G,\gH) \overset{\res}{\longrightarrow} H^1(G_0,\gH)^\Delta \, = \Hom_{\F_\ell[\Delta]}(G_0, \gH)
$$
is an $\F_\ell[\Delta]$-homomorphism $\tc: \; G_0 \to \gH$, with $\tc_g(x) = (gi-i)(x)$ for all $g$ in $G_0$ and $x$ in $X$.  If $G_0$ acts faithfully on $V$, then $\tc$ is injective.  Indeed, if $\tc_g = 0$, then $g$ acts trivially on both $Y$ and $i(X)$, so $g= 1$ on  $V = Y + i(X)$.  Let
$$
Y'' = \ga_{G_0} V = \ga_{G_0} (Y + i(X)) = \ga_{G_0} i(X) = {\rm span} \{ \tc_g(X) \, | \,  g \text{ in } G_0 \}.
$$
Then $Y'' \subseteq Y$ and $0 \to Y/Y''  \to V/Y'' \to X \to 0$ is $\F[G_0]$-split exact.
\end{Rem}

\section{Nuggets} \label{IntNug} \numberwithin{equation}{subsection}
\subsection{Introducing  nuggets}\label{IntNug1}
If the $\go$-module scheme $\cV$ has a subscheme isomorphic to $\cZ_\gl$ with quotient isomorphic to $\Mu_\gl$ and $\Ext^1_{\uD}(\Mu_\gl,\cZ_\gl)=0,$  then $\cV$ admits a filtration with  graded pieces in the reverse order. More generally, when  considering increasing filtrations of $\go$-module schemes, we refer to ``moving $\cZ_\gl$ to the right" and ``moving $\Mu_\gl$ to the left", when modifying filtrations in that way.  To account for the failure of splitting and the existence of exceptional $\go$-module schemes, we  introduce the weight $\w(\cF)$ of a filtration, using several invariants which we now define. 

Let $\cW$ be an $\gl$-primary $\go$-module scheme and $\cF$ a composition series for $\cW$.  Let $\bt_m(\cF)$ be the multiplicity of $\Ml$ in $\gr \cF$ and $\bt_e(\cF)$ that  of $\Zl$.  Then $\bt_m(\cF)$ and $\bt_e(\cF)$ are determined by the Galois module $W=\cW(\overline{\Q})$ if $\ell>2.$  The sum $\bt_m(\cF)+\bt_e(\cF)$ is the number of trivial $\F_{\gl}[G_{\Q}]$-constituents of $W$ when $\ell=2.$

\begin{Not}\label{epsilon}
Set $\epsilon_0(\cW)=\bt_m(\cF) + \bt_e(\cF)$  for any composition series $\cF$ of $\cW.$ This depends only on the Galois module $W.$
\end{Not}

\noindent For a fixed $\lambda$ over ${\ell}$, we  write $\bt^{\lambda}_m(\cF)$ (resp. $\bt^{\lambda}_e(\cF)$) as the multiplicity of $\Ml$ (resp.  $\Zl$) in any composition series for $\cW_{|\Z_\ell}$.  They are local invariants of $\cW$ by Lemma  \ref{Mazur}.   When $\cF$ is prosaic, $\bt^{\lambda}_m(\cF) = \bt_m(\cF)$ and $\bt^{\lambda}_e(\cF) = \bt_e(\cF)$. 

For a simple $\go$-module scheme $\cY,$ put $\w(\cY)=\bt_e^{\lambda}(\cY)\bt_m^{\lambda}(\cY).$
Let ${\bf x}(\cF)$ be the number of exceptional constituents and set 
\begin{equation*}
\alpha(\cF)= \bt_e^\lambda(\cF) +{\bf x}(\cF)\; \;  \text{ and } \; \; \beta(\cF) =\bt_m^\lambda(\cF) +{\bf x}(\cF).
\end{equation*}

\begin{Def}
Let $\cF'=\{0=\cW_0\subset\dots\subset\cW_{s-1}\}$ be a filtration as in (\ref{cF}) and $\cY =\cW_s/\cW_{s-1}.$ Define inductively 
$$
\w(\cF)=\w(\cF')+\w(\cY)+\begin{cases}  
                           \alpha(\cF') & \text{ if } \cY=\Ml, \\
                            \bt_e(\cF')\beta(\cY) &\text{ if } \cY \text{ is exceptional,} \\
                            0 & \text{ if } \cY=\Zl.
                                          \end{cases}
$$
\end{Def}

\begin{lem} \label{weightformula} 
\hfill
\begin{enumerate}[{\rm i)}]
\item For $\cV \subset \cW$, let $\cF_1$ and  $\cF_2$ be filtrations of $\cV$ and $\cW/\cV,$ respectively, and write $\cF_1\cF_2$ for the induced filtration of $\cW.$   Then
$$
\w(\cF_1\cF_2)=\w(\cF_1)+\w(\cF_2)-\bt_e(\cF_1)\bt_m(\cF_2)+\alpha (\cF_1)\bt_m(\cF_2)+\bt_e(\cF_1)\beta(\cF_2).
$$
\item If $\cF$ is  a  prosaic filtration of $\cW$, then $\w(\cF)\le \bt_e(\cF)\bt_m(\cF)$ with equality only if all the $\Zl$ are on the left and all the $\Ml$ on the right.
\item  If $\cF^D$ is the filtration of $\cW^D$ induced by $\cF$ via Cartier duality, then $\w(\cF)=\w(\cF^D).$
\end{enumerate}
\end{lem}

\proof The claims follow easily by induction on the length of a filtration, except perhaps for the second part of (ii).   There,  it suffices to check that if $\cF=\cF' \cY$, with $\cY$ simple, then equality holds for $\cF$ if and only if it holds for $\cF'$ and $\cY=\Zl$. \qed

\begin{Def} A filtration on the $\gl$-primary $\go$-module scheme $\cW$ is {\em special} if:
\begin{enumerate}[i)]
\item $0 \subset \cZ \subset \cW$ with $\cZ \ne 0$ \'etale and $\cW/\cZ = \cM \ne 0$ multiplicative, or
\item $0 \subseteq \cZ \subset \cV \subseteq \cW$, with $\cV/\cZ = \cE$ exceptional and at least one of $\cZ$ or $\cW/\cV = \cM$ not trivial.  Then we have the short exact sequences
\begin{equation}\label{specialeq} 
0 \to \cZ \to \cV \to \cE \to 0 \quad \text{ and } \quad
                  0 \to \cV \to \cW \to \cM \to 0.
\end{equation}
\end{enumerate}
\end{Def}

\begin{lem} \label{weightlem}
Let $ 0\to \cZ\to\cV\to \cX\to 0$ be exact, with $\cZ$ filtered by $\Zl$'s,  $\cX$ connected or exceptional and $\gl \, \cX = 0$. Set $\cZ'=\cZ[\gl]$,  $\overline{\cV}=\cV/\cZ'$  and $\overline{\cZ}=\cZ/\cZ'.$
\begin{enumerate}[{\rm i)}]
\item The sequence $0 \to \overline{\cZ} \to \overline{\cV}\to \cX \to 0$ is split exact.
\item If $\gl \, \cZ=0$ then $\gl \, \cV=0.$

\item If $\cW$ admits a special filtration with $\cZ$ and $\cM$ killed by $\gl$, then $\gl \, \cW = 0$.
\end{enumerate}\end{lem}

\proof  By the snake lemma, $0 \to \cZ' \to \cV[\gl] \stackrel{f}{\to} \cX \to \cZ/\gl \, \cZ.$  If $\cX$ is exceptional, $f$ is surjective by irreducibility of the Galois module $X$.  If $\cX$ is connected, its image in the \'etale group scheme $\cZ/\gl \, \cZ$ is trivial, also proving surjectivity of $f$.   Thus the  subgroup $\cV[\gl]/\cZ'$ of $\overline{\cV}$ is isomorphic to $\cX$  and provides a splitting in (i).  Surjectivity of $f$ further implies the isomorphism $\overline{\cZ} \simeq \cV/\cV[\gl]$, from which (ii) follows.  

In (iii), we may now assume $\cV/\cZ$ is exceptional and $\gl \, \cV = 0.$  The snake lemma for multiplication by $\gl$ gives  $0\to \cV \to \cW[\gl] \to \cM \to \cV.$ Then the Mayer-Vietoris sequence \cite{Sch1} shows that $\Hom_R(\cM,\cV)=\Hom_R(\cM,\cZ)=0$.  The first equality holds  by considering the Galois modules and the second because over $\Z_{\ell},$ $\cM$ is connected, while $\cZ$ is \'etale.  \qed   

\begin{Rem} \label{EllKills}
If $\cV_1$ and $\cV_3$ are annihilated by  $\gl$ and $V_1$ and $V_3$ have no Galois constituents in  common, then the exactness of $ 0 \to \cV_1 \to \cV_2 \to \cV_3 \to 0$ implies that $\cV_2$ is annihilated by $\gl$ as well.
\end{Rem}

\begin{prop} \label{weightprop} 
Suppose that $\cW$ admits a filtration $\cF=\cF_1\cE\cF_2$ with $\cF_i$ prosaic. Then $\w(\cF)\le \bt_e^{\lambda}(\cW)\bt_m^{\lambda}(\cW)+\epsilon_{0}(\cW)$ with equality if and only if $\cF$ is special.
\end{prop}

\proof We apply  Lemma  \ref{weightformula} and find 
\begin{eqnarray*}
\w(\cF)&= & \w(\cF_1)+\w(\cE)+\w(\cF_2)+\bt_e(\cF_1)\bt_m(\cF_2)+\alpha(\cE)\bt_m(\cF_2)+\bt_e(\cF_1)\beta(\cE)\cr
&\le &\bt_e^{\lambda}(\cW)\bt_m^{\lambda}(\cW)-\bt_e^{\lambda}(\cE) \bt_m(\cF_1)-\bt_m^{\lambda}(\cE) \bt_e(\cF_2)-\bt_e(\cF_2) \bt_m(\cF_1)\cr&&\hspace{100 pt}+\bt_e(\cF_1)+\bt_m(\cF_2)\cr
&\le &\bt_e^{\lambda}(\cW)\bt_m^{\lambda}(\cW)+\epsilon_{0}(\cW).
\end{eqnarray*}
By Lemma  \ref{weightformula}, the first inequality is strict unless $\gr \cF_i = [\cZ_\gl^{a_i} \,\Mu_\gl^{b_i}]$ and then $\epsilon_0(\cW)=a_1+b_1+a_2+b_2$ while $\bt_e(\cF_1)+\bt_m(\cF_2)=a_1+b_2.$   Hence the last inequality is strict, unless $b_1=0=a_2.$ \qed

\begin{Def}\label{nug}
An $\go$-module scheme $\cW$ is a {\it nugget} if either it is exceptional or it satisfies the following two properties.
\begin{enumerate}[i)] 
\item $\cW$ has no $\go$-subscheme isomorphic to $\Ml$ and no quotient isomorphic to $\Zl.$ 
\item $\cW$ has a {\em special} filtration $\cF.$ No other filtration has  strictly lower weight. 
\end{enumerate}
If the nugget $\cW$ has no exceptional subquotient, it is called a {\em prosaic nugget}.  We usually keep the filtration $\cF$ implicit. 
\end{Def}
 
The Cartier dual of a nugget $\cW$ is a nugget and Lemma  \ref{weightlem} shows that $\gl \, \cW = 0$.  Let $\cZ'$ and $\cM'$ be $\go$-module subschemes of $\cZ$ and $\cM$, with both $\overline{\cZ} = \cZ/\cZ'$ and $\cM'$ non-zero if $\cW$ is prosaic.  Write $\cW'$ for the pre-image of $\cM'$ in $\cW$ and set $\overline{\cV} = \cV/\cZ'$ and $\overline{\cW}=\cW'/\cZ'$.  Then $0 \subseteq \overline{\cZ} \subseteq \overline{\cV} \subseteq \overline{\cW}$ is a special filtration, with $\overline{\cW}/\overline{\cV} \simeq \cM'$.  Lemma  \ref{weightlem}(ii) and Proposition \ref{weightprop} imply that the subquotient $\overline{\cW}$ is a nugget, referred to as a  {\em subnugget} of $\cW$ by abuse.
A prosaic nugget $\cW$ has a subquotient nugget $\cW'$ with $\gr \cW' = [\Zl \Ml],$ called the {\em core},  which may depend on the chosen special filtration. By Lemma   \ref{ram},  $\Q(W')$ is an elementary $\ell$-extension of $\Q(\mu_{\ell}),$ split over $\ell$ and unramified outside $N_{W'}$. 

\begin{cor}  \label{AllZnug}
If a nugget $\cW$ contains a prosaic $\go$-module subscheme $\cY$ with $N_Y=1$, then $\cY  \simeq \Zl^r.$   
\end{cor}

\begin{proof} 
Since $\cY$ is prosaic, each constituent of $Y$ is a 1-dimensional $\F_\gl[G_\Q]$-module and the  corresponding subquotient  of $\cY$ is isomorphic to $\cZ_\gl$ or $\Mu_\gl$ by Remark \ref{MAZ}. 
Once we show that $\cY$ is \'etale at $\ell$, the claim holds by Lemma \ref{ram}(v). 
Let $\cF'$ be a special filtration of $\cW$ of minimal weight, as in Definition \ref{nug}(ii). 

If $\cW$ is not prosaic,  we have a filtration $\cF = \cF_1\cE\cF_2$ of $\cW$, with $\cY \subseteq \cF_1$.   Since $\cF'$ is special, $\w(\cF')$ attains the upper bound in Proposition \ref{weightprop}.   Minimality of $\w(\cF')$ implies that $\w(\cF) = \w(\cF')$, so $\cF$ also is special.  But then $\cF_1$, and {\em a fortiori} $\cY$, is filtered by $\cZ_\gl$'s.  
When $\cW$ is prosaic,  replace Proposition \ref{weightprop} by Lemma \ref{weightformula}(ii) in this argument.  
\end{proof}	

\begin{prop}\label{FiltByNug} 
Any $\gl$-primary $\go$-module scheme $\cW$ has a filtration 
$$
0\subseteq \cW_0\subseteq \cW_1\subseteq \dots \subseteq \cW_{r-1}\subseteq \cW_r=\cW,
$$
with $\cW_0$ filtered by $\Ml$'s,  $\cW/\cW_{r-1}$ filtered by $\Zl$'s and  $\cW_{i+1}/\cW_{i}$ a nugget for $i=0,..,r-2.$  Such a filtration will be called a {\em nugget filtration}.
\end{prop}

\proof 
Denote by $\mu(\cY)$ any maximal subscheme of $\cY$ filtered by $\Ml$'s. Dividing by $\mu(\cW)$, we may assume  that  $\cW$ has no $\Ml$ submodule. For such $\cW$, we prove the claim by induction, by producing a nugget $\cV \subseteq \cW$ with $\mu(\cW/\cV)=0.$
\begin{enumerate}[{\rm i)}]
\item Suppose there is a subscheme $\cZ$ of $\cW$ such that $\gr \cZ =[\Zl^r]$, with $r \ge 1$ and $\Mu(\cW/\cZ) \ne 0$.  Choose one with $r$ minimal and let $\cV$ be the pullback of $\mu(\cW/\cZ).$  By minimality of $r$ and Lemma \ref{weightformula}(ii), $\cV$ is a prosaic nugget.
\item Next, suppose for all subschemes $\cZ$ of $\cW$ filtered by $\Zl$'s, we have $\mu(\cW/\cZ)=0$.  If there is a subscheme $\cX$ of $\cW$ having a filtration with grading $[\Zl^s \, \cE]$, where $s \ge 0$ and $\cE$ is exceptional, choose $\cX$ to minimize $\bt_e^{\lambda}(\cX)$.  Thus the pullback $\cV$ of $\mu(\cW/\cX)$ has a special filtration and $\mu(\cW/\cV)=0.$  Minimality of $\bt_e^{\lambda}(\cX)$ shows that $\cV$ has no $\cZ_\gl$ quotient, so $\cV$ is a nugget by  Proposition \ref{weightprop}.
\item When the  only simple factors of $\cW$ are $\Zl$'s, we are  done. \qed
\end{enumerate}

\subsection{Prosaic nuggets} 
We generalize \cite[Cor.\! 4.2]{Sch2}, allowing for $\go$-action and several bad primes.  In this section, $\F = \Fl$ and $N$ is prime to $\ell.$  Recall that $\tilde{\ell}=8, \, 9$ or $\ell$ if $\ell=2, \, 3$ or $\ell \ge 5$ respectively.  Let $\varpi(N)$ denote the number of distinct prime factors $p$ of $N.$  When $\ell=2$ or $3,$ set $\varpi_{\ell}(N)=\varpi(N)$ if all $p$ dividing $N$ satisfy $p\equiv \pm 1\bmod{\tilde{\ell}}$ and $\varpi_{\ell}(N)= \varpi(N)-1$  otherwise.  When $\ell \ge 5$, define
$$ 
\varpi_{\ell}(N) = \# \{ \text{primes } p \text{ dividing } N \; | \; p\equiv \pm 1 \bmod{\ell} \}.
$$
Write $p^* = (-1)^{(p-1)/2} p$ for $p$ odd.

\begin{prop} \label{SchoofObst}   With $R=\Z[1/N],$ we have 
$\displaystyle{\dim_\F \Ext^1_{R}(\Ml,\Zl)=\varpi_{\ell}(N).}$  \end{prop}

\proof 
Proceed as in the proof of \cite[Prop.\! 4.1]{Sch2} and use the Mayer-Vietoris sequence of \cite[Cor.\! 2.4]{Sch1}, to obtain the exact sequence
\begin{equation} \label{ExtSeq}
0 \to \Ext^1_{R}(\Ml,\Zl) \to \Ext^1_{R[1/\ell]}(\Ml,\Zl)\to \Ext^1_{\Q_{\ell}}(\Ml,\Zl),
\end{equation}
in which the last two terms may be studied via extensions of Galois modules.
  
Let $L$ be the maximal elementary abelian $\ell$-extension of $F = \Q(\Mu_\ell)$ such that $L/F$ is unramified outside $N$ and split over $\lambda_F$.  Set $G = \Gal(L/\Q)$, $G_0 = \Gal(L/F)$ and $\Delta = \Gal(F/\Q)$.   

By Lemma  \ref{weightlem}(iii), an extension $\cV$ of $\Mu_\gl$ by $\cZ_\gl$ over $R$ is killed by $\gl$.  If $V$ is the associated $\F[G_\Q]$-module, Lemma  \ref{ram} implies that $\Q(V)$ is contained in $L$.  Conversely, if $V$ is an $\F[G]$-module extending $\Mu_\gl$ by the trivial Galois module $\F$, then $V$ arises, by (\ref{ExtSeq}), from an $R$-group scheme $\cV$ as above.  It thus suffices to determine $\Ext^1_{\F[G]}(\Mu_\gl,\F) \simeq H^1(G,\F(-1))$, cf.\! \S \ref{CohomSec}.

Since $\Delta$ has order prime to $\ell$, inflation-restriction shows that 
$$
H^1(G,\F(-1))=\Hom_{\F_\ell}(G_0,\F(-1))^\Delta = \F \otimes \Hom_{\F_\ell}(G_0,\F_\ell(-1))^\Delta.
$$
Let $X$ be the subgroup of $F^\times$ whose elements satisfy:  (i) $x \in F_\lambda^{\times \, \ell}$ and (ii)  $\ord_\gq(x) \equiv 0 \bmod{\ell}$ for all $\gq$ not dividing $N$. By Kummer theory, we have a perfect $\Delta$-pairing $G_0 \times \overline{X} \to \F(1)$, where $\overline{X} = X/F^{\times \, \ell}$.  It follows that $ \Hom_{\F_\ell}(G_0,\F_\ell(-1))^\Delta$ is isomorphic to the $\omega^2$-component $\overline{X}_{\omega^2}$, where $\omega$ is the mod-$\ell$ cyclotomic character.

Let $\overline{Y} = Y/F^{\times \, \ell}$, where $Y$ is the subgroup of $F^\times$ satisfying only (ii) above.  Write $\cC$ for the ideal class group of $F$ and $\overline{U}$ for the image in $\overline{Y}$ of the group of units.  The natural action of $\Delta$ on the prime ideals $\gp$ of $F$ dividing $p$ induces an action on $J_p = \prod_{\gp | p} \Z/\ell\Z$.  Schoof shows that  $\dim_{\F_\ell} (J_p)_{\omega^2} = 1$ if $p \equiv \pm 1 \bmod{\ell}$ and 0 otherwise.  In particular, if $\ell = 2$ or 3, then $\dim J_p = 1$ for all $p$.   

We have the exact sequence of $\F_\ell[\Delta]$-modules
$$
1 \to \cC[\ell] \stackrel{h }{\to} \overline{Y}/\overline{U} \stackrel{i}{\to} \prod_{p | N} J_p \stackrel{j}{\to} \cC/\cC^\ell,
$$
with $i$ induced by $y \leadsto (\ord_\gp y)$ and $j$ by $(c_\gp) \leadsto \prod_\gp \gp^{c_\gp}$.  As for $h$, if the ideal $\ga^\ell = (y)$ is principal, then $y$ is in $Y$ and $h$ is induced by $\ga \leadsto y$.  

If $\ell \ge 5$, the $\omega^2$-component of $\cC[\ell^\infty]$ vanishes by the reflection principle and Herbrand's theorem \cite[Thms.\! 6.17, 10.9]{Wash}.   Hence $\dim_{\F_\ell} (\overline{Y}/\overline{U})_{\omega^2} = \varpi_\ell(N)$.  If $\ell = 2$ or 3, we have $\dim_{\F_\ell}   (\overline{Y}/\overline{U})_{\omega^2} = \varpi(N)$.

Let $U_\lambda$ be the group of local units in the completion $F_\lambda$ and use bars to denote the respective multiplicative groups modulo $\ell^{\rm th}$ powers.  Embedding to the completion induces a map of $\overline{Y} \to U_\lambda F_\lambda^{\times \ell}/F_\lambda^{\times \ell} \simeq \overline{U}_\lambda$ and we have the exact sequence
$$
1 \to \overline{U} \cap \overline{X} \to \overline{X} \to \overline{Y}/\overline{U} \stackrel{\beta}{\to} \overline{U}_\lambda/\overline{U}.
$$

If $\ell = 2$ or 3, then $\overline{U} \cap \overline{X} = 1$ by direct calculation.  Moreover, $\dim_{\F_\ell} \Image \beta = 0$ if $p \equiv \pm 1 \pmod{\tilde{\ell}}$ for all $p$ dividing $N$ and 1 otherwise.  This implies our claim.

If $\ell \ge 5$, then $\dim_{\F_\ell} \overline{U}_{\omega^2} = 1$ and $(\overline{U}_\lambda/\overline{U})_{\omega^2} = 1$ by \cite[Thms.\!  8.13, 8.25]{Wash}.  It follows that the non-trivial elements of $\overline{U}_{\omega^2}$ are not $\ell^{\rm th}$ power locally and so $\overline{U}_{\omega^2} \cap \overline{X} = 1$.  This concludes the proof.  \qed

\begin{cor} \label{core} 
Suppose $0\to \Zl\to\cV\to \Ml\to 0$ is a non-split extension.  Then either $N_V$ is divisible by some prime $p\equiv \pm 1\bmod{\tilde{\ell}}$ or else $\ell\le 3$ and $N_V$ is divisible by at least two primes.
\end{cor}

\begin{Rem}\label{corem} Let $\cV$ be a prosaic nugget  with $\gr \cV = [\cZ \Mu_\gl]$ and $\gr \cZ = [\cZ_\gl \, \cZ_\gl].$ Let $p,\, q$ be primes, with $p$ dividing the conductor $N'$ of the core  and $(q, N') = 1.$  Then generators of inertia at $v \, | \, p$ and $w \, | \, q$ can be put in the form
\begin{equation} \label{sigmavw}
\sigma_v=\left[\begin{smallmatrix} 1 & a_v & b_v \\ 0 & 1 & 1 \\ 0 & 0 & 1 \end{smallmatrix} \right]             \quad \text{and}\quad 
\sigma_w=\left[\begin{smallmatrix} 1 & a_w & b_w \\ 0 & 1 & 0 \\ 0 & 0 & 1 \end{smallmatrix} \right] .
\end{equation}
We have $a_v=0$ because $(\sigma_v-1)^2=0$.  Hence either $\cZ\simeq \Zl^2$ or $a_w\neq 0$ for some $q$ dividing $N_V$.   Since $\Q(Z)/\Q$ is an elementary $\ell$-extension unramified at $\ell,$ Kronecker-Weber implies that any prime ramified in $\Q(Z)$ is 1 mod $\ell.$
Any Frobenius $\Phi_v$ at $v$ normalizes $\sigma_v$ and so acts trivially on $Z$.
\end{Rem}

\begin{lem}\label{pnil} 
Let $S$ contain exactly one prime $p\equiv \pm 1\bmod{\tilde{\ell}}$ and let $\cV$ be a prosaic nugget over $\Z_S.$  Then $\dim V=2$, $N_V=p$ and $\cV$ prolongs to an $\go$-module scheme over $\Z[1/p]$ under any of the following conditions.
\begin{enumerate}[{\rm i)}]
\item $\ell\ge 5$ and $S-\{p\}$ consists of primes  $q\not\equiv \pm 1\bmod{\tilde{\ell}}.$ 
\item $\ell=3,\,S=\{p,q\}$ with $q \equiv 2,\, 5 \, \bmod{9},$ or  $q \equiv 4,\,7 \, \bmod{9}$  and $p^{\frac{q-1}{3}}\not\equiv 1 \, (q).$
\item $\ell=2,$  $S=\{p,q\}$ with  $q^* \equiv 5 \, (8)$ and the Hilbert symbol $(p^*,q^*)_\pi=-1$ for some place $\pi.$
\end{enumerate}
\end{lem}

\proof  If $\dim V>2,$  dualize or pass to a subnugget if necessary, to get  $\gr \cV = [\cZ \, \Mu_\gl]$, with $\gr \cZ = [\Zl \, \Zl].$  Then any core of $\cV$ has conductor $p$ by Proposition \ref{SchoofObst}. 

We prove next that  $\Q(Z)=\Q.$ For (i) and (ii), this follows from  the Remark above.  For (iii), note that  $\Phi_v$ acts trivially on  $\Q(Z)$ but non-trivially on the cubic subfield of $\Q(\Mu_q).$  For (iv), if not,  then  $\Q(Z)=\Q(\sqrt{q^*})$ and $\Q(V)$ is a $D_4$ field whose existence requires $(p^*,q^*)_\pi=1$ for all places $\pi,$ as explained below.

For $v$ over $p$, $Y = (\sigma_v-1)(V)$ is a  Galois  submodule of $Z$ with corresponding subscheme $\cY\simeq \Zl.$ But then the core $\cV/\cY$ is unramified at $p$. \qed
  
\medskip   

We now  suppose $\gl \, | \, 2.$  When writing $\Q(\sqrt{d})$ and its character $\chi_d,$  we assume $d$ is squarefree.  Recall that for $p$ prime, $\chi_d(p)$ is the Legendre symbol\! $\Leg{d}{p}\!.$   Let $D_4^r(d_1,d_2)$ be the set of $D_4$-extensions $M/\Q$  such that
\begin{enumerate}[i)]
\item $\vert \cI_v(M/\Q) \vert \le 2$ at odd  $v$, $\cI_\lambda(M/\Q)^\alpha=1$ for all $\alpha >1$ at even  $\lambda$ and
\item the subfields fixed by the Klein 4-groups in $\Gal(M/\Q)$ are $\Q(\sqrt{d_1})$ and $\Q(\sqrt{d_2})$, with $d_1, d_2$ odd and coprime.
\end{enumerate}
For $d_1,\, d_2$ as in (ii), such an $M$ exists exactly if  $d_1x^2 + d_2 y^2 = 1$ is solvable in $\Q$, i.e.  for all  $\pi,$ the Hilbert symbols $(d_1,d_2)_\pi = 1,$ cf.\! \cite{RR1,JLY}. Let $k = \Q(\sqrt{d_1d_2})$ and let $d_3$ be the product of the odd  $p$ such that some $v$ over $p$ ramifies in $M/k.$

\begin{Not}\label{D4not} Let $D_4(d_1,d_2)\supseteq D_4^{nr}(d_1,d_2)\supseteq D_4^{sp}(d_1,d_2)$ be the following subsets of $D_4^r(d_1,d_2).$   In the first $d_3=1$, in the second $M/\Q$ is unramified over 2 and in the third $M/k$ splits completely over 2. 
\end{Not} 

\begin{lem} \label{twistM} If $M'$ is in $D_4^r(d_1,d_2),$ then 
some twist $M$ of $M'$ is in $D_4(d_1,d_2).$  If $d_1 \equiv d_2 \equiv 1 \, (4)$, we may even arrange that $M$  be in $D_4^{nr}(d_1,d_2).$
\end{lem}

\proof
If $d_3 \ne 1$, adjust its sign so $d_3 \equiv 1 \, (4)$ and let $L =M'(\sqrt{d_3})$.  Because $K = \Q(\sqrt{d_1}, \sqrt{d_2})$ is the maximal abelian subfield of $M'$, we have $\sqrt{d_3} \not\in M'$.  Hence $G = \Gal(L/\Q) \simeq D_4 \times C_2$ and the central involution of $\Gal(M'/\Q)$  may be extended  to $c\in\Gal(L/\Q(\sqrt{d_3}))$.  If $\tau$ generates $\Gal(M'(\sqrt{d_3})/M')$, then the center of $G$ is $\langle c, \tau \rangle = \Gal(L/K)$.  For each prime $p$ dividing $d_3$, there is a place $v$ over $p$ ramified in $M'/K$ and in $K(\sqrt{d_3})/K$. Hence $\cI_v(L/\Q) = \langle c \tau \rangle$ and the subfield $M$ of $L$ fixed by $c \tau$ satisfies our claim, since $c\tau$ is central.  

Suppose $d_1\equiv d_2 \equiv 1 \, (4),$  $M'$ is in $D_4(d_1,d_2)$ and $\lambda \, \vert \, 2$ ramifies in $M'$.   Set $L = M'(i)$ and observe that $\cD_\lambda(L/\Q)$ is abelian. By Lemma  \ref{slamlem}, $g = (-1,L_\lambda/\Q_2)$ restricts nontrivially  to $\Gal(M'/K)$ and to $\Gal(K(i)/K)$, so $g = c \tau$ and $M = L^{\langle c \tau \rangle}$  is in $D_4^{nr}(d_1,d_2).$   \qed

\begin{prop} \label{D41}
Let $V \supsetneq V_1 \supsetneq V_2 \supsetneq 0 $  be semistable $\F[G_\Q]$-modules with $\dim_\F V=3.$ Set $X = V/V_2$, $K = \Q(V_1,X)$ and $L = \Q(V)$.  \begin{enumerate}[{\rm i)}]\item Then $\gcd(N_{V_1}, N_{X}) = 1$ and no prime dividing $N_{V_1}N_{X}$ ramifies in $L/K$. \item If $\Q(V_1) = \Q(\sqrt{d_1})$ and $\sqrt{d_2}$ is in $\Q(X)$, then  $(d_1,d_2)_\pi = 1$ for all $\pi$. 
\end{enumerate}
\end{prop}

\proof
The shape (\ref{sigmavw}) of the generators of inertia at bad places proves (i).  Then $d_1$ and $d_2$ are coprime and  they are odd by Definition\! \ref{GF}. 

Let $\sigma_i$ be an involution  of $G =\Gal(L/\Q)$, non-trivial on $\Q(\sqrt{d_i})$.  By matrix verification, $\sigma_1$ is trivial on $X$, $\sigma_2$ is trivial on $V_1$ and centralizes the elementary 2-group $H = \Gal(L/\Q(\sqrt{d_1}, \sqrt{d_2})),$ while their commutator $c = [\sigma_1,\sigma_2] \ne 1.$  The centralizer of $\sigma_1$ is trivial on $X$ and so fixes $\sqrt{d_2}$.  Using  a commutator identity, this implies that $c \not\in [\sigma_1,H]=\{\, [\sigma_1,h]\, |\,h \in H\}$.  There is a  maximal subgroup  $J$ of $H$ containing $ [\sigma_1,H]$ but not $c.$  Since $J$ is normal in $H$ and $G = \langle \sigma_1,\sigma_2, H \rangle$, $J$ is normal in $G$.  Now $M = L^J$ is in  $D_4^r(d_1,d_2).$ \qed

\begin{cor} \label{D42}
Let $\cV$ be an $\F$-module scheme  with $\gr \cV = [\Zl \, \Ml \,\Ml]$ and $V$ its Galois module.   Then $2$ is unramified in $\Q(X)$ and splits in both $\Q(V_1)/\Q$ and $\Q(V)/\Q(X)$.  If $\Q(V_1) = \Q(\sqrt{d_1})$ and $\sqrt{d_2} \in \Q(X)$, then $d_1 \equiv 1 \, (8)$ and $d_2 \equiv 1 \, (4).$  
If $N_V = |d_1d_2|$ and neither $d_i= 1,$ then $L$ contains a $D_4^{sp}(d_1,d_2)$ field. 
\end{cor}

\proof 
The grading on $\cV$ implies that 2 is unramified in $L$ and splits in $\Q(V_1)$.  Hence $d_2 \equiv 1 \, (4)$ and $d_1 \equiv 1 \, (8)$.  Moreover, even places split in $L/\Q(X)$ by Lemma  \ref{ram}(ii).  Since $d_1 \equiv 1 \, (8)$, they also split in $L/\Q(\sqrt{d_1d_2})$.  If $N_V=|d_1d_2|$, then the field $M$ defined in the proof  above is in $D_4^{sp}(d_1,d_2).$  \qed

\subsection{Invariants of nuggets} 
First we recall a result from \cite[Chap. VII, \S 1]{HBII}. 

\begin{lem}[]\label{BC}
Let  $\chi$ be the character of an irreducible $\F[\Delta]$-module $E$ and $\F_\chi=\F_{\ell}(\chi(g)\ |\ g\in \Delta).$   Write $\dE$ for $E$, viewed as an $\F_{\ell}[\Delta]$-module.  There is an irreducible $\F_{\ell}[\Delta]$-module $X$ such that $\dE=X^a$, with $a = [\F:\F_\chi]$, and
\begin{enumerate}[{\rm i)}]
\item $X\otimes_{\F_{\ell}}\F=\bigoplus E^{\eta}$ is a direct sum of non-isomorphic  conjugate representations, with $\eta$ running over $\Gal(\F_\chi/\F_\ell)$;
\item $(\End_{\F_{\ell}[\Delta]} X)\otimes_{\F_{\ell}}\F=\End_{\F [\Delta]}(X\otimes_{\F_{\ell}}\F)\simeq (\End_{\F[\Delta]} E)^{b}$  with $b = [\F_\chi:\F_\ell]$.
\end{enumerate}
Viewed as an $\F_\ell[\Delta]$-module,  $\widehat{E} = \Hom_\F(E,\F) \simeq \widehat{X}^a$, where $\widehat{X} = \Hom_{\F_\ell}(X,\F_\ell)$, and similarly, $E^*  = \Hom_{\F_\ell}(X,\Mu_\ell) \simeq X^{* \, a}$. 
\end{lem}   

\begin{Not} \label{seldim} 
Let $E$ be an exceptional $\F[G_\Q]$-module, $T$ its set of bad primes and $X$  an irreducible constituent of $\dE$.  Set $F=\Q(E)=\Q(\widehat{E})$ and $\Delta= \Gal(F/\Q)$.   For $S \supseteq T,$ let  $\Lambda_{E}(S)$ be the maximal elementary $\ell$-extension $\Lambda$ of $F$ such that 
\begin{enumerate}[{\rm i)}]
\item $\Lambda/F$ is unramified outside $\{\infty \} \cup (S \backslash T)$ and
\item $\Gal(\Lambda/F) \simeq \widehat{X}^r$ as $\F_\ell[\Delta]$-module.
\end{enumerate}
Let $r_{E}(S)$ be the multiplicity of $\widehat{X}$ in $\Gal(\Lambda_E(S)/F)$ and $\Gamma_E(S) = \Gal(\Lambda_E(S)/\Q).$
\end{Not} 

We introduce  invariants of nuggets over $\Z_S$ which have $\cE$ as subquotient, where $\cE$ is an $\F$-module scheme with Galois module $E$ and $\F = \F_\gl$.  Let $\cZ \simeq  \cZ_\gl^n$ and let  $0 \to \cZ \to \cV \to \cE \to 0$ be an  exact sequence of acceptable $\F$-module schemes over $\Z_S.$ 
Put $G = \Gal(\Q(V)/\Q)$ and let  $[c]$ in $H^1(G, \Hom_\F(E,Z))$ be the obstruction to splitting of the Galois module sequence: 
\begin{equation}\label{ZVE}
0 \to Z \to V \stackrel{\pi}{\to} E \to 0. 
\end{equation}
Remark  \ref{ForD1} and Lemma  \ref{ram} imply that  $L=\Q(V)$ is contained in $\Lambda_E(S)$.

The next two lemmas contain local conditions at $\ell$ and the primes $p$ dividing $N_E$ implied by semistability of $V$. 

\begin{lem} \label{SelGro}
Let $\cI_v(L/\Q) = \langle \sigma_v \rangle$ and $M_v = (\sigma_v-1)(E)$  for $v$ over $p \, | \, N_E$.   The exact sequence 
$
0 \to Z \to \pi^{-1}(M_v) \to M_v \to 0
$ 
consists of trivial $\F[\cI_v]$-modules.  If $\gf_p(V) = \gf_p(E)$, it is $\F[\cD_v]$-split.
\end{lem}

\proof
It is clear that $\pi^{-1}(M_v) = (\sigma_v-1)(V) + Z$ and $\cI_v$ acts trivially because $(\sigma_v-1)^2(V) = 0$.  If $\gf_p(V) = \gf_p(E)$, then $\pi$ induces an isomorphism of the $\F[\cD_v]$-modules $(\sigma_v-1)(V)$ and $M_v$, since they have the same $\F$-dimension.  This gives us the $\F[\cD_v]$-splitting.  \qed

\begin{lem} \label{EllCondition}
Let $0 \to \cZ \to \cV \to \cX \to 0$ be an exact sequence of $\F$-module schemes over $\Z_\ell$ with $\cZ$ \'etale.  Fix $\lambda$ over $\ell$ in $\Q(V)$ and consider the exact sequences of $\cD_\lambda$-modules:
$0 \to Z \to V \stackrel{\pi}{\to} X \to 0$  and  $0 \to Z \to \pi^{-1}(X^0) \to X^0 \to 0.$ 
The first is $\F[\cI_\lambda]$-split and the second is $\F[\cD_\lambda]$-split.
\end{lem}

\proof
The second sequence splits because $\cZ^0 = 0$, so $\pi^{-1}(X^0) = Z + V^0$ is a direct sum.  Now let $j: V \to V^{et}$ be the natural map with kernel $V^0$.  Since the \'etale sequence $0 \to j(Z) \to V^{et} \to X^{et} \to 0$ consists of trivial $\cI_\lambda$-modules, we can find an $\cI_\lambda$-submodule $W$ of $V^{et}$ such that $V^{et} = j(Z) + W$ is a direct sum, with $W \simeq X^{et}$.  It is easy to check that $V = j^{-1}(W) + Z$ is a direct sum and this shows that the first sequence is $\F[\cI_\lambda]$-split. \qed

\medskip

The extension problem (\ref{ZVE}) has   a Selmer interpretation. For a Galois extension $K/\Q$ with 
$\Lambda_E(S)\supseteq K \supseteq F=\Q(E),$ we define
\begin{equation} \label{TSel}
H^1_\cL(\Gal(K/\Q), \widehat{E}) = \ker:     
       H^1(\Gal(K/\Q),\widehat{E}) \stackrel{res}{\longrightarrow}
       \prod_{v | \ell N_E} {\cL_v},
\end{equation}
where $\cL_v = \begin{cases} 
 H^1(\cI_v(K/\Q), \widehat{M}_v) &\text{if } v \, | \, N_E, \\     
 H^1(\cI_v(K/\Q),\widehat{E}) \times H^1(\cD_v(K/\Q), \widehat{E^0}) & \text{if } v \, | \, \ell.\! \end{cases} $

\begin{cor} \label{selcor}
 In {\rm (\ref{ZVE})}, there is a submodule $Z'$ of $Z$ such that the exact sequence 
$0 \to Z/Z' \to V/Z' \to E \to 0$ is $\F[G]$-split and $\dim_\F Z' \le \dim_\F H^1_\cL(G,\widehat{E})$.
\end{cor}

\proof
Apply Corollary \ref{splitcor}(ii) with $X = E$ and $Y = Z$, using Lemmas  \ref{SelGro} and  \ref{EllCondition} for the local conditions over $N_E$  and  $\ell$, respectively.  \qed

\begin{lem} \label{infres}
Let $s_{E}(S)= \dim_\F  H^1_\cL(\Gal(\Lambda_E(S),\widehat{E})$. Then
$$
s_{E}(S) \le \dim_\F   H^1_\cL(\Delta,\widehat{E}) + r_{E}(S)   \dim_\F   \End_{\F[\Delta]} \widehat{E}.
$$
If $\cI_\lambda(F/\Q)$ contains an $\ell$-Sylow subgroup of $\Delta$ for $\lambda\,|\,\ell$, then $H^1_\cL(\Delta,\widehat{E}) = 0$.
\end{lem}

\proof
Let $\Lambda = \Lambda_E(S)$ and $\Gamma = \Gal(\Lambda/\Q)$.  By inflation-restriction, we have 
\begin{equation} \label{infresseq}
0 \to H^1_\cL(\Delta,\widehat{E}) \to H^1_\cL(\Gamma,\widehat{E}) 
\to \Hom_{\F_\ell[\Delta]}(\Gal(\Lambda/F),\widehat{E}),
\end{equation}
since $\Gal(\Lambda/F)$ acts trivially on $\widehat{E}$.  Lemma  \ref{BC} gives us $\F_\ell[\Delta]$-isomorphisms
$$
\Hom_{\F_\ell[\Delta]}(\Gal(\Lambda/F),\widehat{E}) \simeq \Hom_{\F_\ell[\Delta]}(\widehat{X}^r,\widehat{X}^a) \simeq
(\End_{\F_\ell[\Delta]} \, \widehat{X})^{ra}, 
$$
where $r = r_E(S)$.  Since $ab = [\F:\F_\ell]$, Lemma \ref{BC}(iii) now shows that  
$$
\dim_\F \Hom_{\F_\ell[\Delta]}(\Gal(\Lambda/F),\widehat{E}) = r \, \dim_\F \End_{\F[\Delta]} \widehat{E}.
$$

Suppose that $\cI_\lambda(F/\Q)$ contains an $\ell$-Sylow subgroup $P$ of $\Delta$.   Any element $[c]$ of $H^1_\cL(\Delta,\widehat{E})$ restricts to 0 in $H^1(\cI_\lambda(F/\Q), \widehat{E})$, so vanishes on further restriction to $H^1(P,\widehat{E})$.  But then $[c] = 0$ because $H^1(\Delta,\widehat{E}) \stackrel{\res}{\longrightarrow} H^1(P, \widehat{E})$ is injective, cf.\! \cite[Ch. IX, \S 2, Thm.\! 4]{Ser1}.  \qed

\medskip

We introduce two invariants to estimate the dimension  of a non-prosaic nugget. Recall our standard assumption that $A$ is an abelian variety of $\go$-type and $\gl$ is a prime of $\go$ above the prime $\ell$ of good reduction for $A$.

\begin{Def} \label{deficiency}
Let $\gW(\cE)$ be the set of nuggets $\cW$ that are subquotients of $A[\gl^{\infty}]$, have the exceptional $\cE$ as  constituent and  satisfy  $N_W = N_E.$  Put 
$$
\delta_A(\cE) = \max_{\cW \text{ in } \gW(\cE)}(\dim_\F \cW-\dim_\F \cE).
$$
For a fixed $E$ in $\gS_{\gl}(A)$, the {\em deficiency} is given by $\delta_A(E) := \max \, \delta_A(\cE)$, where $\cE$  has Galois module $E$.  We omit $A$ when it is clear from the context.
\end{Def}
 
We say the exact sequence
$
0\to\cU\to\cV\to\cW\to 0
$
is {\em generically split} if the associated exact sequence of Galois modules splits.  In view of (\ref{specialeq}), for any $\cW$ in $\gW(\cE),$  we have exact sequences
\begin{equation} \label{ep2eq}  
0 \to \cZ \to \cV \to \cE \to 0 \quad \text{and} \quad 0 \to \cV/\cZ \to \cW/\cZ \to \cM \to 0, 
\end{equation}
with $\cZ$ and $\cM^D$ constant by Lemma \ref{AllZnug}.

\begin{Def} \label{ep2def} 
Let $\gW_{spl}(\cE)$ consist of those $\cW$ in $\gW(\cE)$ for which both exact sequences in (\ref{ep2eq}) are generically split.  Define $\epsilon_\gl(\cE) = \max \dim_{\F} W-\dim_{\F}  E $ over $\cW$ in $\gW_{spl}(\cE).$  For a fixed $E$ in $\gS_{\gl}(A)$, let  $\epsilon_\gl(E) = \max \epsilon_\gl(\cE)$, taken over $\cE$ with $E$ as Galois module.
\end{Def}

When $\ell$ is odd, generic splitting implies splitting as group schemes, so $\epsilon_\gl(\cE)=0.$  See \S \ref{except} for bounds on $\epsilon_\gl(\cE)$ when $\ell= 2$.  With Notation \ref{seldim}, we have the following bound on the deficiency $\delta_A(E)$.

\begin{prop} \label{crude}
Let  $\Gamma_E = \Gal(\Lambda_E(T)/\Q)$ and $s_E = \dim_\F  H^1_\cL(\Gamma_E,\widehat{E})$.  Then
$
\delta_A(E) \le s_E+ s_{E^*} +\epsilon_\gl(E).
$ 
\end{prop}

\proof
For $\cW$ in $\gW(\cE)$, consider the first exact sequence of (\ref{ep2eq}).  Let $L = \Q(V)$ and $G = \Gal(L/\Q)$.  Since inflation $H^1_\cL(G,\widehat{E}) \to H^1_\cL(\Gamma_E,\widehat{E})$ is injective, $\dim H^1_\cL(G,\widehat{E}) \le s_E$.   Then, by Corollary \ref{selcor}, there is a subspace $Z^\prime$ of $Z$ such that
$
0 \to Z/Z^\prime \to V/Z^\prime \to E \to 0
$
is  $\F[G]$-split exact and $\dim_\F Z^\prime \le s_{E}$.   Write $\cZ_1$ (resp.\! \negmedspace $\cV_1$, $\cW_1$) for the quotient of $\cZ$ (resp.\! \negmedspace ${\cV}$, $\cW$) that corresponds to $Z/Z'$ (resp.\! \negmedspace $V/Z^\prime$, $W/Z'$).  Then $\cW_1$ is a nugget with special filtration $0 \subseteq \cZ_1 \subset \cV_1 \subseteq \cW_1$ and $N_{W_1} = N_E$.  Moreover, 
$
0 \to \cZ_1 \to \cV_1 \to \cE \to 0
$
is generically split.

Passing  to Cartier duals on 
$
0 \to \cV_1/\cZ_1 \to \cW_1/\cZ_1 \to \cM \to 0,
$
we find $\F$-module subschemes $\cM' \subseteq \cM$ and $\cV_1 \subseteq \cW' \subseteq \cW_1$, with $\dim_\F M/ M' \le s_{E^*}$, such that 
$
0 \to \cV_1/\cZ_1 \to \cW'/\cZ_1 \to  \cM' \to 0
$ 
is generically split.  It follows that $\cW'$ is in $\cW_{spl}(\cE)$, so $\dim \cW' - \dim E \le \epsilon_\gl(E)$ by Definition\! \ref{ep2def}.   The claim  now ensues from $\dim \cW \le \dim \cW ' + s_E + s_{E^*}$. \qed

\begin{Rem}\label{onearmed} 
\hfill
\begin{enumerate}[i)]
\item If $\cV$ is a ``one-sided nugget" with
$
0 \to \cZ \to \cV \to \cE \to 0
$
and  $N_V = N_E$, our  proof  gives $\dim Z \le s_E + \epsilon_\gl(E).$ 
\item Since $N_V = N_E$ in the proof above,  Lemma  \ref{SelGro} implies the stronger local condition $\cL_v = H^1(\cD_v,\widehat{M}_v)$ at places $v$ over $N_E$.\end{enumerate}
\end{Rem}

\begin{lem} \label{pretrivial} 
Let $\cW$ be a nugget and $\gf(W)=\sum_p \gf_p(W).$  If $W$ has an exceptional constituent $E$, then $\dim_\F W -\gf(W)\le \dim_\F E  - \gf(E) + \delta (E)$.  If $W$ is prosaic, then $\dim_\F W \le \gf(W) + 1$, with equality only if some core has conductor $p \equiv \pm 1 \, (\tilde{\ell})$.   
\end{lem}
\proof
Suppose the lemma is false and choose a counterexample $\cW$ of 
minimal dimension.   We have exact sequences as in  (\ref{specialeq}).

By definition, $\cZ$ is filtered by copies of $\cZ_\gl$, so $G = \Gal(\Q(Z)/\Q)$ is an $\ell$-group.  Let $\ga_G$ be the augmentation ideal in $\F[G]$ and $r$ the least  integer  such that $\ga_G^r Z=0$.  If $r \ge 2$,  some prime $p$ occurs in the conductor of $\ga_G^{r-2}Z$ by Lemma \ref{AllZnug}.  Let $\sigma_v$ generate inertia at a place $v$ above $p.$ There is an element $z_2$ in $\ga_G^{r-2}Z$ such that $z_1=(\sigma_v-1)z_2 \ne 0$.   Since $z_1$ is in $\ga_G^{r-1}Z$, it generates a trivial Galois module.  Let $\cZ_1$ denote the corresponding $\F$-module subscheme of $\cZ$ and let $\cW' = \cW/\cZ_1$.  Then $\dim_\F \cW'=\dim_\F \cW -1$ and $\gf(W') \le \gf(W)-1$, so
\begin{equation} \label{newW}
\dim_\F W' - \gf(W') \ge \dim_\F W - \gf(W),
\end{equation}
and $\cW'$ would be a smaller  counterexample. Thus $Z$ has trivial action so that  $\cZ\simeq \Zl^a$ is constant of exponent $\ell$.  Upon passing to  Cartier duals, we find similarly that $\cM^D$ is constant and  $\cM\simeq \Ml^b.$ 

Assume  $E$ is non-zero and  set $\overline{\cW}=\cW/\cZ.$
We claim that $N_W = N_{\overline{W}}.$  Otherwise, the conductor exponents of $W$ and $\overline{W}$ differ at some place $v$.  Since Galois acts trivially on $Z,$ there is a non-zero element $z$ in $(\sigma_v-1)(W) \cap Z$.  Let $\cW'$ be the $\F$-module scheme quotient of $\cW$ corresponding to the Galois module $W/\langle z \rangle$.  Then (\ref{newW}) holds for $\cW'$ violating minimality of  $\cW$.  A similar argument with the Cartier dual of the sequence $0  \to \cE \to \overline{\cW} \to \cM \to 0$ implies that $N_{\overline{W}}= N_E.$  So  $N_W = N_E$ and  $\cW$ is not a counterexample by Definition\! \ref{deficiency}. 

If $\cW$ is  prosaic, i.e.\! $E = 0$, we use the argument above, the nugget property and minimality to show that  $\dim_\F Z = \dim_\F M = 1$.  Then $\cW$ is  a core, for which the claim  was established in  Corollary \ref{core}.  \qed

\subsection{Better bounds for  $\delta(E)$} \label{except}
Keep the notation of \ref{seldim} and \ref{ZVE}, with $\ell=2.$  For each $\lambda$ over 2 and group scheme $\cE$,  we have the associated connected, biconnected and \'etale $\cD_\lambda$-modules $E^0$, $E^b$ and $E^{et}.$  Since $\lambda$ is unramified in the elementary 2-extension $L/F$, we have $\cD_\lambda(L/F) = \langle h \rangle$, with $h^2 = 1$.

\begin{lem}  \label{ep2bound}
Let $\gb_\lambda$ be the augmentation ideal in $\F[\cI_\lambda]$.  For $\gl\,|\,2$ in $\go$, we have 
$$
\epsilon_\gl(E) \le \min \{ \dim_{\F}E^{\cD_\lambda}, \dim_{\F}E^{\cI_\lambda} - \dim_{\F}(\gb_\lambda E)^{\cI_\lambda} \}.
$$
If $\cI_\lambda$ acts on $E$ via a non-trivial 2-group, then $\epsilon_\gl(E) \le \dim_{\F}E - 2$.
\end{lem}

\proof
If $\cW$ is in $\gW_{spl}(\cE),$ then the second sequence in (\ref{ep2eq}) is generically split, so there is a Galois submodule $X$ of $W$ with $E_1 = W/X\simeq E$ and the sequence $0 \to \cX \to \cW \to \cE_1 \to 0$ is exact.  Thus $\cX$ is a constant group scheme by Lemma  \ref{AllZnug}.  Taking multiplicative subschemes at $\lambda$, we find that $\cW^m \simeq \cE_1^m$.  

Similarly, we have  $0 \to \cE_2 \to \cW \to \cY \to 0$, with $\cY^D$ constant and $E_2 \simeq E$.  Taking multiplicative subschemes at $\lambda$ gives
\begin{equation} \label{mseq}
0 \to \cE_2^m \to \cW^m \to \cY \to 0.
\end{equation}
Hence 
$
\dim W - \dim E = \dim  Y = \dim  W^m - \dim  E_2^m=\dim  E_1^m - \dim  E_2^m.
$

We have $\dim  E_1^m \le \dim  E^{\cI_\lambda}$ and  $\gb_\lambda E \subseteq E_2^0$ because $\cI_\lambda$ acts trivially on $E_1^m$ and $E_2^{et}$.   Moreover, by \cite{Ray1},  the tame ramification group acts non-trivially on the simple constituents of $E_2^b$ over the strict Henselization, so $(E_2^0)^{\cI_\lambda} \subseteq E_2^m.$ The inequality $\dim E_2^m \ge \dim (\gb_\lambda E)^{\cI_\lambda}$ gives $\epsilon_\gl(E) \le \dim  E^{\cI_\lambda} - \dim  (\gb_\lambda E)^{\cI_\lambda}$.  In particular, when $\cI_\lambda$ acts on $E$ through a non-trivial 2-group, we have $\epsilon_\gl(E) \le \dim  E -2$ .

The isomorphism $\cE_1^m \stackrel{\sim}{\to} \cW^m$ and the surjection $\cW^m \twoheadrightarrow \cY$ in  (\ref{mseq}) yield a surjection of $\cD_\lambda$-modules $\cE_1^m \twoheadrightarrow Y$.  Since $Y$ is a trivial Galois module, this map induces a surjection $E_1^m/\ga_\lambda E_1^m \twoheadrightarrow Y$, where $\ga_\lambda$ is the augmentation ideal in $\F[\cD_\lambda]$.  But $\cI_\lambda$ acts trivially on $E_1^m$ and so $\cD_\lambda$ acts via the group generated by a Frobenius $\Phi.$  Hence 
\begin{eqnarray*}
\dim_\F  W - \dim_\F  E = \dim_\F  Y &\le& \dim_\F  E_1^m/\ga_\lambda E_1^m\\
&=& \dim_\F  \; (E_1^m)^{\langle \Phi \rangle} = \dim_\F  \; (E_1^m)^{\cD_\lambda} \le \dim_\F  E^{\cD_\lambda}.    \qed
\end{eqnarray*} 

\begin{cor} \label{ep2boundcor} 
If $\Delta \subseteq \SL_{\F}(E)$, then
$ 
\epsilon_\gl(E)\le \begin{cases} 
  \dim_\F E &\text{if } \cD_\lambda=1, \\
  \dim_\F E - 1 &\text{if } |\cD_\lambda| = 2, \cI_\lambda=1, \\
  \dim_\F E-2 &\text{otherwise}.
                   \end{cases}
$  \end{cor}

\proof If the claim is false, the lemma implies that $\dim_\F E^{\cD_\lambda} =\dim_\F E-1$.  Then  $\cD_\lambda$ is an elementary 2-group, since $\Delta \subseteq \SL_\F(E).$  By the lemma, we now find that $\cI_\lambda=1$.  Thus $\cD_\lambda$ is cyclic and so has order 2.  \qed

\begin{Def} \label{S-transp} 
We say $\cE$ is {\em $(S \backslash T)$-transparent} if $\Ext^1_{\Z_S}(\cE,\Zl)=0$, where $S \supseteq T$ and $T=T_E$ is the set of bad primes of $E.$   When  $S=T,$ we simply say transparent.
\end{Def}

\begin{lem} \label{Newse0}
Let $E$ be a self-dual exceptional $\F[G_\Q]$-module  and $S \supseteq T_E$.  If $\dim_\F E=2$, $H^1_\cL(\Delta,E) = 0$ and $r_E(S) = 0$, then   
\begin{enumerate} [{\rm i)}]
   \item $\delta(E) = 0$ and $E$ is $(S/T_E)$-transparent if either $|\cD_\lambda(F/\Q)| = |\cI_\lambda(F/\Q)| = 2$ or $|\cD_\lambda(F/\Q)| \ge 3$.
 \item $\delta(E) \le 1$ if $|\cD_\lambda(F/\Q)| = 2$ and $\cI_\lambda(F/\Q) = 1.$
 \item $\delta(E) \le 2$ in all other cases.
\end{enumerate}         
\end{lem}

\proof
Since $E$ is self-dual, it affords a representation whose determinant is the mod 2 cyclotomic character \cite{Rib3}  and so $\Delta$ is contained in $\SL_2(\F)$.  Now use Lemmas \ref{infres}, \ref{ep2boundcor} and Proposition  \ref{crude}.  \qed

\medskip

The restriction $\tc = \res [c] : \Gal(L/F) \to \Hom(E,Z)$ is an $\F[\Delta]$-homomorphism and $\tc_h: E \to Z$ is $\F$-linear, cf.\! Remark  \ref{ForD1}.

\begin{lem} \label{ckill} 
Let $\ga_\lambda$ be the augmentation ideal in $\F[\cD_\lambda(F/\Q)]$. Then $\tc_h$ vanishes on $E^0 + \ga_\lambda E$ and 
$
\dim_\F \tc_h(E) \le \dim_\F \, (E^{et})^{\cD_\lambda(F/\Q)}.
$
\end{lem}

\proof
Since $\tc_h(E^0) = 0$ by Lemma  \ref{EllCondition},  $\tc_h$ factors through $E^{et}$.  Also, $\cI_\lambda$ acts trivially on $E^{et},$ so $\ga_\lambda E^{et} = (\Phi-1)(E^{et})$, for any  Frobenius $\Phi$  in $\cD_\lambda(F/\Q)$.   We know that $\Phi$ acts trivially on $Z$ and $h$ is a power of $\Phi$.  It follows that 
$$
\tc_h(\Phi \bar{e}) = \Phi^{-1}(\tc_h(\Phi \bar{e})) = (\Phi^{-1}(\tc_h))(\bar{e}) = \tc_{\Phi^{-1}(h)}(\bar{e}) = \tc_h(\bar{e})
$$
for all $\bar{e}$ in $E^{et}$.  Hence $\tc_h$ vanishes on $\ga_\lambda E^{et}$ and so it factors through 
$
E^{et}/\mathfrak{a}_\lambda E^{et}.
$  
This last space has the same dimension as $(E^{et})^{\Phi}$.  \qed

\begin{lem} \label{coc}
If the residue degree $f_\lambda(F/\Q)$ is even and $\cD_\lambda(F/\Q)$ acts on $E^{et}$ through a quotient of odd order, then the primes over 2 split completely in $L/F$. 
\end{lem}

\proof
If $f_\lambda(F/\Q)$ is even, then $h$ is a square in $\cD_\lambda(L/\Q)$, say $h = g^2$,  with $g$ chosen to have order a power of 2.  Hence $g$ acts trivially on $E^{et}$ and  $(1+g)(E) \subseteq E^0$.  By Lemma \ref{EllCondition}, we may arrange for the cocycle $c\!: \, G \to \Hom(E,Z)$ to satisfy $c_g(E^0) = 0$.  Then,   
by the cocycle identity, we have
$$
\tc_h(e) = c_{g^2}(e) = ((1+g)c_g)(e) = c_g((1+g^{-1})(e))\,  \in\,   c_g(E^0) = 0
$$
for all $e$ in $E$.  But $\tc$ is injective, so $h = 1$. \qed

\medskip

We now impose the following hypotheses, with notation from \S4.3.

\medskip

\noindent\fbox{\bf D} {\bf 1}. $E$ is 2-dimensional over $\F$, irreducible and self-dual as $\F[G_\Q]$-module.

\hspace{7 pt}{\bf 2}. The generalized Selmer group $H^1_\cL(\Delta,\widehat{E})$ is trivial.

\hspace{7 pt}{\bf 3}. There is an $\F_2[\Delta]$-isomorphism $\Gal(\Lambda/F) \simeq \widehat{X}$, that is  $r_E = 1$.

\hspace{7 pt}{\bf 4}. The primes over $2$ do not split completely in $\Lambda/F$, so  $\cE$ is not biconnected.

\begin{Rem} \label{Ds}
Under {\bf D2}, the cohomological restriction map $[c] \mapsto \tc$ is injective, so in (\ref{ZVE}), $V$  splits if and only if $\tc = 0$.   By {\bf D3}, $L = \Q(V)$ is equal to $F$ or $\Lambda$.  If $L = F$, then $\tc = 0$ and (\ref{ZVE}) splits, while if $L = \Lambda$, $\tc$ induces the isomorphism in {\bf D3} and $\cD_\lambda(L/F) = \langle h \rangle$ has order 2 by {\bf D4}.  By irreducibility of $X$, $h$ generates $\Gal(L/F)$ as $\F_2[\Delta]$-module.  Let $\cZ'$  be the $\F$-module subscheme of $\cZ$ corresponding to
$
Z' = \sum \{\,c_\gamma(E) \,|\, \gamma \in \Gal(L/F) \, \} = \tc_h(E).
$
Remark  \ref{ForD1} shows that the following sequence is generically split:
\begin{equation} \label{barZVE}
0 \to \cZ/\cZ' \to \cV/\cZ' \to \cE \to 0
\end{equation}
\end{Rem}

\begin{lem}\label{bb}
Assume {\rm {\bf D}} and residue degree $f_\lambda(F/\Q) = 2$.  Then $\delta(\cE) \le 1,$ unless  there is a nugget $\cW$ with $\gr \cW = [\cZ_\gl \, \cE \Mu_\gl]$ and $\dim \cE^{et} \ne 1$, in  which  case, $\delta(E) \le 2.$ The latter can only happen if $\dim_\F A[\gl]\ge 6.$
\end{lem}

\proof 
Let $\cV$, $\cW$ be nuggets as in (\ref{ep2eq}) and $L = \Q(V)$.  If $\dim \cE^{et} \le 1$, then $\cD_\lambda$ acts on $E^{et}$ via a subgroup of $\F^\times$ and so the primes over 2 split in $L/F$ by Lemma  \ref{coc}.  By {\bf D4}, $L = F$ and (\ref{ZVE}) splits.  By  Lemma  \ref{ep2bound}, we have 
$$
\dim Z = \dim V - \dim E \le \epsilon_2(E) \le \dim E^{\cD_\lambda} \le 1.
$$ 

If $\dim \cE^{et} = 2,$ Lemma  \ref{ckill} shows  $\dim Z' \le 1$ in the generically split sequence (\ref{barZVE}), which must split as an exact sequence of schemes by Lemma \ref{Mazur}.  Since $\cV$ is a nugget, we have $Z/Z' = 0$, so $\dim_\F \cZ \le 1$ in all cases.  By Cartier duality, $\dim_\F \cM \le 1$.  Hence $\delta(E) \le 2$, with equality only if  $\gr \cW = [\cZ_\gl \, \cE \Mu_\gl]$ for some nugget $\cW$.

When $\dim \cE^{et} = \dim \cE^0 = 1$, we have seen that (\ref{ZVE}) is generically split and so is  $0 \to \cV \to \cW \to \cM \to 0$, by a dual argument.  We have   $\delta(\cE) \le 1$, since  Lemma  \ref{ep2bound} gives
$
\dim \cW - \dim \cE \le \epsilon_2(E) \le \dim E^{\cD_\lambda} \le 1.
$

If $\dim A[\gl] = 4,$ $\dim \cE^{et} = \dim \cE^0 = 1$ since  $\cW = A[\gl]$ has as many $\Mu_\gl$ as $\cZ_\gl.$   \qed

\bigskip

For the final result of this section, we need some facts about representations that respect a flag of $\F_2[G_\Q]$-modules $0 \subset Z \subset V \subset W$, with $\gr W=[\F_2 E \, \F_2].$ We assume $\dim_{\F_2} E = 2$ and $\Delta = \Gal(\Q(E)/\Q) \simeq \SL_2(\F_2)$.  

For $x = (a,b)^t$ a column vector in $E$, consider the $\Delta$-invariant quadratic form $Q(x) = a^2 + ab + b^2$ and define $x^\dagger = (b,a)$.  Then $(x,y) \mapsto x^\dagger y = \det(x,y)$ is the symplectic form on $E$ associated to $Q$ by $Q(x+y) = Q(x) + Q(y) + x^\dagger y$.  Setting 
$$
\iota(x,\delta,c) = \left[ \begin{smallmatrix}
  1 & x^\dagger \delta & c \\ 0 & \delta  & x \\ 0 & 0 & 1
                          \end{smallmatrix}\right],
$$
we have $\iota(x,\delta_1,c_1) \, \iota(y,\delta_2,c_2) = \iota(x+\delta_1y, \, \delta_1\delta_2,\, c_1+x^\dagger\delta_1 y + c_2)$.

Let $\cP=\{\iota(x,\delta,c)\,|\, \delta\in\Delta,\, c\in\F_2\}$ and   
$\cP_1 = \{ \, \iota(x,\delta,Q(x)) \, | \, x \in E, \, \delta \in \Delta \}.$ Then 
$\cP \simeq \cP_1 \times \langle \xi \rangle$, where $\xi = \iota(0,I_2,1)$ is the central involution in $\cP$.  The normal subgroup $H = \{ \iota(x,I_2,Q(x)) \, | \, x \in E \}$ of $\cP_1$ is $\Delta$-isomorphic to $E$ under the action of $\tilde{\delta} = \iota(0,\delta, 0)$ by conjugation. The relation $(\delta x)^\dagger = x^\dagger \delta^{-1},$ implied by the  $\Delta$-invariance of $Q,$  gives
$
\tilde{\delta} \, \iota(x,1,Q(x)) \, \tilde{\delta}^{-1} = 
       \iota(\delta x,1,Q(\delta x)).
$

Let $\tilde{\Delta}= \iota(0,\Delta,0).$ Then  $\cP_1 = H \tilde{\Delta} \simeq \cS_4$  is a Coxeter group, generated by 
$$
\tau_1=\iota\left(0, \left[\begin{smallmatrix} 1 & 1 \\ 0 & 1 \end{smallmatrix}\right] ,0\right), \quad
\tau_2=\iota\left(0, \left[\begin{smallmatrix} 1 & 0 \\ 1 & 1 \end{smallmatrix}\right] ,0\right), \quad
\tau_3=\iota\left( \left[\begin{smallmatrix} 1 \\ 1 \end{smallmatrix}\right] ,  \left[\begin{smallmatrix} 0 & 1 \\ 1 & 0  \end{smallmatrix}\right] , 1 \right),
$$
three involutions whose pairwise products have order 3.

\begin{lem}\label{specialshape}
If $E$   satisfies {\rm {\bf D}}, with $\F=\F_2$ and  $|\cI_\lambda(F/\Q)| = 2$, then $\delta(E) \le 1$.
\end{lem}

\proof
Note that {\bf D1} and {\bf D2} follow from the other two.  In fact,  $r_E =1$ implies $s_E \le 1$ by Lemma \ref{infres}.  Let $\cW$ be a nugget with $N_W = N_E$ as in Definition\! \ref{deficiency}.  Lemma \ref{ep2bound} implies that $\epsilon_2(E)=0$.  It follows from Remark  \ref{onearmed} that the \'etale subscheme $Z$ in (\ref{specialeq}) satisfies  $\dim Z \le 1$, with equality only if $\Q(V)=\Lambda$.  By a dual  argument, we have $\dim M \le 1$, with equality only if $\Q(W/Z)=\Lambda$.  

Assume  $\dim Z = \dim M = 1.$  We build a basis for $W$ reflecting the local structure of $\cW$ at $\lambda$.  Because $|\cI_\lambda(F/\Q)| = 2$, we have $\dim E^0 = 1$ and so $\dim V^0 = 1$.  Let $b_1$ generate $Z,$ $b_2$ generate $V^0.$ Extend to bases $ b_1, b_2, b_3 $ for $V$ and  $b_2, b_4 $ for $W^0$.  Write $\rho_W$ for the matrix representation of $\Gamma = \Gal(\Q(W)/\Q)$ afforded by the basis $b_1, b_2, b_3, b_4$.  The images of the induced representations $\rho_V$ on $V$ and  $\rho_{W/Z}$ on $W/Z$ are both isomorphic to $G = \Gal(\Lambda/\Q)$.  Also, $\rho_V(g_1) = \rho_V(g_2)$ if and only if $\rho_{W/Z}(g_1) = \rho_{W/Z}(g_2)$.

The inertia group $\cI_\lambda(\Q(W)/\Q) = \langle \sigma_\lambda \rangle$ is cyclic of order 2.  Both $W^0$ and $W^{et}$ are unramified $\cD_\lambda$-modules.  The Frobenius $\Phi = \Phi_\lambda$ in $\Gal(\Lambda/F)$  is non-trivial on $\Lambda$ and we have
 $$ 
\rho_W(\Phi) = \iota(  \left[\begin{smallmatrix} 1 \\ 0 \end{smallmatrix}\right],I_2,0) \quad  \text{  and  } \quad 
\rho_W(\sigma_{\lambda})=\tau_1.
$$ 
For each place $v$ over a prime $p$ dividing $N_W = N_E$, let $\sigma_v$  generate $\cI_v(\Q(W)/\Q)$.  Since the conductor exponent $\gf_p(E) = 1$, $\rho_W(\sigma_v)$ is a transvection in $\SL_4(\F_2)$ and  $\rho_V(\sigma_v)$ becomes a transposition under the isomorphism $\Image \rho_V \simeq \cS_4$.   By conjugation, we may produce any transvection in the upper left $3 \times 3$ corner by a suitable choice of  $v$ and so  assume that 
$$
\rho_W(\sigma_v) = \left[ \begin{smallmatrix}
                           1&0&0&0\\0&1&0&0\\0&1&1&s_v\\0&0&0&1
           \end{smallmatrix} \right] .
$$  
Over each bad $q,$ there is a $w$  with $\rho_W(\sigma_w)$ of the same shape for some $s_w$.
 Since $\rho_V(\sigma_v) = \rho_V(\sigma_w)$, we have $\rho_{W/Z}(\sigma_v) = \rho_{W/Z}(\sigma_w)$ and so $s_v = s$ is independent of $p.$  Replacing $b_4$ by $b_4^\prime = b_4 + s  b_2$ preserves $\rho_W(\Phi)$ and $\rho_W(\sigma_\lambda)$, but makes $\rho_W(\sigma_v) = \tau_2$.  Thus the group $\tilde{\Delta},$ generated by $\tau_1$ and $\tau_2,$ is contained in $\Gamma$.   
 
 We claim that $\Gamma\subseteq\cP_1.$ 
But  $\Gamma$ is generated by  its inertia groups, and so by $\Gamma$-conjugates of $\tilde{\Delta}.$ Since $\tilde{\Delta}\subseteq\cP_1$, it suffices to show $\Gamma\subseteq \cP$. Let  $g$ in $\Gamma$ fix $F$, say 
$$
\rho_W(g) = \left[ \begin{smallmatrix} 1& x^t & c \\0&I_2&y \\ 0&0&1        
                   \end{smallmatrix} \right],
$$
with $c$ in $\F_2$, $x^t$ and $y$ in $E.$ Choose $\delta$ in $\Delta$ so that $\delta (1,0)^t = y$ and let $h = \tilde{\delta} \Phi \tilde{\delta}^{-1}.$
Then 
$\rho_{W/Z}(g) = \rho_{W/Z}(h)$ and so  $\rho_V(g) = \rho_V(h).$ Thus, $\rho_W(g)$ and $\rho_W(h)$ agree up to an element of $\langle \xi\rangle$.  Since $\rho_W(h) = \iota(y,I_2,0)$ is in $\cP$, we have $\Gamma\subseteq\cP$.   But $\rho_W(\Phi)$ is not in $\cP_1$, a contradiction.   \qed

\section{General bound} \label{GenBd} \numberwithin{equation}{section} Let $A$ be a semistable abelian variety of $\go$-type. Our aim here is to bound $\epsilon_0(A[\gl])$, the number of one-dimensional constituents  in a composition series for $A[\gl]$ as an $\F[G_\Q]$-module. 
\begin{Def} \label{Nss}
The {\it semisimple conductor} of  $A[\gl]$ is $N_A^{ss}(\gl) =\prod N_E,$
where $E$ runs over the multiset $\gS_\gl(A)$.  Its {\em prosaic conductor} is $N_A^u(\gl) = N/ N_A^{ss}(\gl)$, where $N = N_A^0$ is the reduced conductor of $A$.  When $\gl$ is clear from the context, it may be omitted.  Write  $\Pi_A^u$ for the set of prime factors of $N_A^u.$ 
\end{Def}
The prosaic conductor depends only on $\gl$ and the $\go$-isogeny class of $A$, since this is true of $N_A^{ss}(\gl)$ by Proposition \ref{JH2}, but  it is not the conductor of a  Galois module naturally associated to $A$.

\begin{lem} \label{pretrivialZ} 
Let $A$ have good reduction at a prime $q \neq \ell$ and let $s_0$  be the greatest integer in
$
2 \dim A \cdot \log (1+\sqrt{q}) /\log |\F_\gl|.
$
Let $\cZ$ {\rm(}resp.\! $\cM${\rm)} be an $\go$-module scheme subquotient of $A[\gl^n]$ filtered by copies of $\Zl$ {\rm(}resp.\! $\Ml${\rm)}.  Then 
$$
\len_\go Z  \le  \gf(Z)+s_0 \quad \quad \text{and} \quad \quad \len_\go M  \le  \gf(M)+s_0.
$$
If $\End A = \go,$ then
$
\max\{ \len Z , \, \len M \}  \le  s_1,
$ 
where $s_1$ is the number of isomorphism classes in the $\Q$-isogeny class of $A.$
\end{lem}

\proof 
By duality, we only  prove the assertions about $Z$.  Replacing $A$ by a  quotient, suppose $\cZ \subseteq A[\gl^n]$.  The result holds when $Z$ has trivial Galois action because $A(\Q)[\ell^{\infty}]$ injects into $\tilde{A}(\F_q)$ by specialization and  $s_0$ is the Weil bound for the $\go$-length of the $\gl$-primary component of the reduction $\tilde{A}(\F_q)$ of $A$ modulo $q$. 

Parallel to the proof of  Lemma  \ref{pretrivial}, let $\cZ$ be a counterexample of minimal length, $G$  the $\ell$-group $\Gal(\Q(Z)/\Q)$, $I_G$   the augmentation ideal of $\go_\gl[G]$ and  $r\ge 2$  the  least integer such that $I_G^r Z=0$. There is a prime $p = p_v$ ramified in $I_G^{r-2}Z$ and an element $z_2$ in $I_G^{r-2}Z$, such that $z_1=(\sigma_v-1)z_2\ne 0$ is killed by $\gl$.  Let $\cZ_1\subseteq \cZ$ correspond to the trivial $\go_\gl[G]$-module $Z_1=\langle z_1\rangle$ and let $\overline{\cZ} = \cZ/\cZ_1$.  Then $\len \overline{Z} - \gf(\overline{Z}) \ge \len Z - \gf(Z)$ and $\overline{\cZ}$ is a smaller counterexample.

The stronger bound  uses Faltings's theorem.  Let $\{\cZ_i\}$ be an increasing filtration of subschemes of $\cZ$.  Then  the abelian varieties $A_i=A/\cZ_i$ are non-isomorphic, since the kernel of an isomorphism $A_i \stackrel{\sim}{\to} A_j$ would equal $A_i[\gl^r]$ and so is not \'etale.  \qed

\begin{theo} \label{trivial}
Let $A_{/\Q}$ be an $\go$-type semistable abelian variety, good at $\ell$, and $\gl \, \vert \, \ell$ in $\go$.  Then
\begin{equation*}
  \epsilon_0(A[\gl]) \le \Omega(N_A^u(\gl)) + \Omega_{\ell}(N_A^u(\gl)) 
            +  \sum_{E{\rm \, in\, }\mathfrak{S}_\gl(A)} \delta_A(E).
\end{equation*}
\end{theo}
\proof
Put $\gS =\gS_\gl(A)$ and $\F = \F_\gl$.   Let $\cF$ be a nugget filtration of $A[\gl^n],$  with $\gr \cF = [\cM,\cV_1,\dots,\cV_m,\cZ]$, and each $\cV_i$ a nugget.  Set $\eta_i = 1$ if  $\cV_i$ is a prosaic nugget with a core of prime conductor $p \equiv \pm 1 \, \bmod{\tilde{\ell}}$, and $\eta_i = 0$ otherwise.  If $\cV_i$ is not prosaic, $V_i$ has a unique exceptional $\F[G_\Q]$-module $E_i$ as constituent.
Take the $\go$-length of $A[\gl^n]$  and apply Lemmas \ref{pretrivial} and \ref{pretrivialZ} to obtain
\begin{eqnarray*}
n \dim_\F A[\gl]  &=&  \len_\go Z + \len_\go M + \sum_{i}\dim_\F V_i \\
    &\le& 2s_0 + \gf(Z) + \gf(M)
          + \sum_{V_i \, {\rm prosaic}} (\gf(V_i)+\eta_{ i}) \\
    && \qquad  + \sum_{V_i \, {\rm not \, prosaic} }
           (\gf(V_i) + \dim_\F E_i -\gf(E_i)+\delta_A(E_i)) \\
                \end{eqnarray*}
\begin{eqnarray*}
    &\le&  2s_0 + \gf(Z) + \gf(M) + \sum_{i}\gf(V_i)  +
                \sum_{V_i\ {\rm prosaic}} \eta_i \\
     && \qquad + n \left[\sum_{E\, \text{in}\, \gS} (\dim_\F E-\gf(E) +
                 \delta_A(E)) \right],
\end{eqnarray*}
since  any $E$ appears $n$ times as often in $A[\gl^n]$ as in $A[\gl]$.
By Lemma  \ref{super1}(ii) and the bound (\ref{lr}) on the conductor of $A[\gl^n]$, we have
\begin{eqnarray*}
\gf(Z)+\gf(M) + \sum_{i} \gf(V_i) &\le& \gf(A[\gl^n]) \;\; \le \;\; 
n \,  \Omega(N_A^0) \\
   &\le&  n \left[ \Omega(N_A^u(\gl)) + \sum_{E\in \gS}\gf(E) \right].
\end{eqnarray*}
Clearly $\epsilon_0(A[\gl]) = \dim_{\F} A[\gl] - \sum_{E \in \gS} \dim_\F E$ and $\sum_{\rm prosaic} \eta_i \le n \, \Omega_{\ell}(N_A^u(\gl))$.
Substitute, divide by $n$ and let $n$ go to infinity to finish.  \qed

\begin{cor}   \label{nilpgenbd}
If $N_A$ is the  conductor of $A$ and $\Q(A[\gl])$ is an $\ell$-extension of $\Q(\Mu_\ell)$, then $\gS_\gl(A)$ is empty and 
$2 \, \dim A \le \Omega(N_A) + \Omega_\ell(N_A).$
\end{cor}

\begin{cor} \label{trivialplus}
Assume that $A$ is good outside $S$ and each exceptional constituent $E$ in $\gS_\gl(A)$ is $(S \backslash T_E)$-transparent, as in {\rm Definition} {\rm \ref{S-transp}}.
\begin{enumerate}[{\rm i)}] \item We have $\epsilon_0(A[\gl])=0$ under  one of the following:
\begin{enumerate}[{\rm a.}] 
    \item $\ell\ge 5$ and no prime dividing $N_A^u(\gl)$ satisfies $p \equiv \pm 1 \, (\ell);$ 
    \item $\ell = 3$ and  $N_A^u(\gl) =q^a$ with $q \not\equiv \pm 1 \, (9);$      
    \item   $\ell = 2$ and $N_A^u(\gl) =q^a$ with $q^* \equiv 5 \, (8).$  
\end{enumerate}
\item We have $\epsilon_0(A[\gl])\le 2a$ under one of the following: 
\begin{enumerate}[{\rm a.}] 
       \item  $\ell = 3$ and   $N_A^u(\gl) =p^a q^b$ with either $q \bmod{9}$ in 
             $\{ 2,5\},$ or with $q \bmod{9}$ in $\{ 4,7\}$ and 
             $p^{\frac{q-1}{3}}\not\equiv 1\, (q);$
       \item $\ell = 2$ and $N_A^u(\gl) =p^a q^b$ with $p^*\equiv 1 \, (8),$   $q^*\equiv 5 \, (8)$ and some Hilbert symbol $(p^*,q^*)_\pi=-1.$
\end{enumerate}
\end{enumerate} 
\end{cor}

\proof 
By $(S \backslash T)$-transparency, each non-prosaic nugget is itself an exceptional.  Then, as in the proof of the last theorem,
$$
n \,\epsilon_0(A[\gl]) \le 2 s_0 + \sum_{V_i \; {\rm prosaic}} \dim V_i.
$$
By Corollary \ref{core}, there is  no prosaic nugget in (i).   In (ii), all are two-dimensional, of conductor $p$ by Lemma \ref{pnil}, so there are at most $na$ prosaic nuggets.  \qed

\section{Mirages}\label{PQR} 
\numberwithin{equation}{subsection}

\subsection{Introducing mirages}
Let $A_{/\Q}$ be a semistable abelian variety of $\go$-type with good reduction at $\ell.$ Let $\gl$ be a prime of $\go$ above $\ell$ and $\F=\F_\gl.$  Recall that $B$ is an object of $\gI^\gl_A$ if there is an $\go$-isogeny from $A$ to $B$ with $\gl$-primary kernel.  In this section, all   isogenies  are $\go$-linear, with $\gl$-primary kernels.

\begin{Def} \label{mirage}
A {\em mirage} $\gC$ associates to each $B$ in $\gI^\gl_A$ a set of $\F$-module subschemes of $B[\gl]$ such that $\varphi(\gC(B_1))\subseteq\gC(B_2)$ for each isogeny $\varphi\!:\,B_1\to B_2.$   Call $B$  {\em obstructed} if $\gC(B)=\{0\}$ and  $\gC$ {\em unobstructed} if no $B$ in $\gI_A^\gl$ is  obstructed.
\end{Def}

\begin{prop} \label{Faltings} 
If $\gC$ is  unobstructed on $\gI_A^\gl$, then there is a  $B$ in $\gI_A^\gl$ and a filtration
$
0 \subset \cW_1 \subset \dots \subset \cW_s = B[\gl^r],
$
with $\cW_{i+1}/\cW_i$ in $\gC(B/\cW_i)$ for all $i $.
\end{prop}

\proof  
Set $A_0 = A$ and construct inductively the abelian variety $A_{n}=A_{n-1}/\kappa_{n}$ with $\kappa_{n}$ chosen in $\gC(A_{n-1}).$ Write $\cK_n$ for the kernel of the induced isogeny from $A$ to $A_n$.  By Faltings \cite{Falt}, we may find an isomorphic pair $B=A_m$ and $B'=A_n$ with $m < n$.  This produces an endomorphism $\alpha$ of $B,$  whose kernel $\cW = \cK_n/\cK_m$ admits a filtration as above.  Since $\alpha$ is in $\End B = \go$ and $\cW$ is killed by a power of $\gl$, we have $\cW=B[\alpha]=B[\gl^r]$, with $\alpha \go=\gl^r.$    \qed

\medskip

Here is a particularly simple illustration of the use of mirages.

\begin{cor}\label{2Es} If  $E, E_1,..,E_r$ are the distinct irreducible Galois constituents of $A[\gl]$, then a subquotient of $A[\gl^{\infty}]$ is a non-split extension of $E_i$ by $E$ for some $i \ge 1$. 
\end{cor}

\proof For any $B$ in $\gI_A^\gl,$ consider the mirage $\gC(B)$ consisting of the $\F$-module subschemes of $B[\gl]$ whose Galois submodules have {\em no} constituent isomorphic to $E.$  Thanks to Propositions  \ref{JH1} and \ref{Faltings}, we may assume that $A$ is obstructed.  Choose a subscheme $\cV\subseteq A[\gl]$ with $\gr V=[E \dots E E_i]$ having the smallest number of  $E$'s. If the quotient $W$ corresponding to the last two terms were split, $\cW$  would contain  an $\F$-submodule scheme $\cE'$ with $E'\simeq E_i$. This violates the minimality  and thereby the obstruction.  \qed

\medskip

Let $C$ be a covariant functor from $\gI^\gl_B$ to the category of $\gol$-modules, such that $C(A)$ is a pure $\gol$-submodule of $\T_\gl(A)$ for all $A$ in $\gI^\gl_B$.   Let $\varphi_* = C(\varphi)$ be the map induced by an $\go$-isogeny $\varphi:  A \to A^\prime.$    Denote the image of $C(A)$ in $A[\gl^n]$ by
\begin{equation} \label{funcmir}
C^{(n)}(A)=(C(A)+\gl^n\T_\gl(A))/\gl^n\T_\gl(A)
\end{equation}
and set $\overline{C}(A) = C^{(1)}(A) \subseteq A[\gl]$.  We create a mirage by letting $\gC(A)$ be the set of all simple $\F$-module subschemes of $A[\gl]$ whose Galois module is contained in $\overline{C}(A)$.  We say that $C$ is obstructed if  $\gC$ is obstructed. 

\begin{lem} \label{FaltC}
If $C$ is unobstructed, then $C(A) = \T_\gl(A)$ for all $A\in \gI^\gl_B$.
\end{lem}

\proof 
We first show that if
$
A_1 \stackrel{\varphi}{\longrightarrow} A_2
\stackrel{\psi} {\longrightarrow} A_3 
$ 
is a chain of $\go$-isogenies such that $\ker \varphi_* \subseteq C^{(n_1)}(A_1)$ and $\ker \psi_* \subseteq C^{(n_2)}(A_2)$, then $\ker (\psi\varphi)_* \subseteq C^{(n_1+n_2)}(A_1)$.  The kernel of $\varphi$ is annihilated by $\gl^k$ for some $k \le n_1$.  There is a quasi-inverse $\go$-isogeny  $\varphi^\prime:  A_2 \to A_1$, such that the induced maps $(\varphi \varphi^\prime)_*$ and $(\varphi^\prime \varphi)_*$ are multiplication by $\gl^k$ on $\T_\gl(A_2)$ and on $\T_\gl(A_1)$, respectively.  Hence, 
$$
C^{(n_2)}(A_2) = \gl^k C^{(n_2 + k)}(A_2) = \varphi_*\varphi'_*(C^{(n_2 + k)}(A_2)) \subseteq \varphi_*(C^{(n_2 + k)}(A_1)).
$$ 
If $x$ lies in $\ker (\psi \varphi)_*$, then $\varphi_*(x)$ is in $\ker \psi_* \subseteq C^{(n_2)}(A_2)$, so we can find $y$ in $C^{(n_2+k)}(A_1)$ satisfying $\varphi_*(x) = \varphi_*(y)$.  Hence $x$ is in $y + \ker \varphi_*  \subseteq C^{(n)}(A_1)$ for all $n \ge \max \{n_2+k,n_1\}$.

Next, as in the proof of Proposition \ref{Faltings}, we may find an endomorphism of some $A$ in $\gI^\gl_B$ whose kernel $\cW = A[\gl^r]$  is the kernel of the composition of a suitably long chain of isogenies as above.  Hence $\cW \subseteq C^{(n)}(A)$ for $n$ sufficiently large.  Thus $\rk_{\gol} C(A) = \rk_{\gol} \T_\gl(A).$  The ranks on both sides are $\go$-isogeny invariants.  Therefore, by purity, $C(A') = \T_\gl(A')$ for all $A'$ in $\gI^\gl_B$.  \qed

\medskip 

The toric space $M_t(A,v,\gl)$ and finite space $M_f(A,v,\gl)$ described in \S \ref{Prelims} will be used to build mirages.  Let  $\cP$ be a  set of places of $\overline{\Q}$ with exactly one $v$ over each  bad prime $p$ of $A.$  For any subset $\cP'$ of $\cP$, let 
$$
M_t(A,\cP',\gl) = \langle M_t(A,v,\gl) \, | \, v \in \cP' \rangle^{sat},
$$
where the {\em saturation} of an $\gol$-submodule $X$ of $\T_\gl(A)$ is the pure submodule
$$
X^{sat} = ( k_\gl \otimes X) \cap \T_\gl(A),
$$ 
with $k_\gl$  the field of fractions of $\go_\gl$.  For $\go$-isogenies $\varphi\!:\, A \to A'$, we have the desired functoriality $\varphi_*(M_t(A,v,\gl)) \subseteq M_t(A',v,\gl)$.  If $C(A)$ contains $M_t(A,\cP',\gl)$, then the same holds for all $B$ in $\cI^\gl_A$ by purity.  In view of (\ref{rank}), we have
\begin{equation} \label{sizeC} 
 \max_{\; v \in \cP'} \; \tau_{p_v} \le  \rk_{\gol} M_t(A,\cP',\gl) \le  \sum_{v \in \cP'} \tau_{p_v}.
\end{equation} 
The  following  lemma can provide a better lower bound when $\cP$ is suitably  chosen.

\begin{lem} \label{sizeC3} 
Let $X$ be a proper pure $\gol$-submodule of $\T_\gl(A)$ and  $p$  a prime of bad reduction for $A.$   Then we can find a place $v$ above $p$ in $L_\infty=\Q(A[\gl^{\infty}])$ such that $X+M_t(A,v,\gl)$ contains $X$ properly. 
\end{lem} 

\proof 
Let $G=\Gal(L_\infty/\Q)$ and pick some place $w$ over $p$. If the claim is false, we have  $\sigma(M_t(A,w,\gl))=M_t(A,\sigma(w),\gl)\subseteq X$ for all $\sigma$ in $G,$ so $X$ contains the $\gol[G]$-submodule $Y$ of $\T_\gl(A)$ generated by $M_t(A,w,\gl).$ But Tate's conjecture, proved by Faltings, asserts that $\End_{\Z_{\ell}[G]}  (\T_{\ell}(A))=\End A \otimes_{\Z} \Z_{\ell}.$ Thus $\End_{\gol[G]} (\T_{\gl}(A))=\gol$ and the semisimplicity of $\T_\gl(A)$ implies that $X=Y=\T_\gl(A).$  \qed

\subsection{Mirages in the prosaic case} \label{UnipMir}  
Let $G$ be a 2-group, $\F$ a finite field of characteristic 2 and $W$ an $\F[G]$-module.  For any subgroup $H$ of $G$, let $\ga_H$ be the augmentation ideal in $\F[H]$, with $\ga = \ga_G$.  If $H = \langle g_j \, | \, 1 \le j \le n \rangle$, the identity 
$$ 
g_1 g_2 - 1 = (g_1-1) + (g_2 - 1) + (g_1-1)(g_2-1)
$$
shows that
$  
\ga_H =  \langle g_j-1 \, | \, 1 \le j \le n \rangle. 
$ 
For $k \ge 0$, we consider the filtration 
\begin{equation} \label{filt} 
W_k = \{ x \in W \, | \, \ga^k x = 0 \} = \{ x \in W \, | \, \ga x \in W_{k-1} \}.
\end{equation} 
Then $0 = W_0 \subset \dots \subset W_j \subset \dots \subset W_m = W$ for some $m \ge 0$, with proper inclusions along the way.  Within the appropriate ranges of $k$ and $k'$, we have 
\begin{equation} \label{expW} 
\ga^k \, W_{k'+k} \subseteq W_{k'}.
\end{equation} 
Thus $G$ acts trivially on $W_{k+1}/W_k$, has exponent two on $W_{k+2}/W_k$ and exponent dividing four on $W_{k+4}/W_k$.  In particular, $W_1 = W^G$. 

\begin{lem} \label{Wbd}
Let $H = \{ h \in G \, | \, (h-1)(W_{k+2}) \subseteq W_k \}$.  Then $\overline{G} = G/H$ is elementary abelian, say of rank $r$, and 
$
\dim W_{k+2}/W_{k+1} \le  r  \, \dim W_{k+1}/W_k.
$
\end{lem}

\proof
We have an injective $\F$-linear map $\psi:W_{k+2}/W_{k+1} \to \Hom_{\F_2}(\overline{G},W_{k+1}/W_k)$ induced by $\psi(x)(g) = (g-1)(x)$ for $x$ in $W_{k+2}$ and $g$ in $G$.    \qed 

\begin{lem} \label{abG} Assume that the maximal quotient $\gG$ of $G$ acting faithfully on $W_3$ is abelian.   If either {\rm (i)} $W_2^{\langle g \rangle} = W_1$ for some involution $g$ in $\gG$, or {\rm (ii)} $\gG$ is elementary abelian and $\dim W_2/W_1 = 1,$ then $W_3 = W_2$.
\end{lem}

\proof
(i) If $x$ is in $W_3-W_2$, we can find $h$ in $\gG$ such that $y = (h-1)(x)$ is not in $W_1$ and so $z = (g-1)(y) \ne 0$.  But $(g-1)(x)$ is in $W_2^{\langle g \rangle} = W_1$ and so is  fixed by $h$.  Hence,
$
 0 = (h-1)(g-1)(x) = (g-1)(h-1)(x) = (g-1)(y) = z.  
$

(ii) For some $g$ in $\gG,$  $g-1$ has rank one on $W_2$ and  $\dim W_2^{\langle g \rangle} = \dim W_1.$ Now (i) applies because $W_1 \subseteq  W_2^{\langle g \rangle}.$ \qed
 
\begin{lem} \label{WProps}  If   $g^2=1$ on $W_{3}$ and  $W_2^{\langle g\rangle}=W_1,$ then $g$ acts trivially on $W_{3}/W_{1}.$    
\end{lem}

\proof 
We have $(g-1)(W_3) \subseteq W_2^{\langle g \rangle} = W_1$, so $g$ is trivial on $W_3/W_1$.  \qed

\begin{lem} \label{4steps}
If $g = \left[
    \begin{smallmatrix}1&x&a&b\\0&1&c&d\\0&0&1&y\\0&0&0&1\end{smallmatrix}
                               \right]$ 
has order two, then $cx=0, cy=0$ and $dx+ay=0$.
\end{lem}

\noindent For the rest of this subsection, we assume:

\medskip

\setlength{\hangindent}{30 pt}\noindent \fbox{\bf M1}
$W$ is a  Galois submodule of $A[\gl]$ with $G = \Gal(\Q(W)/\Q)$ a 2-group generated by involutions and $W_k$ is given by (\ref{filt}).

\setlength{\hangindent}{0 pt}

\medskip

\noindent Let $\chi_d$ denote the quadratic character of $\Q(\sqrt{d})$.  

\begin{lem} \label{step1} 
Suppose $C(A)$ contains $M_t(A,v,\gl)$ and $A$ is obstructed for $C$.  Then $\overline{M}_t(A,v,\gl) \cap W_1 = 0$ and $p_v$ does not ramify in $\Q(W_2)$. Assume further  that $\Q(W_2) = \Q(\sqrt{d_1}) $ is a quadratic field.
\begin{enumerate}[{\rm i)}] 
   \item Then $W_2$ contains all submodules $U$ of $W$ such that $\Q(U) \subseteq \Q(\sqrt{d_1})$.
   \item If $\chi_{d_1}(p_v) = -1$, then $\overline{M}_t(A,v,\gl) \cap W_2 = 0$ and $p_v$ does not ramify in $\Q(W_3)$.  
   \item If 2 ramifies in $\Q(\sqrt{d_1})$, then  $\cW_2/\cW_1\simeq \Zl^r$ as  group schemes.
\end{enumerate}
\end{lem}
\proof 
By definition, $G$ is trivial on $W_1$.  Thus $X = \overline{M}_t(A,v, \gl) \cap W_1$ is a Galois module and then $X = 0$ because $A$ is obstructed.    Since  $(\sigma_v-1)(W_2)$ is contained in $X$, we see that $p_v$ does not ramify in $\Q(W_2)$.
\begin{enumerate}[{\rm (i)}]
\item  $G$ acts on $\Q(\sqrt{d_1})$ via $\langle g\rangle$ for some involution $g$ in $G$.  Hence $G$ acts trivially on $(g-1)(U)$.  We deduce that $(g-1)(U) \subseteq W_1$ and so $U \subseteq W_2$. 
\item Any Frobenius $\Phi_v$ in $G$ restricts to a generator of $\overline{G} = \Gal(\Q(\sqrt{d_1})/\Q)$.  But $\overline{M}_t(A,v,\gl)$ is a $\cD_v$-module, so $Y = \overline{M}_t(A,v,\gl) \cap W_2$ is a $G_\Q$-module and then $Y = 0$ because $A$ is obstructed.   Since  $(\sigma_v-1)(W_3)$ is contained in  $Y$, we see that $p_v$ does not ramify in $\Q(W_3)$. 
\item the involution $\sigma_\lambda$ (see Remark  \ref{slam}) restricts to a generator of $\overline{G}$ and $\sigma_\lambda$ acts trivially on the multiplicative component $W_2^m$ at $\lambda$.  Hence $W_2^m$ is contained in  $W_1$.  It follows that $\cW_2/\cW_1$ is \'etale at 2.  Since $G_\Q$ acts trivially, $\cW_2/\cW_1$ is isomorphic to a direct sum of copies of $\cZ_\gl$ globally.  \qed 
\end{enumerate}  

\medskip

Let $\cP^u = \{ v \in \cP \; | \; p_v \in \Pi_A^u\}$,  where $\Pi_A^u$ is the set of prime divisors  of the prosaic $\gl$-conductor  $N_A^u(\gl)$ as in Definition \ref{Nss}.  Note that $N_W$ divides $N_A^u(\gl)$.  

\begin{prop} \label{step2} 
If $C(A)\supseteq M_t(A,\cP^u,\gl)$, $A$ is obstructed for $C$ and $W_1 \subsetneq W_2,$ then $\Q(W_2) = \Q(i)$.   Moreover:
\begin{enumerate}[{\rm i)}] 
 \item the odd primes ramified in $\Q(W_3)$ are $1 \bmod{4}$;
 \item $K=\Q(W_3/W_1)$ is a totally real elementary 2-extension unramified at 2;
 \item $\Q(W_3)/K$ is unramified at odd places.
\end{enumerate}
\end{prop}

\proof
Lemma  \ref{step1} implies $\Q(W_2)$ is unramified at odd places, so $\Q(W_2) = \Q(i)$.  By Lemma  \ref{WProps}, we find that $g=\sigma_\infty$ and $g=\sigma_\lambda$ act trivially on $W_3/W_1$.  Hence $K$ is totally real and unramified over 2.  If $p_v$ ramifies in $\Q(W_3)$, then $p_v \equiv 1 \, (4)$ by Lemma  \ref{step1}(ii).   Furthermore, $p_v$ already ramifies in $K$.  Otherwise, $\sigma_v$ acts trivially on $W_3/W_1$,  so $(\sigma_v-1)(W_3) \subseteq \overline{M}_t(A,v,\gl)\cap W_1=0$ by Lemma  \ref{step1} making $\sigma_v$  trivial on  $\Q(W_3)$.  The necessarily odd primes that  ramify in $K/\Q$, cannot ramify further in $\Q(W_3)/K$ by Lemma  \ref{ram}.  \qed

\begin{cor} \label{step3}
Assume that $K$ is a quadratic field $\Q(\sqrt{d_2})$.  Then $\Q(W_3)$ is in $D_4(-1,d_2).$   Let $n$ be maximal such that $W_{n-1} \ne W_n$ and $\Gal(\Q(W_n)/\Q)$ is generated by two elements.  If $q_w \equiv 3 \, (4)$ and  $\chi_{d_2}(q_w) = -1$ for some $w$ in $\cP^u$, then $\overline{M}_t(A,w,\gl) \cap W_n = 0, \, W_n \subsetneq W_{n+1}$ and $q_w$ does not ramify in $W_{n+1}$.
\end{cor}

\proof
Fix $v$ in $\cP$ such that $p_v$ divides $d_2$.  The group  $\Gal(\Q(W_3)/\Q)$ is generated by $\sigma_\infty$ and involutions $\sigma_{v'}$ with $v'$ in $\cP^u$.   If $\sigma_{v'}$ is not trivial on $W_3$, we show that $\sigma_{v'} = \sigma_v$ on $W_3$.  Indeed, $\sigma_v$ and $\sigma_{v'}$ agree on $K$.  For $x$ in $W_3$,  it follows that  $y= \sigma_v(x) - \sigma_{v'}(x)$ becomes trivial in $W_3/W_1$, so $y$ is in $W_1.$  Now 
$$ 
y = (\sigma_v - 1)(x) - (\sigma_{v'}-1)(x) \in \overline{C}(A) \cap W_1 = 0
$$ 
because $A$ is obstructed.  Hence $\Gal(\Q(W_3)/\Q) = \langle \sigma_\infty, \sigma_v \rangle$.  From matrix representations for $\sigma_\infty$ and $\sigma_v$ with respect to the filtration on $W_3$, one easily sees that $\Q(W_3)$ is in $D_4(-1,d_2)$.  

By Burnside's theorem, $\Gal(\Q(W_n)/\Q) = \langle \sigma_\infty, \sigma_v \rangle$. Thus $\Gal(\Q(W_n)/\Q)$ is dihedral and $\tau = \sigma_\infty \sigma_v$ generates the cyclic subgroup of index 2. The fixed field of $\tau$ is $\Q(\sqrt{-d_2})$. 

Suppose the hypotheses on $q_w$ hold.  Then the restriction of a Frobenius $\Phi_w$ to $\Q(W_n)$ generates the same subgroup of $\Gal(\Q(W_n)/\Q)$ as $\tau$.  Since $M_t(A,w,\gl)$ is a $\cD_w$-module, $\tau$ preserves $Y = \overline{M}_t(A,w,\gl) \cap W_n$.  If $Y \ne 0$, then $\tau$ has a non-zero fixed point $y$ in $Y.$  It follows that $\sigma_\infty(y) = \sigma_v(y)$ and so 
$$
z = (\sigma_\infty-1)(y) = (\sigma_v-1)(y)
$$ 
is fixed by $\sigma_\infty$ and $\sigma_v$.  Hence $z$ is a rational point in $\overline{M}_t(A,w,\gl)$.  But then $z=0$ because $A$ is obstructed.  From this, it follows that both $\sigma_v$ and $\sigma_\infty$ fix $y$.  Hence $y$ is a rational point in $Y$.  Since $A$ is obstructed, we conclude that $Y=0$, so $W_n \subsetneq W_{n+1}$ and $q_w$ is unramified in $W_n$.  \qed

\begin{theo}\label{pqrthm}
Let $A$ be a  semistable $\go$-type abelian variety with odd  conductor $N$.  Then $2 \dim A \le \Omega(N)$ if  $\Gal(\Q(A[\gl])/\Q)$ is a $2$-group for some $\gl \vert 2$ in $\go$ and either:
{\rm i)} all prime factors of $N$ are $3 \bmod{4},$ or

\hspace{14pt} {\rm ii)} at least two primes divide $N$, one $p \equiv 1 \,(4)$ and  for all other $q \, | \, N,$

\hspace{27pt}  we have $q\equiv 3 \,(4)$ and $\chi_p(q) = -1$. 
\end{theo} 

\proof 
Let $C(B) = M_t(B,\cP,\gl)$ for all $B$ in $\gI^\gl_A$.  If $A$ is obstructed for the associated mirage, Corollary  \ref{FaltC} implies that $C(A) = \T_\gl(A)$ and so (\ref{sizeC}) gives
$$
2 \dim A = d \, \rk_{\gol} \T_\gl(A) = d \, \rk_{\gol} C(A) \le   \Omega(N_A).
$$
Suppose $A$ is obstructed and consider the filtration (\ref{filt}) of $W=A[\gl]$.  Since $\cP$ is not empty, Lemma  \ref{step1} shows that $W_2\supsetneq W_1,$ so $\Q(W_2) = \Q(i)$ by Proposition \ref{step2}(i).  Since at least one prime $q \equiv 3 \, (4)$ divides $N_A$, Lemma  \ref{step1}(ii) shows that $W_3\supsetneq W_2$ and the odd primes ramifying in $\Q(W_3)$ are $1 \bmod{4}$.  In case (i), we now have $\Q(W_3) = \Q(i)$.  But then $W_3 = W_2$ by Lemma  \ref{step1}(i).

Assume (ii) holds and $\Q(W_3)\supsetneq \Q(i).$ 
By Proposition \ref{step2},  $K = \Q(W_3/W_1) = \Q(\sqrt{p}).$ With $n$ as in Corollary \ref{step3}, $\Q(W_3)$ is in $D_4(-1,p)$ and $W\supsetneq W_n.$ By Burnside, there is a quadratic field $\Q(\sqrt{d_3})$ in $\Q(W_{n+1})$ but not in $\Q(i,\sqrt{p}).$ Thus  some $q$ ramifies in $\Q(W_{n+1})$ and contradicts Corollary \ref{step3}.  \qed

\medskip

\noindent For the rest of this subsection, we assume: 

\medskip

\setlength{\hangindent}{30 pt}\noindent \fbox{\bf M2}
$A$ is $(\go,N)$-paramodular, $W = A[\gl]$, $L = \Q(W)$ and $G = \Gal(L/\Q)$ is a $2$-group.   In particular, $A$ has good ordinary reduction at 2.  Recall that $N_A^0 = N$ is the reduced conductor of $A$.

\setlength{\hangindent}{0 pt}\noindent
 
\begin{prop} \label{morepqr}
Assume $N_A^0 = pqr$ for primes $p,q,r$ with $p \equiv -q \equiv 5 \, (8)$ and $r \equiv 7 \, (8)$.  Then $\chi_p(r)=1$.  Moreover, $\chi_q(p) = 1$ or $\chi_q(r) = 1$.
\end{prop}

\proof
By Lemma  \ref{sizeC3}, we choose $\cP$ so that  $C(A) = \M_t(A,\cP,\gl)^{\rm sat}$ has   $\go_\gl$-rank three.  Suppose $A$ is obstructed for the associated mirage and let $W = A[\gl]$.  Since $\overline{C}(A) \cap W_1 = 0$, we have $\dim_\F W_1 =1$.

Proposition \ref{step2} and its Corollary show that $\Q(W_2) = \Q(i)$, $\Q(W_3/W_1) = \Q(\sqrt{p})$ and $\Q(W_3)$ is in $D_4(-1,p)$.  We have $\dim W_k = k$ for $k = 1,2,3$ by Lemma  \ref{Wbd} and so $W = W_4$.  Because 2 ramifies in $\Q(W_2)$ and is inert in $\Q(W_3/W_1)$, we find that $\gr W_3 = [\Mu_\gl \cZ_\gl \cZ_\gl]$.  Hence $\gr W = [\Mu_\gl  \cZ_\gl \cZ_\gl  \Mu_\gl]$, forcing 2 to split in $\Q(W_4/W_2)$.  But the conductor of $W_4/W_2$ divides $qr$ and so $\Q(W_4/W_2) = \Q(\sqrt{-r})$.  By Corollary \ref{step3}, we have $\chi_p(r) = 1$.

Let $\Phi_q = {\rm Frob}_w$ at the place $w$ in $\cP$ over $q$.  Suppose, contrary to our claim, that $\chi_q(p)  = \chi_q(r) = -1$.  Then $\Phi_w$ admits a matrix representation as in (\ref{4steps}) with $c,x,y$ all non-zero, so $\ker (\Phi_w - 1) = W_1$.  Since the $\Phi_w$-module $\overline{M}_t(A,w,\gl)$ is 1-dimensional over $\F$, $W_1 = \overline{M}_t(A,w,\gl)$ and so $A$ is not  obstructed.   \qed

\begin{prop} \label{firstpqa}
If $q\! \equiv\!  3 \, (4)$ and $N_A^0 \! = \!  pq^a,$  then $a\!  = \! 2$, $p \equiv 1 \, (4)$ and $\chi_p(q)\!  =\!  1$.
\end{prop}

\proof
If $v$ is a place over $p$, then $M_f(A,v,\gl) = \T_\gl(A)^{\cI_v}$ is a pure $\go_\gl$-submodule of $\T_\gl(A)$ of rank 3.  Suppose $A$ is obstructed for the mirage associated to $C(A) = M_f(A,v,\gl)$.  Then $\overline{C}(A) \cap W_1 = 0$,  so $\dim_{\F} W_1 = 1;$ thus $W = W_1 \oplus  \overline{M}_f(A,v,\gl)$.  Now $\cI_v$ acts trivially on $W$, so $L/\Q$ is unramified at $p$.  It follows that the maximal elementary 2-extension of $\Q$ inside $L$ is contained in $\Q(i,\sqrt{q})$.  Hence $G = \langle \sigma_\infty, \sigma_w \rangle$, where $w$ is a place over $q.$ 

Since the Hilbert symbol $(-1,q)_q = -1$, there is no $D_4(-1,q)$ field and so $G$ is abelian.  Lemmas \ref{Wbd} and \ref{abG} now imply that $\dim W_2 = 3$,  and $\Q(W_2) = \Q(i,\sqrt{q})$.  In particular, $\sigma_w$ is not trivial on $W_2$ and so $(\sigma_w-1)(W_2)=W_1$.

If $a=1$, then $\dim (\sigma_w-1)(W) = 1$, so $(\sigma_w-1)(W) = W_1$ and we find that $(\sigma_\infty-1)(\sigma_w-1)(W) = 0$.  Because $\sigma_w$ and $\sigma_\infty$ are commuting involutions, it follows that $\ga_G^2 \, W = 0$.  But then $W = W_2$, a contradiction.  Hence  $a = 2$.  Finally, by Theorem \ref{pqrthm}, we have $p \equiv 1 \, (4)$ and $\chi_p(q) = 1$. \qed

\begin{prop} \label{secondpq}
If $q \equiv 5 \, (8)$ and $N_A^0 = pq^2$, then $p^* \equiv 1 \, (8)$ and $\chi_p(q) = 1.$
\end{prop}

\proof
We have $p^* \equiv 1 \, (8)$ by Theorem \ref{trivial}.  Fix a place $\lambda$ over 2 to define the multiplicative component $\T_\gl(A)^m$, which has  $\gol$-rank 2 because $A$ is ordinary at 2.   By Lemma  \ref{sizeC3}, we can choose $v$ over $p$ to guarantee that the $\gol$-rank of 
$$
C(A) = \{M_t(A,v,\gl) + \T_\gl(A)^m\}^{sat}
$$
is 3.  Assume that $A$ is obstructed for $C$ and let $W = A[\gl]$.  Then $\overline{C}(A) \cap W_1 = 0$, so $\dim W_1 = 1$.  Moreover, the $\F$-module scheme associated to $W_1$ is $\cW_1 \simeq \cZ_\gl$ and $\Q(W_2)$ is unramified at $p$.  Choose $\sigma_\lambda$ as in Remark  \ref{slam}.  Then  $(\sigma_\lambda-1)(W_2)$ is contained in $W_1^m = 0$. Thus $\sigma_\lambda$ fixes $W_2$ pointwise 
and  $\Q(W_2)$ is unramified at 2.

It follows that $\Q(W_2) = \Q(\sqrt{q})$ and so $\dim W_2 = 2$ by Lemma  \ref{Wbd}.  Moreover, $\gr \cW_2 = [\cZ_\gl \, \cZ_\gl]$ because 2 is inert in $\Q(W_2)$. Let $V$ be any  Galois submodule of $W$ containing $W_2$ with $\dim_\F V = 3$.  Then $\Gal(\Q(V)/\Q) \simeq D_4$ and $\gr \cV = [\cZ_\gl \, \cZ_\gl \,  \Mu_\gl]$.  Since 2 splits in $\Q(V/W_1)/\Q$, we have $\Q(V/W_1) = \Q(\sqrt{p^*})$, whence $\chi_p(q) = 1$ by Lemma  \ref{D41}.    \qed

\begin{prop} \label{pmmir1}
Let $N_A^0 = pqr^a$ with  $p,q,r$  prime and $q^* \equiv r^* \equiv 5 \, (8).$ Assume $K = \Q(\sqrt{p^*},\sqrt{q^*},\sqrt{r^*})$ has no quadratic extension   unramified outside $\infty$ and split over 2.
\begin{enumerate}[{\rm i)}]
\item If  $p^* \equiv 1 \, (8)$ and $1\le a\le 2,$    then $(p,q,r) \equiv (1,5,5) \bmod{8}$.
\item If $p^* \equiv 5 \, (8)$ and $a=2,$ then $p \equiv q \equiv r \equiv 5 \, (8)$.
\end{enumerate}
\end{prop}

\proof We refer to  the filtration (\ref{filt}) of $\cW = A[\gl].$  Let $A$ be obstructed for the mirage
$
C(A) = \{M_t(A,v,\gl) + \T_\gl(A)^m\}^{sat}
$
of  $\go_\gl$-rank three, as in the proof of Proposition \ref{secondpq}, so $\cW_1 \simeq \cZ_\gl$ and $\Q(W_2) \subseteq \Q(\sqrt{q^*},\sqrt{r^*})$.  By Lemma  \ref{Wbd}, $\dim W_2 \le 3$ with equality only if $[\Q(W_2):\Q]=4$.   A nugget filtration of $\cW$ must have one of the following gradings:
$$
{\alpha :} \; \;  [\cZ_\gl\cZ_\gl\Mu_\gl\Mu_\gl ],  \quad \quad {\beta\rm :}  \;\;  [\cZ_\gl\Mu_\gl\Mu_\gl\cZ_\gl ], \quad \quad {\gamma\rm :} \;\; [\cZ_\gl\Mu_\gl\cZ_\gl\Mu_\gl].
$$
If $0 \subsetneq \cV_1 \subsetneq \cV_2 \subsetneq \cV_3 \subsetneq \cW$ is the corresponding flag, then  $\cV_1 = \cW_1$ and  $\cV_2 \subseteq \cW_2$. 

Let $X\subseteq W$ be a Galois submodule with $\Q(X)/\Q$ abelian. If $p$ ramifies in $\Q(X)$, then $(\sigma_v-1)(X)$ is a Galois submodule of $W$, violating the obstruction.

Suppose $\alpha$ or $\beta$ holds.  Then $\cV_3$ is a nugget and by Corollary \ref{D42} (or its dual), we find that $K(V_3)/K$ is an elementary 2-extension, unramified at finite places and split over 2. Our assumption now implies that $K(V_3) = K,$ so $\Q(V_3)/\Q$ is  abelian.  But $p$ ramifies in $\Q(V_3)$ by Lemma  \ref{pretrivial}, a contradiction. 

Assume $\gamma$ holds. Then $\Q(V_2) = \Q(\sqrt{q^*r^*})$ because $p$ is unramified and 2 splits.   Since $\cW/\cV_2$ is a nugget,  we have  $\Q(W/V_2) = \Q(\sqrt{d_3})$, with $d_3 = p^*$ in case (i) and $d_3 = p^*r^*$ in case (ii).  Let $g$ in Lemma  \ref{4steps} generate the relevant inertia group, to conclude that $\Q(V_3/V_1)$ is unramified at $p,q,r$ and so $\Q(V_3/V_1)\subseteq \Q(i)$.  

If $\Q(V_3/V_1) =\Q$, then $V_3 = W_2$ and $\Q(W_2) = \Q(\sqrt{q^*},\sqrt{r^*})$.  Because $W^m$ and $W^{et}$ are unramified $\cD_\lambda$-modules,  2 is unramified in $\Q(W)$ and $f_\lambda(\Q(W)/\Q) = 2$.   But then $K(W) = K$ leads to a contradiction as above. 

If $\Q(V_3/V_1) =\Q(i),$ Lemma  \ref{D41} shows that $(-1,q^*r^*)_\pi = (-1,d_3)_\pi = 1$ for all $\pi$ in $\{p,q,r\}$.  Hence $p \equiv q \equiv r \equiv 1 \, (4)$ and the claim  ensues. \qed

\medskip

We sketch a more easily tested version of the previous proposition.

\begin{prop} \label{pmmir2} 
Let $N_A^0 = pqr^a$ with $p,q,r$  prime  and $q^* \equiv r^* \equiv 5 \, (8).$  
\begin{enumerate}[{\rm i)}]
\item If $p^* \equiv 1 \, (8)$, $1\le a\le 2$ and none of $D_4^{nr}(p^*,q^*)$, $D_4^{nr}(p^*,r^*)$, $D_4^{sp}(p^*,q^*r^*)$ exists, then $(p,q,r) \equiv (1,5,5) \bmod{8}$.
\item If  $p^* \equiv 5 \, (8),$ $a=2$  and  none of  $D_4^{nr}(p^*q^*,r^*)$, $D_4^{nr}(p^*r^*,q^*)$, $D_4^{nr}(q^*r^*,p^*)$ exists, then $p \equiv q \equiv r \equiv 5 \, (8)$.
\end{enumerate}
\end{prop}

\proof 
In case $\alpha$ or $\beta$ of the previous proof, Corollary \ref{D42} implies that the current conditions suffice. Case $\gamma$ leads to a quadratic extension $L/K$, unramified outside infinity and split over 2, such that $L/\Q$ is Galois, with group $D_4 \times C_2$.  This descends to a $D_4$-extension $M/\Q$, such that $M/k$ is cyclic of order 4, unramified outside infinity and split over 2, with $k = \Q(\sqrt{p^*q^*})$ or $k = \Q(\sqrt{p^*r^*})$ in (i) and $k = \Q(\sqrt{p^*q^*r^*})$ in (ii).  \qed

\subsection{Mirages with exceptionals}
In this subsection,  $A$ is a  semistable abelian variety of $\go$-type with good reduction at $\ell=2$, $\gl$ is a prime over 2 in $\go$ and $\F = \F_\gl$.  Let $E$ be an exceptional constituent of $A[\gl]$ and $T=T_E$ the set of bad primes of $E$.  We need a variant of Notation \ref{seldim}.

\begin{Not} \label{seldimc} 
Let $X$ be an irreducible component of $E$ as $\F_2[\Delta]$-module.  For $S\supseteq T$, let $\Lambda^{cr}_{E}(S)$ be the maximal elementary 2-extension $\Lambda$ of $F$ such that
\begin{enumerate}[{\rm i)}]
\item $\Lambda/F$ is unramified outside $\{2, \infty \} \cup (S \backslash T)$,
\item for $\lambda \, \vert \, 2$, the ramification groups $\cI_\lambda(\Lambda/\Q)^\alpha=0$ when $\alpha  > 1$ and
\item $\Gal(\Lambda/F) \simeq X^{* \, r}$ as an $\F_2[\Delta]$-module.
\end{enumerate}
Let $r^{cr}_{E}(S)$ be the multiplicity of $X^*$ in $\Gal(\Lambda^{cr}_E(S)/F)$ and $\Gamma^{cr}_E(S) = \Gal(\Lambda^{cr}_E(S)/\Q)$. Note that $X^* \simeq \widehat{X}$, so  $\Lambda^{cr}_E(S)$ contains $\Lambda_E(S)$. 
\end{Not}

\begin{Rem}\label{MVErem}
Assume $V$ is a semistable $\F[G_\Q]$-module, $T_E \subseteq T_V \subseteq S$ and  
\begin{equation} \label{MVE}
0 \to \F^n \to V \stackrel{\pi}{\to} E \to 0
\end{equation}
is exact. Then $ \Q(V)\subseteq\Lambda^{cr}_E(S)$ by  Remark \ref{ForD1} and Lemma   \ref{ram}.  
\end{Rem}

Suppose $F\subseteq K \subseteq \Lambda^{cr}_E(S)$, with $K$ Galois over $\Q$.  At  bad places $v$ of $E$, let $M_v = (\sigma_v-1)(E)$ and $\cL^{cr}_v = H^1(\cI_v(K/\Q),M_v^*)$.  Set
$$
H^1_{\cL^{cr}}(\Gal(K/\Q),E^*) = \ker\! :    
       H^1(\Gal(K/\Q),E^*) \stackrel{res}{\longrightarrow}
       \prod_{v | N_E}  \cL^{cr}_v.
$$

\begin{lem} \label{crSplitting}
If  $G=\Gal(\Q(V)/\Q)$ and $H^1_{\cL^{cr}}(G,E^*) = 0$, then {\rm (\ref{MVE})} splits.
\end{lem}

\proof
Let $M = \ker \pi \simeq \F^n$ in (\ref{MVE}) and let $v$ lie over a bad prime of $E$.  Since $V$ is semistable, $\sigma_v$ acts trivially on $\pi^{-1}(M_v) = (\sigma_v-1)(V) + M$.  Hence the sequence $0 \to M \to \pi^{-1}(M_v) \to M_v \to 0$ is $\F[\cI_v]$-split.  Corollary \ref{splitcor}(i), with $X = M^*$, $Y = E^*$, $\overline{Y}_i = M_v^*$ and $\overline{V}_i = \pi^{-1}(M_v)^*$ now implies that $X' = X^G = M^*$.  We conclude by duality from Lemma  \ref{makesplit}(i). \qed

\medskip

The following hypothesis will be used to create mirages of the form (\ref{funcmir}).    
\vspace{5 pt}

\setlength{\hangindent}{28 pt}
\noindent \fbox{\bf M3} There is an odd order subgroup $H$ of $G_{\infty} = \Gal(\Q(A[\gl^{\infty}])/\Q)$ such that $E^H=0$ for all $E$ in  $\gS_\gl(A).$

\setlength{\hangindent}{0 pt} 

\begin{lem}
If $\Gal(\Q(E)/\Q)$ is solvable for all $E$ in $\gS_\gl(A)$, then {\bf M3} holds.
\end{lem}

\proof
Since $\Q(A[\gl^\infty])$ is a pro-2 extension of the field $\Q(A[\gl]^{ss})$ generated by the points of all the exceptional Galois $\gol$-modules $E$ in $\gS_\gl(A)$, $G_\infty$   is solvable.  The profinite version of Hall's theorem provides a subgroup $H$ of maximal odd order in $G_\infty$.  Fix $E$ and let $\overline{H}$ be the projection of $H$ to $\Delta_E = \Gal(\Q(E)/\Q).$  Then $\overline{H}$ has maximal odd order in $\Delta_E$.   A minimal normal subgroup  $N$ of $\Delta_E$ is a $p$-group. Since  $E$ is irreducible and $\Delta_E$ acts faithfully, we have $E^N = 0$ and so $p$ is odd.  Hence $\overline{H}$ contains a conjugate of $N$.  It follows that $E^H = E^{\overline{H}} = 0$.     \qed

\medskip

Since $H$ has odd order, the central idempotent
$
e_H=\frac{1}{|H|}\sum_{h\in H} h
$ 
gives a  natural $H$-splitting $M = M^H \oplus (1-e_H) M$  for any $\gol[H]$-module $M$.  Define
$$
D_H = D_H(A) = (1-e_H)\T_\gl(A) = \ker\!:  \T_\gl(A) 
\stackrel{e_H}{\longrightarrow} \T_\gl(A).
$$ 
Then $D_H$ is a pure $\gol$-submodule of $\T_\gl(A)$.  For $\gol$-linear isogenies $\varphi: \, A \to A'$, the functorial property $\varphi(D_H(A)) \subseteq D_H(A')$ holds.  By projection to $A[\gl]$, we obtain an $\F[H]$-module $\overline{D}_H = (1-e_H)A[\gl]$, such that
\begin{equation} \label{rankD}
 \dim_{\F} \overline{D}_H = \dim_{\F} A[\gl] - \epsilon_0(A[\gl]),
\end{equation}
where $\epsilon_0(A[\gl])$ is the number of one-dimensional constituents of $A[\gl]$. Under {\bf M3}, $\dim A[\gl]^H=\epsilon_0(A[\gl])$  and any $E$ in $\gS_\gl(A)$  which is a submodule of $A[\gl]$ lies in $\overline{D}_H.$

\begin{Def} \label{GT} 
We say $E$ is $(S\backslash T_E)$-{\em fissile} if, for every semistable $\F[G_\Q]$-module $Y$ such that $T_Y \subseteq S$, the exact sequence $0 \to \F \to Y \to E \to 0$ splits.  We say {\em fissile} if $S=T_E$.
\end{Def}

Recall that $N_A^u = N_A^u(\gl)$ is the  prosaic $\gl$-conductor of $A$ and $\Pi^u_A$ is the set of primes dividing $N_A^u$, as in Definition\! \ref{Nss}.  See Definition \ref{S-transp} for $p$-transparency.

\begin{theo} \label{bothsplit1}
Assume that {\bf M3} applies and that all $E$ in $\gS_\gl(A)$ are fissile.  Then $
\epsilon_0(A[\gl]) \le \Omega(N_A^u)$, if one of the following holds:
\begin{enumerate}[{\rm i)}]
  \item all primes in $\Pi_A^u$ are $3 \bmod{4}$,  or  \vspace{2 pt}
  \item exactly one $p$ in $\Pi_A^u$ is $1 \bmod{4}$, every subquotient $\cE$ of $A[\gl]$ with Galois module $E$ is $p$-transparent and $\chi_p(q) = -1$ for all other $q$ in $\Pi_A^u$.
\end{enumerate}
\end{theo}

\proof
Recall that $\cP^u$ contains one place of $\Q(A[\gl^\infty])$ for each prime  of  $\Pi_A^u$.  Let $C(A) = (M_t(A,\cP^u,\gl) + D_H)^{sat}$.  If the associated mirage is unobstructed, then Corollary \ref{FaltC} and  (\ref{sizeC}) imply that
$$
\dim_{\F} A[\gl] = \dim_{\F} \overline{C}(A) \le \dim_{\F} \overline{D}_H + \dim_{\F} \overline{M}_t(A,\cP^u,\gl)
\le \dim_{\F} \overline{D}_H + \Omega(N_A^u) 
$$
and our claim follows from (\ref{rankD}).  We therefore assume that $A$ is obstructed. 

Let $X$ be an $\F[G_\Q]$-submodule of $A[\gl]$ of minimal length with exactly one exceptional constituent.  Then we have a filtration $0 \subseteq W \subset X$, with $X/W \simeq E$ in $\gS_\gl(A)$ and $\Q(W)/\Q$ a 2-extension.  Moreover, $W \ne 0$ or else $E$ is a Galois submodule of $\overline{D}_H$ and $A$ is unobstructed.

The corresponding $\F$-module scheme $\cW$ admits a filtration with quotients isomorphic to $\Zl$ or $\Ml$ and conductor $N_W$ dividing $N_A^u.$  Because $X$ is minimal and $E$ is fissile, there is a place $w$ in $\cP^u$ ramified in $\Q(X)$ and unramified in $\Q(E).$   For all such $w$, the action of $\sigma_w$ on $E$ is trivial, so
\begin{equation} \label{RamInX}
 0 \ne (\sigma_w-1)(X) \subseteq \overline{M}_t(A,\cP^u, \gl) \cap W. 
\end{equation}

Consider the filtration (\ref{filt}) on $W$.  If $W = W_1$, then $(\sigma_w-1)(X)$ is a Galois module, violating the obstruction.  Hence $W_1 \subsetneq W_2$ and Proposition \ref{step2} gives $\Q(W_2) = \Q(i)$.  Assuming (i), we have $W=W_2$ by Proposition \ref{step2}(ii).  But then (\ref{RamInX}) violates Lemma  \ref{step1}(ii) and we are done.  From now on, we therefore assume that (ii) holds.  

For each $k$, we have the exact sequence of $\F$-module schemes
\begin{equation} \label{cansplit}
0 \to \cW/\cW_k \to \cX/\cW_k \to \cE \to 0.
\end{equation}
Suppose $W = W_2$.  Then $\cW/\cW_1 \simeq \cZ_\gl^a$ by Lemma  \ref{step1}(iii). By (\ref{RamInX}) and Lemma  \ref{step1}(ii), any odd prime ramified in $\Q(X)$ but not in $\Q(E)$ is $1 \bmod{4}.$  Thus $T_X \subseteq \{ p \} \cup T_E$.  Depending on whether or not $p$ divides $N_E$,  we may use fissility or $p$-transparency on (\ref{cansplit}) with $k = 1$ to contradict minimality of $X$.  Hence $W_3$ contains $W_2$ properly.
By Proposition \ref{step2}(i), $p$ is the only odd prime that may ramify in $\Q(W_3)$, and so $\Q(W_3)$ is in $D_4(-1,p)$.

Let $W_n$ be as defined in Corollary \ref{step3}.  Since that Corollary and (ii) preclude the existence of a prime ramified in $\Q(W_{n+1})$ but unramified in $\Q(W_n)$, we have $W=W_n$.  Now (ii),  Corollary \ref{step3} and (\ref{RamInX}) imply that $T_X \subseteq \{p\} \cup T_E$.  In fact, $p \not\in T_E$ and $p$ must ramify in $\Q(X/W_{n-1})$.  Otherwise, we contradict the minimality of $X$ by using fissility on (\ref{cansplit}) with $k = n-1$.

Let $v$ be the place over $p$ in $\cP^u$.  Because $G_\Q$ acts trivially on $W/W_{n-1}$, we know that $Y = (\sigma_v-1)(X) + W_{n-1}$ is a $G_\Q$-module.  We  claim $Y = W$.  If not, let $W'\supseteq Y$ be  a Galois submodule of codimension 1 in $W.$   Since $\sigma_v$ acts trivially on $X/W'$, the bad primes of $X/W'$ are  in  $T_E$. Then we contradict the minimality of $X$, thanks to  the splitting of  
$
0 \to W/W' \to X/W' \to E \to 0
$
implied by fissility.  Hence $Y = W$ and so  $(\sigma_v-1)(W) = (\sigma_v-1)(W_{n-1}) \subseteq W_{n-2}$.  It follows that $\Q(W/W_{n-2}) = \Q(i)$.  The argument used to prove Lemma  \ref{step1}(iii) shows that $\cW/\cW_{n-1}$ is a direct sum of copies of $\cZ_\gl$.    Minimality of $X$ is contradicted now by applying $p$-transparency to (\ref{cansplit}) with $k=n-1$.  \qed

\medskip

We impose the following assumption for the rest of this subsection.  

\medskip

\setlength{\hangindent}{28 pt}
\noindent \fbox{\bf M4} $A$ is $\go$-paramodular,  $\gS_\gl^{all}(A) = \{\F, \F, E\},$ $\dim_{\F_\gl} E = 2$ and $H^1(\Delta,E) = 0$. Let $T$ be the set of bad primes of $E$ and $S$ that of $A$. 
\setlength{\hangindent}{0 pt}
\medskip

By Proposition \ref{JH2},  $E^* \simeq E$.

\begin{prop} \label{punchline}
Assume {\bf M4}, $E$ absolutely irreducible and $r_E^{cr}(T) = 1$.  If one of the following holds, then $N_A^u(\gl) > 1$:
\begin{enumerate}[\em i)]
\item $F = \Q(E)$ is the maximal real subfield of $\Lambda_E^{cr}(T)$; or
\item $\lambda$ does not ramify in $F$, $|\cD_\lambda(F/\Q)| \le 2$ and $\lambda$  ramifies in $\Lambda_E^{cr}(T)$.
\end{enumerate}
\end{prop}

\proof
For $B$ in $\cI^0_A$, let $\gC(B)$ consist of the $\F$-module subschemes of $B[\gl]$ isomorphic to $\Mu_\gl$ or $\cZ_\gl$ and let $A$ be obstructed for this mirage.  Then there is a filtration $0 \subset \cE \subset \cV \subset A[\gl]$ whose Galois modules have the grading $[E \,\F\,\F]$.  Moreover $V$ does not split and so $L = \Q(V)$ contains $F$ properly, since $H^1(\Delta,E)=0$.

If $N_A^u(\gl) = 1$, we have $L = \Lambda_E^{cr}(T)$, since Remark \ref{MVErem} gives the inclusion and then the argument in Remark \ref{Ds} gives equality because $r_E^{cr}(T) = 1$.  Set $G = \Gal(L/\Q)$ and $H = \Gal(L/F)$.  As in the proof of Lemma \ref{infres},  inflation-restriction, the vanishing of $H^1(\Delta,E)$ and Lemma \ref{BC}(iii) imply that
$$
\dim_\F H^1(G, E) = \dim_\F \Hom_{\F_2[\Delta]}(H,E) = 
              r_E^{cr}(T) \, \dim_\F \End_{\F[\Delta]}E = 1.
$$

Define $W$ by the exact sequence $0 \to E \to A[\gl] \to W \to 0$, so $\gr W =[\F \, \F]$.  If $\Q(W) = \Q$, then $\Q(A[\gl])$ is contained in $L$ by Remark \ref{MVErem}, whence $\Q(A[\gl]) = L$.  By Lemma \ref{makesplit}(i), there is a submodule $W'$ of $W$ whose preimage $V'$ in $A[\gl]$ fits into a {\em split} exact sequence of Galois modules
$
0 \to E \to V' \to W' \to 0
$
wherein $\dim_\F W/W' \le 1$.  But then $W' \ne 0$ and so splitting of this last sequence contradicts the obstruction.  Hence $\Q(W)$ is a {\em quadratic} field.

But $N_W = 1$ and so $\Q(W) = \Q(i)$.  If (i) holds, let $\tau$ be a complex conjugation in $G$. If (ii), let $\tau = \sigma_\lambda$ as in Remark \ref{slam}, since $\cD_\lambda(L/\Q)$ is a 2-group.  Thus $\tau$ is an involution, trivial on $F$, but not on $V$ nor $W$.  This contradicts Lemma \ref{4steps}. \qed

\medskip

	We defined $r_E$   in Notation \ref{seldim} and fissile in Definition \ref{GT}.

\begin{prop} \label{mir1}
Assume {\bf M4},  2 ramifies in $F = \Q(E) $ and $r_E(S)=0,$  where $S$ is the set of  bad places for  $A$ and  $T$ that for $E$.
If $p^*\equiv 1\, (8)$ and all  $q_i^*\equiv  5\, (8),$ then {\em none} of the following occurs: 
\begin{enumerate}[{\rm i)}]
\item $E$ is fissile and $N_A^u = p^a$ or $q_1^a q_2^b$, 
\item $E$ is $q_3$-fissile and $N_A^u = p^a q_3^b$, 
\item $r_E^{cr}(T) = 1$, $N_A^u = p^a$ {\rm (}resp.\! $N_A^u = q_1^a q_2^b${\rm )}  and the primes $v\,|\,p$ {\rm (}resp.\! $v\,|\,q_1$ or $v\,|\,q_2${\rm )} do not split completely in $\Lambda_{E}^{cr}(T)/F$.
\item $N_A^u = p^a q_3^b$, $r_E^{cr}(T\cup \{q_3\}) = 1$ and the primes $v\,|\,p$ do not split completely in $\Lambda_{E}^{cr}(T\cup \{q_3\})/F$.
\item $F$ is the maximal totally real subfield of $\Lambda_{E}^{cr}(T)$ and $N_A^u = p^a$ with $p \equiv 7 \, (8)$.
\end{enumerate}
\end{prop}

\proof
For each $B$ in $\gI_A^\gl$, let $\gC(B)$ consist of all subschemes of $B[\gl]$ isomorphic to $\Ml$ or an $\cE$.  Since 2 ramifies in $\F,$ $\cE_{|  \Z_2}$  is biconnected or $\gr \cE_{|  \Z_2} = [\Ml \, \Zl]$.  Hence $B[\gl^s]$ also has as many $\cZ_\gl$'s as $\Mu_\gl$'s globally.  By Proposition \ref{Faltings}, if $\gC$ is not obstructed, then for some $B,$  no subquotient of $B[\gl^r]$ is isomorphic to $\Zl$, a contradiction.   Thus, assume that $A$ is obtructed.  

By Lemma \ref{Newse0}(i), $A$ is $\Pi_A^u$-transparent because $r_E(S) = 0$ and 2 ramifies in $F$.  This leads to a filtration   
$
0 \subset \cV_1 \subset \cV \subset \cW=A[\gl] 
$ 
with grading $ [\cZ_\gl \, \Mu_\gl \, \cE]$ and $\cV$ not split.  Corollary \ref{core} gives $N_V=p$  or $q_1q_2.$    Write $\cX = \cW/\cV_1$ and $L = \Q(X)$.  By Lemma \ref{super2}, $N_X$ is squarefree and so $\gcd(N_E,N_X/N_E) = 1$.

If an involution $\tau$ is trivial on $E$ but not on $V$, Lemma \ref{4steps} shows that $\tau$ is trivial on $X$.  By choosing $\tau = \sigma_w$ at places $w$ that divide $N_V$  but not $N_E$,  we deduce that $N_X = N_E$ in (i), (iii) and (iv), while $N_X$ divides $q_3N_E$ in (ii).  

In cases (i) and (ii), fissility provides a Galois submodule $E'$ of $X$ isomorphic to $E$.
Such an $E'$ is also available in cases (iii) and (iv).  Otherwise, we have $L = \Lambda_E^{cr}(T)$ in (iii), while $L= \Lambda_E^{cr}(T\cup \{q_3\})$) in (iv).   Any Frobenius $\Phi_v$ in $\cD_v(\Q(W)/F)$ centralizes the generator $\sigma_v$ of inertia in $\cD_v(\Q(W)/F)$.  They are represented by matrices the form
$$
\sigma_v = \left[\begin{smallmatrix} 1 & 1 & * & * \\
                                                         0 & 1 & 0 & 0 \\
                                                         0 & 0 & 1 & 0 \\
                                                         0 & 0 & 0 & 1        \end{smallmatrix} \right]
 \quad      and  \quad
 \Phi_v = \left[\begin{smallmatrix} 1 & * & * & * \\
                                                         0 & 1 & \alpha & \beta \\
                                                         0 & 0 & 1 & 0 \\
                                                         0 & 0 & 0 & 1        \end{smallmatrix} \right],
$$
with $(\alpha,\beta) \ne (0,0)$ by our assumption that $v$ does not split in $L/F$.  But then the matrices do not commute.

Such an $E'$ also is available in (v), where we have $\Q(V) = \Q(\sqrt{-p})$, and we may use  $\tau = \sigma_\infty$ to see that $K = \Q(X)$ is totally real.  But $T_X = T_E = T$, so $K \subseteq \Lambda_E^{cr}(T)$ and therefore $K = F$.  

In all cases, we now have a filtration $\gr \cW = [\cZ_\gl \, \cE' \,  \Mu_\gl]$.  Thanks to the $\Pi_A^u$-transparency of $\cE'$, there is an exceptional $\F$-module subscheme $\cE''$ of $A[\gl],$ violating the obstruction. \qed
 
\begin{lem}\label{mirpar}
Suppose  {\bf M3}, {\bf M4}  and $N_A^u = p^a q$ with $p,q$ primes  not dividing $N_E$. If $E$ is $p$-fissile and all    $\cE$  are  $\{p,q\}$-transparent, then  $p^* \equiv 1 \, (8)$ and $\chi_{p^*}(q) = 1.$
\end{lem}

\proof
Fix $w$ over $q$, write $M_t = M_t(A,w,\gl)$ and let $A$ be obstructed for the mirage associated to $C(A) = (M_t + D_H)^{sat}$.  Then $E$ is not  a Galois submodule of $A[\gl]$.  Let  $V$ be a Galois submodule  of $A[\gl]$ with $\gr V = [\F \, E]$.  If $q$ ramifies in $\Q(V)$, then $\overline{M}_t = (\sigma_w-1)(V)$ is the 1-dimensional Galois submodule of $V$ and  $A$ is unobstructed.  Hence $N_V$ divides $pN_E$.  But then $V$ is split by $p$-fissility, a contradiction. 

We thus  have a filtration $0 \subset \cW_1 \subset \cW_2 \subset \cW = A[\gl]$, in which   $\gr W = [\F \, \F \, E]$ and $E$ cannot move to the left.  The  $\{p,q\}$-transparency implies that $\cW_2$ is a nugget with  $\gr \cW_2 = [\cZ_\gl\,  \Mu_\gl].$  Since $E$ is $p$-fissile, $q$ ramifies in $\Q(W/W_1)$ and so $q$ is unramified in $\Q(W_2)$.  Hence $\Q(W_2) = \Q(\sqrt{p^*})$, with $p^* \equiv 1 \, (8)$.   
Now $(\sigma_w-1)(W/W_1) = W_2/W_1$ and so $W_2 = W_1 + \overline{M}_t$  is a trivial $\cD_w$-module.  Therefore $\chi_{p*}(q) = 1.$  \qed

\begin{prop} \label{3.59}
Suppose {\bf M4},  $N_A^u=  p$ and  $\cD_v$ acts irreducibly on $E$ for  $v\, |\,p$ in $\Q(A[\gl^\infty]).$
Unless $\cW \simeq \Mu_\gl$ and $\cE_{|\Z_2}$ is \'etale, assume that    all exact sequences  
\begin{equation} \label{WVE}
0 \to \cW \to \cV \to \cE \to 0
\end{equation}
of $\F$-module schemes  over $\Z_{T}$ with $W \simeq \F$ are generically split. Then $p \equiv 1 \, (4). $
\end{prop}

\proof 
The irreducibility of $E$ as a $\cD_v$-module and normality of the cyclic 2-group $\cI_v$  imply that $E^{\cI_v} = E$ and so $p$ is unramified in $\Q(E).$ Let $H$ be a cyclic odd Hall subgroup of $\cD_v(\Q(A[\gl^\infty]).$ Then {\bf M3} holds because $E^H=E^{\cD_v}=0.$ 
Since $M_t = M_t(A,v,\gl)$ is a pure $\cD_v$-module of $\go$-rank one,   $\overline{M_t} \cap \overline{D}_H = 0$ and $C(A) = M_t+D_H$ is a pure $\gol$-submodule of $\T_\gl(A)$ of rank 3.  Assume $A$  is obstructed for the associated mirage. 

Suppose $A[\gl]\supset \cV$  for an $\F$-module subscheme  as in (\ref{WVE}),  defined over $\Z_S$ with $S = T \cup \{p\}$.  In Lemma \ref{super1}, we have $\overline{\delta} = 0$, so $f_p(V) = f_p(E) = 0$.  Hence $\cV$ extends to an $\F$-module scheme over $\Z_T$.  By   obstruction, the generic splitting assumption implies that $\cW \simeq \Mu_\gl$ and $\cE_{|\Z_2}$ is \'etale.   Hence $\gr A[\gl]=[\Mu_\gl \, \cE \, \Mu_\gl]$.  Since $E^*\simeq E$, the splitting assumption allows us to move $E$ to the right, creating Galois submodules $X\supset X_1$  with $\gr X = [\F \, \F].$ Since $\overline{M}_t\subseteq X \cap \overline{C}(A)$   and $A$ is obstructed, $X=X_1+\overline{M}_t$ is not  a trivial Galois module and $p$ is unramified in $\Q(X).$ Thus $\Q(X) = \Q(i)$ and  $\cD_v$ acts on $X$ as an odd order group. This implies our conclusion by Lemma  \ref{step1}(ii).  \qed

\section{Small irreducibles and their extensions}\label{Irr} 
The goal of this section is to make the criteria obtained earlier testable, by reducing the study of large Galois extensions to that of more tractable cyclic extensions of smaller fields with precisely controlled conductors. The Bordeaux tables, Maple and Magma then helped with the numerical verifications. 

\numberwithin{equation}{section} \numberwithin{equation}{subsection}
\subsection{Extensions of $E$ by $\F$} \label{2DX}  
Let $\cE$ be a simple $\F$-module scheme whose Galois module $E$ is semistable, self-dual and 2-dimensional over $\F.$    Let $F = \Q(E)$,  $\Delta = \Gal(F/\Q)$ and $\ell = \ch(\F) \not\in T,$  where $T$ is the  set of bad places of $E$. Assume also that $E$ remains irreducible  as $\F_\ell[\Delta]$-module (cf.\! Lemma \ref{BC}).     As in  \cite{Rib3},  $\Mu_\ell \subseteq F,$ since $\det(\rho_E) = \omega$ is the mod-$\ell$ cyclotomic character, and $\Delta \simeq \rho_E(G_\Q)$ is conjugate to a subgroup of
$$
{\rm R}_2(\F) := \{ M \in \GL_2(\F) \, | \, \det M \in \F_\ell^\times\}.
$$
There are {\em transvections} in $\Delta$, i.e. elements $g$ such that \! $\rk_\F(g-1)=1$.   

\begin{lem} [\cite{Rib3},\cite{Suz}]  We have $\Delta={\rm R}_2(\F)$ unless:
\begin{enumerate}[{\rm i)}]
\item  $\ell = 2$ and  $\Delta=D_m\subseteq \SL_2(\F)$, with  $\F$ minimal such that $|\F|\equiv \pm 1\bmod{m},$ or 
\vspace{5 pt}
\item  $\ell = 3,$  $\Delta=\langle \left[\begin
{smallmatrix} 1 & \hspace{6pt} 0 \\ 0 & -1 \end{smallmatrix}\right] ,
 \left[\begin{smallmatrix} 1 & 1 \\ 0 & 1 \end{smallmatrix}\right],
 \left[\begin{smallmatrix} 1 & 0 \\ i & 1 \end{smallmatrix}\right]\rangle$ with  $i\in \F_9$ and $i^2=-1.$ Then $\Delta\cap\SL_2(\F_9)$ is isomorphic to $\SL_2(\F_5).$ 
 \end{enumerate} 
\end{lem}

\begin{lem} \label{nocohom} 
We have $H^j(\Delta,E)=0$  for all $j \ge 0,$ unless $\ell=2$ and $|\F| \ge 4.$
\end{lem}

\proof
Each $\Delta$  contains a non-trivial normal subgroup $\Gamma$ of order prime to $\ell$.  Since $E$ is irreducible, $E^\Gamma=E^\Delta = 0$.  We have  $H^{k}(\Gamma, E) = 0$ for $k\ge 1$ and conclude thanks to the inflation-restriction sequence for $j\ge 1$:
$$
\hspace{20 pt} 0=H^j(\Delta/\Gamma, E^\Gamma) \to H^j(\Delta,E) \to H^j(\Gamma,E)^{\Delta/\Gamma}=0. \hspace{40 pt} \qed
$$

\medskip

In this subsection,  assume $H^1(\Delta,\widehat{E}) = 0$.   For $S \supseteq T$, suppose $\cV$ is an $\F$-module scheme over $\Z_S$ such that the following sequence is {\em not} generically split
\begin{equation} \label{Vseq}
0 \to \cV_1 \to \cV \to \cE \to 0, 
\end{equation}
where $\cV_1 = \cZ_\gl$ or $\F = \F_2$ and  $\cV_1 = \Mu_2.$  Set $L =\Q(V)$ and $G = \Gal(L/\Q)$.  Then $V$ affords a matrix representation:
\begin{equation} \label{rhoV}
\rho_V(g) = \left[\begin{matrix} 1 & x_g \hspace{4pt} y_g \\
0 & \rho_E(g) \end{matrix}\right]  \in \GL_3(\F),
\end{equation}
where $(x_g,y_g)$ is viewed as an element of $\widehat{E} = \Hom_\F(E,\F) \simeq \F \oplus \F.$ The  class $[c]$ in $H^1(G,\widehat{E}),$ associated  to (\ref{Vseq}) does not vanish, even when restricted to
$$
H^1(\Gal(L/F),\widehat{E})^\Delta = \Hom_{\F_\ell[\Delta]}(\Gal(L/F),\widehat{E}).
$$
By irreducibility of $E$ over $\F_\ell$, $\res[c]: \, \Gal(L/F) \to \widehat{E}$ is an isomorphism of $\F_\ell[\Delta]$-modules and so $G$ is a semidirect product
$$
G\simeq \Gal(L/F)\rtimes \Gal(F/\Q) \simeq \left[\begin{smallmatrix}  
        1 & \widehat{E} \\ 0 & I \end{smallmatrix}\right]  \rtimes 
         \left[\begin{smallmatrix} 1 & 0 \\  0 &  \Delta \end{smallmatrix}\right].
$$

\medskip

We describe a subfield $F_1$ of $F$ and an extension $L_1/F_1,$ such that $L$ is the Galois  closure of $L_1/\Q.$ Since any $\ell$-Sylow subgroup $P$ of $\Delta$ fixes a line in $E$ pointwise, assume $P$ is contained in
$
\Delta_1 = \Delta \cap \left(\begin{smallmatrix} 1 & * \\ 0 & * \end
{smallmatrix}\right).
$
Let
$$
F_1 = F^{\Delta_1}, \; G_1 = \Gal(L/F_1) = G \cap \left[\begin
{smallmatrix} 1 & * & *  \\ 0 & 1 & * \\ 0 & 0 & * \end{smallmatrix}
\right],  \;
N_1 = G \cap \left[\begin{smallmatrix} 1 & 0 & *  \\ 0 & 1 & * \\ 0 &  
0 & * \end{smallmatrix}\right] \; \text{and} \; L_1 = L^
{N_1}.
$$
Then $\det\!: \Delta_1/P \simeq \F_\ell^\times$, $N_1$ is normal in $G_1$ and $G_1/N_1 \simeq \F$.  If $F_1 \subset L_2 \subseteq L_1$, the Galois closure of $L_2/\Q$ is $L$ by irreducibility of $E$ as $\F_\ell[\Delta]$-module.  

Lemma \ref{super1} shows that $\cV$ extends to a finite flat group scheme over $\Z_{S'},$
where 
$$
S' = T \, \cup \, \{p_v \in S \backslash T \; | \; E\text{ is not irreducible as an }\F[\cD_v(F/\Q)]\text{-module}  \},
$$
and so we tacitly assume $S=S'.$   Write $\gc_\infty$, $\gc_\ell$ and $\gc_p$ for the semilocal components of the ray class conductor of $L_1/F_1$ at the places over $\infty$, $\ell$ and $p \ne \ell$ respectively.    

\begin{lem} \label{MakeL1}
We have the following bounds on the conductor $\gc(L_1/F_1)$.
\begin{enumerate}[{\rm i)}]
\item  $\gc_p$ divides $p$ if $p$ is in $S \backslash T$ and $\gc_p = 1$ for other $p \ne \ell$.  
\item $\gc_\infty = 1$ unless $F_1$ is totally real, when $\gc_\infty$ is the product of its infinite places.
\item $\gc_\gl = 1$  when $\cV_1 = \cZ_\gl.$ 
\end{enumerate} 
\end{lem}

\proof 
Let $v$ be a prime of $L$ ramifying in $L_1/F_1.$    If $p_v \ne \ell$ is in $S \backslash T,$  then $v$ is tame, with conductor exponent one.
If $v$ lies over $T \cup \{\ell\}$, Lemma \ref{ram} implies that $\cI_v(L/F) = 0.$  Thus $\cI_v(L/F_1)$ contains an element $\sigma_v$ of order $\ell$, not trivial on $L_1,$ such that $\rho_E(\sigma_v) \ne 1.$  It follows that   
$
\rho_V(\sigma_v) = \left[\begin{smallmatrix} 1 & x & y \\ 0 & 1 & a \\ 0 & 0 & 1 \end {smallmatrix}\right]
$
with $xa \neq 0$ and so $(\sigma_v-1)^2(V) \ne 0.$

If $v$ lies over $T$, this contradicts semistability.  If $v$ lies over $\ell$, then $\sigma_v$ acts wildly on $E$, ruling out the possibility that $\cE_{|\Z_\ell}$ be biconnected.  Hence $\gr \cE_{|\Z_\ell} = [\Mu_\gl \, \cZ_\gl]$ in the filtration induced by our fixed basis for $V$ and  $\gr \cV_{|\Z_\ell} = [\cZ_\gl \, \Mu_\gl \, \cZ_\gl]$.   But inertia acts tamely on $\Mu_\gl$, contradicting $x \neq 0$.

Suppose $\ell=2$ and $F_1$ has a complex place, whence $F$ is totally complex.  If $\sigma_v$ is complex conjugation for $v$ lying over a real place of $F_1$, then $\rho_E(\sigma_v) \ne 1$.  But $\sigma_v$ fixes $F_1$ and so $\rho_V(\sigma)$ is upper triangular. If $v$  ramified in $L_1/F_1$, we would have the same  contradiction as  for $v$ over $T$, since $\sigma_v^2 = 1$.  \qed

\medskip

For the rest of this  section, assume  $\F = \F_2$, so $\ell=2,$  $\Delta \simeq \SL_2(\F_2)$ and $F_1$ is a cubic field.  Moreover, $E\simeq \widehat{E}$ as Galois modules, $H^1(\Delta,E) = 0$  and $\Gal(L/\Q)\simeq \cS_4.$  

Define the prime $\lambda_1 \, | \, 2$ in $F_1$ according to the factorization of $(2)\cO_{F_1}$:
\begin{equation} \label{Over2}
(2)\cO_{F_1} = 
\begin{cases}
    \lambda_1^3              &\text{if } e_{\lambda_1}(F_1/\Q) = 3, \\
    \lambda_1^2 \lambda_1'   &\text{if } e_{\lambda_1}(F_1/\Q) = 2, \\
    \lambda_1 \lambda_1'     &\text{if } f_{\lambda_1}(F_1/\Q) = 2.
\end{cases}
\end{equation}

\begin{lem} \label{gentran}  
If $\cV_1 = \Mu_2$ in {\rm (\ref{Vseq})}, then $\gc_2(L_1/F_1)$ divides $4.$ It even divides $\lambda_1^2$ if: \quad \rm (i) $2$ ramifies in $F/\Q$ \quad or  \quad  (ii) \ $f_\lambda(F/\Q) = 2$ and $\cE^m \ne 0$ over $\Z_2$.  
\end{lem}

\proof  
Conductors of small extensions of $\Q_2$ may be found by direct calculation or in the Tables of  \cite{JR}, where the  last entry of {\em Galois Slope Content} is at most 2 exactly when the higher ramification bound  in Lemma \ref{FontaineBound} holds.

Assume $\lambda$ ramifies in $L/F$.  If $\cE_{| \Z_2}$ is biconnected, $\cD_\lambda(L/\Q) \simeq \cS_4$ and $\cI_\lambda(L/\Q) \simeq \cA_4$.  By \cite{JR} for sextics over $\Q_2$, we  have $\ord_2(d_{L_1/\Q}) = 6$, whence $\gc_2(L_1/F_1) = \lambda_1^2$ by the conductor-discriminant formula.  The end of the proof of Lemma \ref{condL0} gives an explicit description of the completion $L_\lambda$.

When $\cI_\lambda(L/\Q)$ is a 2-group, Lemma \ref{FontaineBound} implies that $\cI_\lambda(L/\Q)_2 = 1$.  Passing between lower numbering for subgroups and upper numbering for quotients, we find that $\cI_\lambda(L_1/F_1)_2 = 1$ and so the conductor exponent of $L_1/F_1$ at $\lambda$ is 2 by \cite[Ch.\! XV,\S2,Cor.\! 2]{Ser1}. 

Now assume (i) with $e_\lambda(F/\Q) = 2$, or (ii).   Fix the primes $\lambda$ and $\lambda'$ of $L$ over $\lambda_1$ and $\lambda_1'$.  Since $\lambda'_1$ splits in $F_1/\Q$, $\cD_{\lambda'}(L/\Q)$ is contained in $G_1$.  The non-trivial action of $\cD_{\lambda'}(F/\Q)$ on our basis for $V$ implies that $\gr V =  [\Mu_2 \, \Mu_2 \, \cZ_2]$ at $\lambda'$.  But then $\cI_{\lambda'}(L/\Q) \subset N_1$, so is trivial on $L_1$.  Thus $\lambda'_1$ is unramified in $L_1/F_1$ and $\gc_2(L_1/F_1)$ divides $\lambda_1^2$.  \qed

\begin{lem} \label{2inert}
Let $L_1$ be a sextic field whose Galois closure $L$ is an $\cS_4$-field with  $F$ as its $\cS_3$-subfield.
If $L/F$ is unramified over 2 and one of the following holds for $\lambda \, | \ 2$, then $\vert \cD_\lambda(L/F) \vert = 2$.
\begin{enumerate} [{\rm i)}]
\item $e_\lambda(F/\Q) = 2$ and there are exactly 3 primes over 2 in $L_1$, or
\item $e_\lambda(F/\Q) =1$, $f_\lambda(F/\Q) = 2$  and there are exactly 2 primes over 2 in $L_1.$
\end{enumerate}
\end{lem}

\proof
A 2-Sylow subgroup $G_1$ of $G = \Gal(L/\Q) \simeq \cS_4$ cuts out the cubic subfield $F_1$ and the subgroup $N_1$ generated by the two transpositions in $G_1$ cuts out $L_1$.  The subgroup $\kappa$ of $G$ generated by the even involutions cuts out $F$.  Write $\lambda$, $\lambda'$ for primes of $L$ over $\lambda_1$,  $\lambda'_1$ respectively.  Then $\cD_{\lambda'}(L/\Q)$ is contained in $G_1$ because $\lambda'_1$ is split in $F_1/\Q$. 

If $\cD_\lambda(L/F) = 1$, then $\cD_\lambda(L/\Q)$ has order 2 and is not trivial on $F$, so it is generated by a transposition.  Thus $\cD_\lambda(L/F_1)=\cD_\lambda(L/L_1)$ and the two primes above 2 in $F_1$ split into 4 in $L_1.$ 
 
Suppose there is a residue extension over 2 in $L/F$, so $\cD_{\lambda'}(L/\Q)$ has order 4.  

\begin{enumerate}[i)]
\item If $e_{\lambda'}(F/\Q) = 2$, then $\cI_{\lambda'}(L/\Q)$ is generated by a transposition $\sigma_{\lambda'}$.  Hence $\cD_{\lambda'} = N_1$ and $\lambda_1'$ splits in $L_1/F_1$.  Let  $\lambda=\gamma(\lambda'),$ with $\gamma \in G$ of order 3.  Then $\cD_\lambda(L/F) \cap N_1 = 1$ and so $\lambda_1$ is inert in $L_1/F_1$. 
\item  If $f_{\lambda'}(F/\Q) = 2$, then $\cD_{\lambda'}(L/\Q)$ is cyclic, generated by a Frobenius and so $\cD_{\lambda'}(L/\Q) \cap N_1$ is generated by the unique even involution in $N_1$.   It follows that $\lambda'_1$ is inert in $L_1/F_1$.    Conjugating by $\gamma$, we find that $\cD_\lambda(L/\Q) \cap N_1 = 1$.  Hence $\lambda_1$ also is inert in $L_1/F_1$.   \qed
\end{enumerate}

\begin{Rem} \label{DoFFE}
Let $A$ be a hypothetical $(\go,N)$-paramodular variety with $\vert \F_\gl \vert = 2$ and $\gS_\gl(A) = \{E\}$, where $\dim_{\F_2} E = 2$.  Then  $N_E$ is a squarefree divisor  of $N$ and the $\cS_3$-field $F = \Q(E)$ can be constructed by class field theory  or Magma, as a cyclic cubic over $\Q(\sqrt{\pm N_E})$, ramified only over $2\infty$, with $F_1/\Q$ as cubic subfield.  

Let $S$ contain the set $T$ of primes dividing $N_E$.  If no quadratic $L_1/F_1$ satisfies the bounds in Lemma \ref{MakeL1},  then $r_E(S) = 0$ and  extensions  (\ref{Vseq}) over $\Z_S$ with $\cV_1 \simeq \cZ_2$ are generically split.  Lemma \ref{Newse0} controls the deficiency $\delta_A(E)$ in Theorem \ref{trivial}.  When $r_E(S) = 0$ and $2$ ramifies in $F$, $E$ is $(S \backslash T)$-transparent.  When $r_E(S) = 0$ and 2 is unramified in $F$, with residue degree 2, we get only $\delta_A(E) \le 1$, due to the fickle nature of group schemes $\cE$ corresponding to $E$.   If only one quadratic extension $L_1/F_1$ satisfies the conductor bound, then $r_E(T) = 1$, as required by {\bf D}3 in \S \ref{except}.  Lemma \ref{2inert} serves for testing {\bf D}4. 

Now retain the bounds in Lemma \ref{MakeL1}(i),(ii) at odd places, but invoke the weaker bounds on $\gc_2(L_1/F_1)$ in Lemma \ref{gentran}.  When no quadratic $L_1/F_1$ exists, $r_E^{cr}(S) = 0$ and $E$ is $(S \backslash T)$-fissile.  Further, under \ref{gentran}(ii), the splitting required in Lemma \ref{3.59} holds.  Finally, $r_E^{cr}(S) = 1$ if exactly one quadratic $L_1/F_1$ exists.
\end{Rem}

\subsection{Extensions of $E_2$ by $E_1$}\label{EE} 
\numberwithin{equation}{subsection}
Here  $\vert \F_\gl \vert = 2$ and $\gS_\gl^{all}(A)=\{E_1,E_2\}$ with $E_i$  two-dimensional non-isomorphic Galois modules and $\text{cond}(E_i) = N_i$.
By Corollary \ref{2Es}, we may assume that $W=A[\gl]$ is a non-split extension: 
\begin{equation} \label{EWE}
0 \to E_1 \to W \to E_2 \to 0.
\end{equation}
For $i=1,\,2$, set $F_i = \Q(E_i)$, $F=\Q(E_1,E_2)=F_1F_2$ and $\Delta_i =\Gal(F/F_{3-i}) \simeq \Gal(F_i/\Q) \simeq \SL_2(\F_2)$, so $\Gal(F/\Q)= \Delta_1 \times \Delta_2$.  Let $L = \Q(W)$ and $G = \Gal(L/\Q)$.

\begin{lem}
$L$ contains $F$ properly.
\end{lem}

\proof  As $\Delta_1$-modules,  $E_2 \simeq \F_2^2$  and  $\cH=\Hom_{\F_2}(E_2,E_1) \simeq E_1^2.$  Assume $L=F$ and use inflation-restriction for exactness  of the sequence
$$
 H^1(\Gal(F_2/\Q), \cH^{\Delta_1}) \to H^1(\Gal(F/\Q), \cH) \to H^1
(\Delta_1, \cH).
$$
The first term vanishes since $\cH^{\Delta_1} \simeq \Hom_{\F_2[\Delta_1]}(E_2,E_1) = 0$ and the last because $H^1(\Delta_1,E_1)=0$. Thus the middle term is trivial and (\ref{EWE})  splits.   \qed

\medskip

Let $\rho_i$ be the Galois representation afforded by $E_i$, fix a basis $w_1,w_2$ for $E_1$ and extend by $w_3,w_4$ to a basis for $W$.  Then $G$ admits a representation of the form
$$
\rho: g \mapsto \left[ \begin{matrix} \rho_1(g) & B(g) \\ 0 & \rho_2(g) \end{matrix} \right]
$$
Conjugation by $\left[ \begin{smallmatrix} X & Z \\ 0 & Y \end{smallmatrix} \right]$ on 
$
\left[\begin{smallmatrix} I  & B \\ 0 & I \end{smallmatrix} \right]
$ 
in $\Gal(L/F)$ yields 
$
\left[\begin{smallmatrix} I & XBY^{-1} \\ 0 & I \hspace{8 pt} \end{smallmatrix} \right].
$
Since $M_{2 \times 2}(\F_2)$ is $\F_2[\Delta_1  \times \Delta_2]$-irreducible, $\rho$ maps onto the parabolic subgroup indicated above.  

Let $H$ be the 2-Sylow subgroup of $G$ whose image under $\rho$ is the group of all unipotent upper triangular matrices.  Its fixed field $K = L^H$ is the compositum of the cubic fields $K_i = K \cap F_i$.    Let $J$ be the subgroup of $H$ with $B(g)$ upper triangular.  Then $L_0 = L^J$ is a quadratic extension whose Galois closure over $\Q$ is $L$.  Let $\gc_2$ be the 2-part of the ray class conductor of $L_0/K$.  Write $N_1N_2Q$ for the  Artin conductor of $W$  and let $Q_0$ be the part of $Q$ prime to $N_1N_2$.

\begin{lem} \label{condL0}
The extension $L_0/K$ is unramified outside $2\infty Q_0$ and ramifies at $\infty$ only if $F$ is totally real.  Assume that ${\cE_2}$ is biconnected at 2.
\begin{enumerate}[{\rm i)}]
\item If ${\cE_1}$ is biconnected at 2, then $2\cO_K =  
(\lambda_K \lambda_K')^3$ and $\gc_2$ divides $(\lambda_K \lambda_K')
^2$.
\item If  ${\cE_1}$ at 2 is a non-split extension of $\cZ_2$  
by $\Mu_2$, then
$\gc_2$ divides $\lambda_K^2$, where $\lambda_K$  is the unique prime  
of $K$ with $[K_{\lambda_K}\!:\!\Q_2] = 6$. The other prime over 2 splits.
\end{enumerate}
\end{lem}

\proof
At an archimedean or odd place $v$ of $L$ ramified in $L/K$, the inertia group $\cI_v(L/K)$ inside $\Gal(L/K)$  is generated by an involution $\sigma$, as in Lemma  \ref{4steps}. If $v$ ramifies in $L_0/K$, then $\sigma$ does not fix $L_0$, so $c=1$ and thus $x=y=0.$  Hence $v$ does not ramify in $F/\Q$.  It follows that either $v$ lies over $Q_0$ or $v$ is archimedean and $F$ is totally real. 

Let $\pi$ be a root of $x^3-2$ in $\overline{\Q}_2$. In (i), we have $2\cO_K = (\lambda_K \lambda_K')^3$, since $F_i \otimes \Q_2  \simeq \Q_2(\Mu_3, \pi)$. The bound on $\gc_2$ is  in \cite[Prop.\! 6.4]{Sch1}.

\vspace{2 pt}

For (ii), let  $\lambda'$ be a prime of  $L$ with $K_{\lambda'} = \Q_2(\pi)$.  Then the completion of $K_1$ at $\lambda'$ is $\Q_2$ and $\cD_{\lambda'}(F_1/K_1) = \Gal(F_1/K_1)$ has order 2.  Hence the connected component $E_1^0$ at $\lambda'$ is the subspace $\langle w_1 \rangle$.   For any $\sigma$ in $H=\Gal(L/K)$, we have $(\sigma-1)(w_3) \in E_1$.  If, in addition, $\sigma$ is in $\cD_\lambda(L/K)$, then $(\sigma-1)(W) \subseteq W^0$ because $W^{et}$ is 1-dimensional.  Hence $(\sigma-1)(w_3)$ is in
$
W^0 \cap E_1 = E_1^0 = \langle w_1 \rangle,
$
so $\sigma$ is  in $J$ and $\lambda'$ splits in $L_0/K$. 

Suppose  $\lambda$ over 2 in $L$ ramifies in $L_0/K$ and  let $W^0$ be the connected component at $\lambda$. Then $K_\lambda=\Q_2(\pi, \sqrt{d})$ with $d\in\{-1,3,-3\}$  and $\Q_2(W^0)$ is the unique $\cS_4$-field $M$ over $\Q_2$ satisfying  
Fontaine's bound  (cf.\! \cite{JR}). Explicitly,
$$
M  = \Q_2(\zeta, \pi, \sqrt{1+2\pi^2},\sqrt{1+2\zeta\pi^2}),
$$
with $\zeta$  a primitive cube root of unity.  Further,  $\cD_\lambda\simeq \cS_4\times \cS_2$ if 2 ramifies in $F_1$ and $\cD_\lambda\simeq \cS_4$  otherwise.

We use tildes for the completions of various fields at $\lambda$.  If $d \equiv 3 \bmod{4}$, then $\tilde{F} = \tilde{K}(\zeta) = \Q_2(\pi, \zeta, i)$ and $\tilde{L} = M(i)$.  The abelian conductor exponents of $\tilde{L}_0/\tilde{K}$ and $\tilde{L}_0\tilde{F}/\tilde
{F}$ are equal since $\tilde{F}/\tilde{K}$ is unramified.  Local class field theory or the conductor-discriminant formula implies that every quadratic extension of $\tilde{F}$ inside $\tilde{L} = \tilde{F}(\sqrt{1+2\pi^2},\sqrt{1+2\zeta\pi^2})$ has conductor exponent 2.  Hence $\gc_2$ divides $\lambda_K^2$, where  
$\lambda_K$ lies below $\lambda$ in $K$.

If $d=-3$, then $\tilde{L} = M$, $\tilde{K} = \Q_2(\zeta, \pi)$ and the same method applies.  \qed

\subsection{Wherein $A[\gl]$ is irreducible and $\vert \F_\gl \vert = 2$} \label{E4}
Let $A$ be a semistable $(\go,N)$-paramodular abelian variety with  odd $N$ and $S$ the set of primes dividing $N$.  Suppose $\gl$ has a totally positive generator and $E = A[\gl]$ is irreducible.  We  may assume $A$ is minimally $\go$-polarized. Then $\cE$ is Cartier self-dual by Lemma \ref{flip}(ii) and $E$ is symplectic.  

  Let $F = \Q(E)$ and $G = \Gal(F/\Q)\subseteq \cS_6,$ via the induced action  on the set $\Theta^-$ of six odd theta characteristics  \cite[\S 2, \S 4]{BK3}.  When $N$ is not a perfect square, $G$ can only be $\cS_5$, $\cS_6$ or the wreath product $\cS_3\wr\cS_2$.  The last is a group of order 72, isomorphic to the orthogonal group O$_4^+(\F_2)$, as in \cite[p.\! 409]{Pol}.  
The field  $F$ is the Galois closure of a subfield $K$ of degree 5 or 6  fixed by the stabilizer of an odd theta.  Let $d_K$ be the absolute discriminant of $K$.  See \cite[Prop.\! 4.1]{BK3} for $\ord_2(d_K)$.  For odd $p$, \cite[Prop.\! 3.11]{BK3} gives $\ord_p(d_K)$ in terms of the toroidal dimensional of the special fiber, cf.\! Notation \ref{torp}.

\begin{prop} \label{2-disc} We have $\ord_p(d_{K})= t_p,$ unless $t_p=2$ and $ [K\!:\!\Q] = 6,$ when $2\le\ord_p(d_{K})\le 3.$   In addition, $\ord_2(d_{K}) \le 6$  if $ [K\!:\!\Q] = 6$ and $\ord_2(d_{K}) \le 4$  if $ [K\!:\!\Q] = 5.$
 \end{prop}

\begin{lem} \label{bdres}
No prime over 2 in $K$ has residue degree 5.  If $[K\!:\!\Q]=5$, then  
no prime over 2 in $K$ has residue degree 3.
\end{lem}

\proof
If the residue degree $f_\lambda(K) = 5$, then 2 is unramified in $K$ and so in  $F.$   Hence $\cE_{\vert \Z_2}$ is ordinary and the connected-\'etale sequence implies that the order of $\cD_\lambda(F/\Q)$ divides 48, a contradiction.
 
If $f_\lambda(K) = 3$, then $e_\lambda(K) \le 2$ and again $\cE_{\vert \Z_2}$ is ordinary. Any element $\Phi$ of order 3 in  $\cD_\lambda(F/\Q)$ is fixed point free on $E^0$ and on $E^{et}$, so also on $E.$  But if $[K\!:\!\Q]=5$, then $\Phi$  fixes three odd thetas pointwise and the difference of any two corresponds to a fixed point for the action of $\Phi$ on $E$. \qed

\begin{Rem} \label{bs} Let $N$ be  the conductor of $E$ and let $K$ be a sextic field, as above.  Write  $N=N_1N_2^2N_3^2, $ 
where the involutions generating inertia above $p\,|N_i$  are the product of  $i$ transpositions   in $\cS_6.$ 
Then the discriminant of $K$ divides $2^6 N_1N_2^2N_3^3$.  The discriminant of its  {\em twin field} \cite{Rob}, whose representation is twisted by the outer automorphism of $\cS_6,$ divides $2^8 N_1^3 N_2^2 N_3$, because one has weaker control at primes over 2, while products of 1 and 3 transpositions are switched.  When both discriminants exceed 200,000, a totally complex field $K$ {\em might} exist and lie beyond the tables in \cite{BT}.   The conductors $N< 1000 $ for which this issue arises are $5^2\!\cdot\!N_1$ with $29\le N_1\le 39$, $7^2\!\cdot\!N_1$ with $11\le N_1\le 19$ and $11^2\!\cdot\!7,$ all with $N_2 = 1$.   The solvable case $\cS_3 \wr \cS_2$ does not occur by class field theoretic calculation. John Jones kindly  verified, with his targeted searches, that no such $\cS_6$ field exists either.  
\end{Rem}

\section{The Modified Conjecture} \label{Mod}  \numberwithin{equation}{section}
Frank Calegari \cite{Cal} was kind enough to point out an oversight in Conjecture \ref{wcon}.   As a result, we propose the following modification.

\begin{Def} \hfill
\begin{enumerate}[ i)]
\item The abelian variety  $B/K$ is a {\em QM abelian variety} or {\em QM by $D$ abelian variety}  if $\End_K B$ is an order in the  non-split quaternion algebra $D/\Q.$   \vspace{2 pt}

\item A cuspidal, nonlift Siegel paramodular newform $f$ of genus $2$, weight $2$ and level $N$ with rational Hecke eigenvalues will be called a {\em suitable paramodular form of level $N$}. 
\end{enumerate}
\end{Def}

\begin{Con}   \label{modcon}
Let $\mathcal{A}_N$ be the set of isogeny classes of abelian surfaces $A/\Q$ of conductor $N$ with $\End_\Q A=\Z$,
let $\mathcal{B}_N$ be the set of isogeny classes of QM abelian fourfolds $B/\Q$ of conductor $N^2$ and let $\mathcal{P}_N$ be the set of suitable paramodular forms of level $N$, up to nonzero scaling.   There is a bijection $\mathcal{P}_N \leftrightarrow \mathcal{A}_N\cup \mathcal{B}_N$ such that $$L(C,s)=L(f,s,\text{spin})\text{ if } C\in \mathcal{A}_N\text{ and } L(C,s)=L(f,s,\text{spin})^2\text{ if  }C\in  \mathcal{B}_N.$$
\end{Con}
 
Let $B/K$ be a QM by $D$ abelian variety over a number field $K$.   Since $D$ acts on the tangent space, the dimension of $B$ is even, say $\dim B=2d$.  Up to isogeny over $K$, we may assume that $\cO$ is a maximal order and henceforth do so. The actions of $\cO$ and $G_K=\Gal(\overline{K}/K)$ on the Tate module $\T_\ell(B)$ commute. If $\ell$ is a sufficiently large prime, then $\cO_\ell = \cO \otimes \Z_\ell$ is isomorphic to $M_2(\Z_\ell)$ and so  the commutant of $\cO_\ell$ in $\End(T_\ell)\simeq M_{4d}(\Z_\ell)$ is isomorphic to $M_{2d}(\Z_\ell)$.  Thus, the $\ell$-adic representation of $G_K$ on $\T_\ell(A)$ factors:
\begin{equation}\label{subalg}
\rho_\ell\!: \, G_K  \hookrightarrow \GL_{2d}(\Z_\ell)\hookrightarrow \GL_{4d}(\Z_\ell).
\end{equation}
When $\gp$  is a prime of $K$ of good reduction, the Euler polynomial of $A$ at $\gp$ is the characteristic polynomial of a Frobenius Frob$_\gp$ as a matrix in $M_{4d}(\Z_\ell)$. Since $\rho_\ell(\text{Frob}_\gp)$ lands in the smaller subalgebra $M_{2d}(\Z_\ell)$, the Euler factor is a perfect square. A similar argument applies for primes of bad reduction by acting on the part of the Tate module fixed by inertia. Hence $L(B,s)=M(s)^2$ for a Dirichlet series $M(s)$ with coefficients in $\Z$.  See \cite{Chi} for a related argument.

Given a suitable paramodular form $f,$ an abelian surface $A$ with $\End_\Q A=\Z$ and a QM abelian fourfold $B$ cannot simultaneously satisfy $L(A,s)=L(f,s,\text{spin})$ and $L(B,s)=L(f,s,\text{spin})^2$.  Indeed, since the coefficients of a Dirichlet series are known from the function it represents, the $L$-series of an abelian variety determines the trace of Frobenius on the Galois representation of its Tate module.  If $A$ and $B$ both correspond to $f$, then the Galois representations of $A^2$ and $B$ are equal and so $A^2$ and $B$ are isogenous by Faltings.   But they have different endomorphism algebras, namely $M_2(\Z)$ and $\cO$, respectively.   

Let $S_C$ be the set of bad primes and $N_C$ the conductor of an abelian variety $C$.

\begin{prop} \label{mod1}
Let $B/\Q$ be a QM  abelian fourfold.  Then $N_B=n^4 m^2$, where the support of $n$ is $S_B$ and $m$ divides $\gcd\{30,n\}$.  If $B$ is semistable, then $m = 1$.
\end{prop}

\begin{prop} \label{mod2}
Let $A$ be a QM abelian surface over an imaginary quadratic field $K$.  Then all bad primes of $A$ are potentially good and $N_A=\gn^4 \gm^2$, where the support of $\gn$ is $S_A$ and  $\gm$ divides $\gcd\{6,\gn\}$.  
\end{prop}

\begin{proof}
We now assume that $A$ is an abelian variety of dimension $2d$ over a local field $F$ with maximal ideal $\gp$ and residue characteristic $p$.  In addition, $A$ has endomorphisms over $F$ by a maximal order $\cO$ of $D$.   Let  $V=\T_\ell(A)\otimes \Q_\ell$, where $\ell$ is a rational prime different from $p$ and let $I$ be the inertia subgroup of $\Gal(F(A[\ell^\infty])/F)$.   
Recall \cite[\S4]{Gro} that the conductor exponent of $A$ is given by
\begin{equation}
\gf_\gp(A)=\epsilon(A)+\delta(A[\ell]),
\end{equation}
where $\epsilon(A)=\dim_{\Q_\ell}(V/V^{I})$ is the tame conductor exponent, $\delta(A[\ell])$ is the wild or Swan conductor and each term is independent of $\ell$. 

First we treat the tame conductor.  Let $\cA^0$ be the connected component of the closed fiber of the N\'{e}ron model for $A$ over the ring of integers of $F$.  If $t$, $u$ and $a$ denote the dimensions of the toric subgroup, the unipotent subquotient and the abelian variety quotient of $\cA^0$, then $t+a+u = 2d$ and $\epsilon(A) = t + 2u$.   Note that $\epsilon(A)$ and $t$ are multiples of 4, since $\cO$ acts on $V/V^I$ and on the character group of the torus, which is a free $\Z$-module of rank $t$.   Hence $u$ is even.  This leaves very few possibilities for surfaces and fourfolds with bad reduction:
$$
\begin{array}{| c || c | c | c | c |}
\hline
& t & u & a & \epsilon(A) \\
\hline
\text{surface} & 0 & 2 & 0 & 4 \\
\hline 
\text{fourfold} & 4 & 0 & 0 & 4 \\
\cline{2-5}
                      & 0 & 2 & 2 & 4 \\	 
\cline{2-5}
                      & 0 & 4 & 0 & 8 \\
\hline          
\end{array}
$$
The contributions from $\epsilon(A)$ account for $n$ and $\gn$ in the Propositions.  In particular, when $A$ is semistable (i.e. $u = 0$), there is no wild ramification in $F(V)/F$ and no further contribution to the conductor. 

Pass to an extension of $F$ over which $A$ becomes semistable.  Then the dimensions of the torus and abelian variety in the bad fiber become $\tilde{t} \ge t$ and $\tilde{a} \ge a$, with $\tilde{u} = 0$.  Since $\tilde{t}$ also is a multiple of 4, $\gp$ is potentially good when $u = 2$ in the table above.

For the wild conductor, let $L/F$ be a finite Galois extension with higher ramification groups $G_i$ in Serre's lower numbering \cite[IV]{Ser1}, where $G_0$ is the inertia subgroup of $G = \Gal(L/F)$. If $W$ is an $\F_\ell[G]$-module, then $\delta(W)$ is defined by 
\begin{eqnarray}
\delta(W)&=&\sum_{i=1}^\infty \frac{|G_i|}{|G_0|}\dim_{\F_\ell} (W/W^{G_i})
\end{eqnarray}
and is known to be an integer.   Take $L = F(A[\ell])$ and $G = \Gal(L/F)$.  When $\ell\!>\!>0$, we have  $\cO/\ell\cO\simeq M_2(\F_\ell)$.  Since the actions of $G$ and $\cO$ on $A[\ell]$ commute, we have inclusions analogous to \eqref{subalg}:  
\begin{equation}\label{subalg1}
G \hookrightarrow \GL_{2d}(\F_\ell)\hookrightarrow \GL_{4d}(\F_\ell),
\end{equation}
leading to an $\F_\ell[G]$-representation $W$ of  dimension $2d$ such that $\delta(A[\ell])=2\delta(W)$.   

Let $A$ be a fourfold, so $\dim W = 4$ and $|\GL_{4}(W)|=(\ell^4-1)(\ell^3-1)(\ell^2-1)(\ell-1)\ell^6$.  Taking $\ell$ to be a primitive root mod $p$ shows that if $G_1 \ne 0$, then $p\le 5$.   Similarly, for QM abelian surfaces $B$, wild ramification necessitates $p \le 3$.  
\end{proof}

\begin{Rem}
If $B$ in Proposition \ref{mod1} is paramodular, the level of the corresponding paramodular form $f$ is $n^2 m$.  If $A$ in Proposition \ref{mod2} is Bianchi modular, then the level of $f$ is $N_{K/\Q}(\gn^2 \gm)$.
\end{Rem}

The only QM abelian fourfolds known to us are  Weil restrictions $B=R_{K/\Q}(A)$ of QM abelian surfaces $A$ over imaginary quadratic fields $K,$ with $A$ not isogenous to its conjugate and Bianchi modular \cite{Sch}.  Since \cite{Mil} implies that $N_B=d_K^4 N_{K/\Q}(N_A)$, those examples  satisfy  Conjecture \ref{modcon}, thanks to \cite{BDPS}. 

Conjecture \ref{modcon} presents two tasks not yet dealt with:
\begin{enumerate}[i)]\item eliminate QM abelian fourfolds $B$ of certain conductors;\item if a suitable paramodular form exists and there is no corresponding abelian surface, find a corresponding QM abelian fourfold.  
\end{enumerate}

As to (i), Proposition \ref{mod1} severely limits  the conductor of a QM abelian fourfold and thus the level  of a corresponding paramodular form.  The computational evidence in our appendix deals only with  semistable abelian varieties of odd conductor at most 1000.  Recall  Schoof's beautiful result:

\begin{theo}[\cite{Sch2, Sch3, Sch4}] \label{sc29} 
Any $\Q$-simple semistable abelian variety with good reduction outside the odd square-free integer $n \le 23$ is isogenous to  $J_0(n).$   
\end{theo}
It follows that  for odd $n\le27$, semistable abelian surfaces of conductor $n^2$ with $\End A=\Z$ and semistable $QM$ fourfolds of conductor $n^4$ do not exist. According to \cite{PoYu2}, there is no suitable paramodular form of square level at most $1000$, so that neither $n = 29$ nor $31$ should occur.

\appendix

\section{How conductors are ruled out} \label{DATA}   \numberwithin{equation}{section} 
Assume  $\go$ has a prime $\gl \, \vert \, 2$ of degree one.   For each odd integer $N < 1000$, we considered all Galois structures available for $A[\gl]$ when $A$ is  a semistable $(\go,N)$-paramodular abelian variety of reduced conductor  $N.$ We examine in detail the various possibilities for $\gS_\gl(A),$ the multiset of irreducible constituents of $A[\gl]$ of dimension at least 2 over $\F_2$. We say that a given multiset is ruled out for  $N$ if no such $A$ exists.

\begin{enumerate}[I.]
\item $\gS_\gl(A)$ is empty.  Then  $\Gal(\Q(A[2])/\Q)$ is a $2$-group.  All odd $N<1000$, except those marked {\bf u} in Table \ref{BadCond}, are  ruled out by Corollaries \ref{nilpgenbd}, \ref{trivialplus}(ii)b and the criteria in \S \ref{UnipMir}. The Jacobian  of  $ y^2=(x^2+2x+5)(x^4+2x^3+3x^2-2x+1),$ is of type  {\bf u} 
and has the least  conductor $N=1649=17.97$ known to us.
\vspace{5 pt}

\item $\gS_\gl(A)=\{E\}$, with $\dim E = 2$.  Set $F = \Q(E)$ and denote residue and ramification degree at $\lambda \, | \, 2$ by $f_\lambda$ and $e_\lambda$ respectively.  See Remark \ref{DoFFE} regarding the relevant invariants.
\begin{enumerate}[A.]
\item 2 ramifies in $F$.  We have $\delta(E) = 0$ for 92 cases with $e_\lambda(F/\Q) = 3$ and 64 cases with $e_\lambda(F/\Q)=2$.   Then Theorem \ref{trivial} rules out $N_E$ and $qN_E$ when $q  \equiv \pm 3 \bmod{8}$ is prime.   Corollary \ref{trivialplus}(i)c rules out $q^2 N_E$ when $r_E(T\cup\{q\})=0$.  Nine more cases with $e_\lambda(F/\Q) = 3$ are ruled out by Proposition \ref{punchline}(i).  We also use Lemma \ref{specialshape} and Theorem \ref{trivial} to eliminate $N_E$ when $e_\lambda(F/\Q) = 2$
\item  2 is unramified in $F$.  Lemma  \ref{selfdual} precludes  existence if $f_\lambda(F/\Q) = 3$. Proposition \ref{punchline}(ii) rules out 431 and 503, the only conductors available for $E$ if $f_\lambda(F/\Q) = 1$.   Finally, consider $f_\lambda(F/\Q) = 2$.  Rule out  $N_E$  by Lemma \ref{Newse0}(ii) and Theorem \ref{trivial} for the 24 cases with $r_E(T)=0$ and by Lemma \ref{bb} for the 7 cases with  $r_E(T)=1.$  Use Proposition \ref{3.59} for some $p N_E$. 
\end{enumerate}

\vspace{ 5pt}

\item $\gS_\gl(A)=\{E_1,E_2\}$ with $\dim E_i=2$ and $F_i = \Q(E_i)$.    
\begin{enumerate}[A.]
\item Suppose at least one of the $F_i$ is ramified at 2.  When $N=N_{E_1} N_{E_2},$ we   use Lemma \ref{condL0} and Corollary \ref{2Es} to eliminate all but  three cases.    For those, we know examples labeled ``Prym" in Table 2.  One rules out $N=pN_{E_1}N_{E_2},$ with  $p$  inert in  $F_1$ and $F_2,$ by Remark \ref{Groth}, even when $E_1\simeq E_2.$  This happens for $N=3\!\cdot\!11^2,\, 5\!\cdot\!11^2$ and $3\!\cdot\!11\!\cdot\! 19.$ 

\item Suppose $F_1$ and $F_2$ are unramified at 2.  Locally at 2, the associated group schemes must be  \'etale or multiplicative and so are Cartier duals. The only case available is $N=713=23\cdot31$, for which we know two isogeny classes of Jacobians in Table \ref{EX}.  
\end{enumerate}

\vspace{5 pt}

\item  $\gS_\gl(A)=\{E\}$ with $\dim E=4.$  The criteria in \S \ref{E4} are given in terms of a {\em stem field}, i.e. a subfield $K$ of $\Q(E)$ whose Galois closure is $\Q(E).$
\begin{enumerate}[A.]
\item Thirteen quintic fields are candidates for $K.$  Each has Galois group $\cS_5$ and is determined by its conductor $N_E.$ The single case of the form $qN_E\le 1000$ is $3\!\cdot\!277$ and is ruled out by Remark \ref{Groth}.

\item  The only candidates for a sextic  $K$ are three  $\cS_3\wr\cS_2$-fields and four $\cS_6$-fields. 
\end{enumerate}
\end{enumerate}

\medskip

Table \ref{BadCond} lists all odd integers $N$ which were not eliminated by our criteria and are {\em not} conductors of known semistable surfaces.  

\begin{Not}\label{key}\hfill
\begin{enumerate}[{\rm i)}] \item {\bf u} means  $A[\gl]$ is prosaic; that is $\gS_\gl^{all}(A)=\{\F, \F, \F, \F\},$ with $\F = \F_2$.
\item  boldface integers are the Artin conductors of any two-dimensional constituents.
\item {\bf q} means $\gS_\gl^{all}(A)=\{ E\}$ for an irreducible, symplectic $E $ with  $\Gal(\Q(E)/\Q)$ isomorphic to $\cS_5$ and quintic stem field for $\Q(E).$ 

\item  {\bf wr72} or  ${\bf \cS_6}$ means $\gS_\gl^{all}(A)=\{ E\}$ for an irreducible  symplectic $E,$ with sextic stem field for $\Q(E)$ and $\Gal(\Q(E)/\Q)\simeq\cS_3\wr\cS_2$ or ${\bf \cS_6}$ respectively.
\end{enumerate}
\end{Not}

\begin{table}[h]\begin{caption}[BC]{Hypothetical Semistable Odd Conductors Not Eliminated}\label{BadCond}\end{caption}
\begin{tabular}{|c|c||c|c||c|c||c|c||c|c||}
\hline
N&WHY&N&WHY&N&WHY&N&WHY&N&WHY\\
\hline
415 & {\bf 83}  &  613 & {\bf q}  & 687& {\bf 229} &847& ${\bf 11,11}$  & 921& {\bf 307}  \\
417& {\bf 139} &  615 & {\bf u}  & 695 & {\bf 139} & 849 & {\bf 283} & 927 & {\bf wr72} \\
531& {\bf 59}   &  629 & {\bf 37}&697 & {\bf u}  &  853 & {\bf q}   & 957 & {\bf 11} \\
535& {\bf 107}&  637&{\bf 91}  &  735&{\bf u} & 859& {\bf 859} & 961 & {\bf 31,31}\\
547&{\bf q}  & 645&{\bf 43} & 747& {\bf 83} &873&{\bf u}  & 963 & {\bf 107} \\
 571&{\bf 571} &  649& {\bf 59} & 749& {\bf 107} &885&  {\bf 59} & 969&{\bf u} \\
581& {\bf 83}&657&{\bf u} &   767& {\bf 59}& 897&{\bf u}& 985&{\bf 197} \\
591&{\bf 197} & 663&{\bf u}& 775&{\bf u} & 903&{\bf 43}  & 989&{\bf 43} \\
595&{\bf u}& 669&{\bf 223} & 777&{\bf u},\,{\bf 37}&913& {\bf 83} &991& {\bf q}  \\
599 & {\bf q}  & 677& {\bf q}  &841& ${\bf 29,29}$ &917&{\bf 131} & 993&{\bf 331}  \\
\hline
\end{tabular}\smallskip\end{table}

Odd square  conductors $N < 841$ do not appear in Table \ref{BadCond} by Theorem \ref{sc29}.  We could not eliminate  semistable  $W$'s of Artin conductor $N=657,$ 775 and 847 because $W\simeq B[2]$ for a  {\em non}-semistable surface $B$   of that conductor  in Table 2.  

\noindent For most other conductors $N$   in Table 1, there are semistable abelian surfaces $B$ whose conductor is a proper multiple of $N$ with $W\simeq B[2]$.    

Only  903 and 969 in  Table \ref{BadCond} should  be conductors of surfaces under our conjectures and data in \cite{PoYu2}.   There should also exist 4-dimensional abelian varieties with $\go=\Z[\sqrt{2}]$ and reduced conductors 637, 645 and 927 and a 6-dimensional abelian variety of reduced conductor 991 with $\go$  the maximal order of the cubic field of discriminant 148. Taking $\gl$ as the prime of degree one over 2 in those cases   is  consistent with the corresponding entries in Table \ref{BadCond}.

\section{Abelian surfaces of odd conductor $< 1000$}\label{jacs} 
Table \ref{EX} gives one member of each isogeny class of paramodular  abelian {\em surfaces} of odd conductors below 1000 known  to us and found  by purely  {\em ad hoc} methods.  The  ``INFO" column uses Notation \ref{key}. In general, our methods rule  out all other possibilities for  $\gS_\gl^{\rm all}(A)$. Most examples are semistable,  except for those labeled ``notSS."  If a polynomial $f(x)$ is given, then the  surface is the Jacobian of the curve $y^2=f(x).$  Analogous  tables for even conductors will  appear in \cite{PoYu2}.

Let $C_{\!/\Q}$ be a curve and $\cC$ a global  integral model over $\Z.$  We have  {\em mild}  reduction at $p$ if $\cC$  is bad at $p$, but  the  N\'eron model of $J(C)$ is not.  Assume  that $C$ is given by the {\em non}-minimal model 
$$
C\!:\,  y^2 +  ( m a_1 s + m^2 a_3)y= ma_0 s^3+m^2 a_2 s^2+m^3 a_4 s+m^4 a_6,$$
where $m$ is an integer, $s$ is a quadratic polynomial in $\Z[x],$ \'etale  mod $m,$ $a_0$ is an integer prime to $m $ and all other $a_i$ are linear in $\Z[x].$  If the discriminant of $C$ is $m^{22}n$ with $n$ prime to $m,$ then the prime divisors of $m$ are   of mild reduction. The  converse can be  deduced by strong approximation from  \cite{Liu}. In  Table \ref{EX}, such a curve is indicated by the symbol ``mild@$m$" and the conductor of its Jacobian is in the  first column. 
  
If $X$ is a curve of genus three with a degree two cover of a genus one curve $C$, then the kernel Prym$(X/C)$ of the natural projection  $\pi\! : J(X)\to J(C)$ is an abelian surface with (1,2)-polarization.  Its conductor is  the quotient of that of $J(X)$ by that of $J(C).$  The surfaces of conductors 561,\, 665,\, 737 are such Pryms. They  are {\em not} $\Q$-isogenous to Jacobians and will be described in  a note \cite{BK4} on abelian surfaces of polarization (1,2).  

Let $E$ be an elliptic curve, defined over  $k=\Q(\sqrt{d}),$  of conductor $\mathfrak{c} $  and not isogenous  to its conjugate. Then the Weil restriction $S=R_{k/\Q}E$ is a surface of paramodular type with conductor $d^2N_{k/\Q}(\mathfrak{c})$ (see \cite{Mil2}). 
The surfaces of conductors $657=3^2\! \cdot\!73$ and $775=5^2\!\cdot\!31$  are  Weil restrictions of    curves  defined over $\Q(\sqrt{-3})$ and $\Q(\sqrt{5}),$ respectively.  It is expected that  elliptic curves over real fields should correspond to  parallel weight 2 Hilbert modular forms. In the recent preprint \cite{JLR}, such Hilbert modular eigenforms over real quadratic fields are lifted to paramodular forms when Remark \ref{rem}(iii) applies.  This supports our conjecture, with expected level, for Weil restrictions of ``Hilbert modular elliptic curves."  For imaginary quadratic fields, work of Cremona and his students combined with   \cite{Tay2} suggests modularity there as well. 

\vfill\pagebreak
 
\begin{center}
\begin{table}[h]\begin{caption}[EX]{Paramodular Abelian Surfaces of ODD Conductor $<1000.$}\label{EX}\end{caption}
\begin{tabular}{|c|c|c|c|}
\hline
 N &EQUATION& INFO\\
\hline
249&$x^6 + 4x^5 + 4x^4 + 2x^3 + 1$&{\bf 83}\\
277&$x^6 + 2x^5 + 3x^4 + 4x^3 - x^2 - 2x + 1$&{\bf q}\\
295&$x^6 - 2x^3 - 4x^2 + 1$&{\bf 59}\\
349&$x^6 - 2x^5 + 3x^4 - x^2 - 2x + 1$&{\bf q}\\
353&$ x^6 + 2x^5 + 5x^4 + 2x^3 + 2x^2 + 1$&{\bf wr72}\\
389&$x^6 + 2x^5 + 5x^4 + 8x^3 + 8x^2 + 4x$&{\bf 389}\\
427&$x^6 - 4x^5 - 4x^4 + 18x^3 + 16x^2 - 16x - 15$&{\bf 61}\\
461&$x^6 + 2x^5 - 5x^4 - 8x^3 + 11x^2 + 10x - 11$&{\bf q}\\
523&$x^6 - 2x^5 + x^4 + 4x^3 - 4x^2 - 4x$&{\bf 523}\\
555&$x^6 + 6x^5 + 5x^4 - 16x^3 - 8x^2 + 12x$&{\bf 37}\\
561&  PRYM&${\bf 11,51}$\\
587a&$-3x^6 + 18x^4 + 6x^3 + 9x^2 - 54x + 57$&${\bf \cS_6}$, mild@3\\
587b&$x^6 + 2x^4 + 2x^3 - 3x^2 - 2x + 1$&${\bf \cS_6}$\\
597&$x^6 + 4x^5 + 8x^4 + 12x^3 + 8x^2 + 4x$&{\bf q}\\
603&$ x^6 - 4x^5 + 2x^4 + 4x^3 + x^2 - 4x$&{\bf 67}\\
623&$ -224x^6 - 1504x^5 - 4448x^4 - 7200x^3 - 6080x^2 - 2048x$&{\bf 89}, mild@8\\
633&$24x^6 + 40x^5 + 28x^4 + 80x^3 + 52x^2 - 32x$&{\bf 211}, mild@2\\
657& WEIL RESTRICTION & {\bf u}, notSS\\
665 & PRYM & ${\bf 19,35}$\\
691&$x^6 + 2x^5 - 3x^4 - 4x^3 + 4x$&{\bf 691}\\
709&$-4x^5 - 7x^4 - 4x$&{\bf 709}\\
713a&$x^6 + 2x^5 + x^4 + 2x^3 - 2x^2 + 1$&${\bf 23,31}$\\
713b&$ x^6 - 2x^5 + x^4 + 2x^3 + 2x^2 - 4x + 1$&${\bf 23,31}$\\
731&$x^6 - 6x^4 + 4x^3 + 9x^2 - 16x - 4$&{\bf 43}\\
737&  PRYM & ${\bf 11,67}$\\
741&$ x^6 - 6x^5 + 9x^4 - 4x^2 + 12x$&{\bf 19}\\
743&$x^6 - 2x^4 - 2x^3 + 5x^2 - 2x + 1$&${\bf \cS_6}$\\
745&$x^6 + 2x^4 - 2x^3 + x^2 + 2x + 1$&{\bf wr72}\\
763&$4x^5 + 9x^4 - 6x^2 + 1$&{\bf 763}\\
775& WEIL RESTRICTION & {\bf u}, notSS\\
797&$x^6 + 4x^3 - 4x^2 + 4x$&{\bf q}\\
807&$ x^6 - 4x^5 + 2x^4 + 8x^3 - 3x^2 - 8x - 4$&{\bf 269}\\
847 &$x^6 - 2x^5 + 5x^4 - 4x^3 + 4x - 8$&${\bf 11,11}$, notSS\\
893a&$5x^6 - 40x^5 + 30x^4 - 510x^3 - 195x^2 - 1690x - 1295$&${\bf \cS_6}$, mild@5\\
893b&$x^6 - 2x^4 - 2x^3 - 3x^2 - 2x + 1$&${\bf \cS_6}$\\
901&$-7x^6 - 140x^5 - 532x^4 - 966x^3 - 504x^2 + 1596x + 2065$&{\bf 53}, mild@7\\
909&$x^6 - 2x^4 + 5x^2 - 4x$&{\bf 101}\\
925&$4x^5 + 8x^4 + 4x^3 - 3x^2 - 2x + 1$&{\bf 37}\\
953&$x^6 - 2x^5 + 5x^4 - 6x^3 + 2x^2 + 1$&{\bf wr72}\\
971&$ x^6 + 4x^5 - 8x^3 + 4x$&{\bf q}\\
975&$x^6 + 4x^5 - 6x^3 + 8x^2 - 4x - 3$ &{\bf u}, notSS\\
997a&$x^6 + 2x^5 + x^4 + 4x^2 + 4x$&{\bf 997}\\
997b&$x^6 - 4x^4 - 8x^3 - 8x^2 - 4x$&{\bf q}\\
\hline
\end{tabular}\smallskip\end{table}
\end{center} 

\vfill

\end{document}